\documentclass[a4paper,11pt]{article}
\usepackage[utf8x]{inputenc}
\usepackage{amssymb}
\usepackage{amsthm}
\usepackage{graphicx}

\newcommand{\Rm}{{\mathbb R}}

\newcommand{\eps}{\varepsilon}

\newtheorem{theorem}{Theorem}[section]
\newtheorem{lemma}[theorem]{Lemma}
\newtheorem{proposition}[theorem]{Proposition}

\title{On the Energy Subcritical, Non-linear Wave Equation with Radial Data for $p \in (3,5)$}
\author{Ruipeng Shen}

\begin{document}

\maketitle

\section{Introduction}
In this paper we will consider the energy subcritical, non-linear wave equation in $\Rm^3$ with radial initial data.
\begin{equation}
\left\{\begin{array}{l} \partial_t^2 u - \Delta u = \pm |u|^{p-1} u, \,\,\,\, (x,t)\in \Rm^3 \times \Rm;\\
u |_{t=0} = u_0 \in \dot{H}^{s_p} (\Rm^3);\\
\partial_t u |_{t=0} = u_1 \in \dot{H}^{s_p-1}(\Rm^3).\end{array}\right.
\label{eqn}
\end{equation}
Here $3 <p <5$ and
\[
 s_p = \frac{3}{2} - \frac{2}{p-1}.
\]
The positive sign in the non-linear term gives us the focusing case, while the negative sign indicates the defocusing case.
The following quantity is called the energy of the solution. The energy is a constant in the whole lifespan of the solution,
as long as it is well-defined.
\begin{equation}
 E(t) = \frac{1}{2} \int_{\Rm^3} \left(|\partial_t u (x,t)|^2 + |\nabla u (x,t)|^2\right) dx \mp
\frac{1}{p+1} \int_{\Rm^3} |u(x,t)|^{p+1} dx. \label{defenergy}
\end{equation}
Please note that the energy could be a negative number in the focusing case.
\paragraph{Previous Results in the Energy-critical Case} In the energy-critical case, namely $p=5$, the initial data is in the energy space
$\dot{H}^1 \times L^2$. This automatically guarantees the existence of the energy by the Sobolev embedding. This kind of wave equations have been
extensively studied.
In the defocusing case, M. Grillakis (See \cite{mg1, mg2}) proved the global existence and scattering of the solution with any $\dot{H}^1 \times L^2$
initial data in 1990's. In the focusing case, however, the behavior of solutions is much more complicated. The solutions may scatter,
blow up in finite time or even be independent of time. Please see \cite{secret, kenig} for more details. In particular,
a solution independent of time is usually called a ground state or a soliton.
This kind of solutions are actually the solutions of the elliptic equation $-\Delta W(x) = |W(x)|^{p-1} W(x)$.
We can write down all the nontrivial radial solitons explicitly as below. The letter $\lambda$ here is an arbitrary positive parameter.
\begin{equation}
 W(x) = \pm \frac{1}{\lambda^{1/2}} \left(1 + \frac{|x|^2}{3 \lambda^2} \right)^{-1/2}. \label{radial solitons}
\end{equation}
\paragraph{Energy Subcritical Case} We will consider the case $3< p <5$ in this paper, thus $1/2 <s_p <1$.
In this case the problem is critical in the space $\dot{H}^{s_p}(\Rm^3) \times \dot{H}^{s_p -1}(\Rm^3)$,
because if $u(x,t)$ is a solution of (\ref{eqn}) with initial data $(u_0,u_1)$, then for any $\lambda>0$, the function
\[
 \frac{1}{\lambda^{3/2 - s_p}} u \left(\frac{x}{\lambda}, \frac{t}{\lambda}\right)
\]
is another solution of the equation (\ref{eqn}) with the initial data
\[
 \left(\frac{1}{\lambda^{3/2 - s_p}} u_0 (\frac{x}{\lambda}),
\frac{1}{\lambda^{5/2 - s_p}} u_1 (\frac{x}{\lambda})\right),
\]
which shares the same $\dot{H}^{s_p} \times \dot{H}^{s_p -1}$ norm as the original initial data $(u_0,u_1)$.
These scalings play an important role in our discussion of this problem.
\begin{theorem} \emph{\textbf{(Main Theorem)}} \label{main theorem}
Let $u$ be a solution of the non-linear wave equation (\ref{eqn}) with radial initial data
$(u_0,u_1) \in \dot{H}^{s_p} \times \dot{H}^{s_p -1} (\Rm^3)$ and a maximal lifespan $I$ so that
\begin{equation}
 \sup_{t \in I} \|(u, \partial_t u)\|_{{\dot{H}^{s_p}}\times{\dot{H}^{s_p -1}}} < \infty. \label{add1}
\end{equation}
Then $u$ is global in time $(I = \Rm)$ and scatters, i.e.
\[
 \|u(x,t)\|_{S(\Rm)} < \infty,\; \hbox{or equivalently}\; \|u(x,t)\|_{Y_{s_p}(\Rm)} < \infty.
\]
This is actually equivalent to saying that there exist two pairs $(u_0^+,u_1^+)$ and $(u_0^-,u_1^-)$
in the space $\dot{H}^{s_p} \times \dot{H}^{s_p-1}$ such that
\[
 \lim_{t \rightarrow \pm \infty} \left\|\left(u(t) - S(t)(u_0^\pm, u_1^\pm),
\partial_t u(t) - \partial_t S(t)(u_0^\pm, u_1^\pm)\right)\right\|_{\dot{H}^{s_p} \times \dot{H}^{s_p-1}} = 0.
\]
Here $S(t)(u_0^\pm,u_1^\pm)$ is the solution of the Linear Wave Equation with the initial data $(u_0^\pm, u_1^\pm)$.
\end{theorem}
\noindent Please refer to section 2 for the definition of the $S$ and $Y_{s}$ norms.
\paragraph{Remark on the Defocusing Case} As in the energy-critical case, we expect that the solutions always scatter. Besides the radial condition,
the main theorem depends on the assumption (\ref{add1}), which is expected to be true for all solutions. Unfortunately, as far as the author knows,
no one actually knows how to prove it without additional assumptions.
\paragraph{Remark on the Focusing Case} In the focusing case, the solutions may blow up in finite time. (Please see theorem
\ref{Non-positive energy implies blow-up}, for instance) Thus the assumption (\ref{add1}) is a meaningful
and essential condition rather than a technical one. The main theorem gives us the following rough classification of the radial solutions.
\begin{proposition}
Let $u(t)$ be a solution of (\ref{eqn}) in the focusing case with a maximal lifespan $I$ and radial
initial data $(u_0,u_1) \in \dot{H}^{s_p} \times \dot{H}^{s_p-1}(\Rm^3)$. Then one of the following holds for $u(x,t)$.
\begin{itemize}
\item (I) (Blow-up) The $\dot{H}^{s_p} \times \dot{H}^{s_p -1}$ norm of $(u(t),\partial_t u(t))$ blows up, namely
\[
 \sup_{t \in I} \|(u(t), \partial_t u(t))\|_{{\dot{H}^{s_p}}\times{\dot{H}^{s_p -1}}} = + \infty.
\]
\item (II) (Scattering) If the upper bound of the $\dot{H}^{s_p} \times \dot{H}^{s_p -1}$ norm above is finite instead, namely,
the assumption (\ref{add1}) holds, then $u(t)$ is a global solution (i.e $I = \Rm$) and scatters.
\end{itemize}
\end{proposition}
\paragraph{Main Idea in this Paper} The main idea to establish theorem \ref{main theorem} is to use the compactness/rigidity argument, namely to show
\begin{itemize}
\item (I) If the main theorem failed, it would break down at a minimal blow-up solution, which is almost periodic modulo scalings.
\item (II) The minimal blow-up solution is in the energy space.
\item (III) The minimal blow-up solution described above does not exist.
\end{itemize}
\paragraph{Step (I)} The method of profile decomposition used here has been a standard way to deal with both the wave equation
and the Schr\"{o}dinger equation. Thus we will only give important statements instead of showing all the details.
The other steps, however, depend on the specific problems. One could refer to \cite{bahouri} in order to understand what is the profile
decomposition, and to \cite{kenig, kenig2} in order to see why the profile decomposition leads to the existence of a minimal blow-up solution.
\paragraph{Step (II)} We will combine the method used in my old paper \cite{shen1} and a method used in C.E.Kenig and F.Merle's paper
\cite{km} on the supercritical case of the non-linear wave equation in $\Rm^3$. The idea is to use the following fact.
Given a radial solution $u(x,t)$ of the equation
\[
 \partial_t^2 u(x,t) - \Delta u(x,t) = F(x,t)
\]
in the time interval $I$,
if we define two functions $w, h: \Rm^+ \times I \rightarrow \Rm$, such that $w(|x|,t) = |x|u(x,t)$ and $h(|x|,t)=|x| F(x,t)$, then
$w(r,t)$ is a solution of the one-dimensional wave equation $\partial_t^2 w(r,t) - \partial_r^2 w(r,t) = h(r,t)$. This makes it convenient
to consider the integral
\[
 \int_{r_0 \pm t}^{4r_0 \pm t} |\partial_t w (r, t_0 +t) \mp \partial_r w(r, t_0 +t)|^2 dr.
\]
as the parameter $t$ moves.
\paragraph{Step (III)} Given an energy estimate, all minimal blow-up solutions are not difficult to kill except for the soliton-like solutions
in the focusing case. As I mentioned earlier, this kind of solutions actually exist in the energy-critical case.
The ground states given in (\ref{radial solitons}) are perfect examples.
In the energy subcritical case, however, the soliton does not exist at all. More precisely,
none of the solutions of the corresponding elliptic equation is in the right space $\dot{H}^{s_p}$.
This fact enables us to gain a contradiction by showing a soliton-like minimal blow-up solution must be a real soliton,
which does not exist, using a new method introduced by Thomas Duyckaerts, Carlos Kenig and Frank Merle. They classified all radial
solutions of the energy-critical, focusing wave equation in their recent paper \cite{secret} using this "channel of energy" method.
\paragraph{Remark on the Supercritical Case} Simultaneously to this work, Thomas Duyckaerts, Carlos Kenig and Frank Merle \cite{dkm2}
proved that similar results to those in this paper also hold in the supercritical case $p>5$ of the focusing wave equation,
using the compactness/rigidity argument, a point-wise estimate on "compact" solutions obtained in
the paper \cite{km} and the "channel of energy" method mentioned above.
\section{Preliminary results}
\subsection{Local Theory with $\dot{H}^{s_p} \times \dot{H}^{s_p-1}(\Rm^3)$ Initial Data}
In this section, we will review the theory for the Cauchy problem of the nonlinear wave equation (\ref{eqn}) with initial data in the critical space
$\dot{H}^{s_p} \times \dot{H}^{s_p-1}(\Rm^3)$. The same local theory works in both the focusing and defocusing cases.
It could be also applied to the non-radial case.
\paragraph{Space-time Norm} Let $I$ be an interval of time. The space-time norm is defined by
\[
 \|v(x,t)\|_{L^{q} L^{r} (I \times \Rm^3)} = \left(\int_I \left(\int_{\Rm^3} |v(x,t)|^r dx\right)^{q/r} dt\right)^{1/q}.
\]
This is used in the following Strichartz estimates.
\begin{proposition} \emph{\textbf{Generalized Strichartz Inequalities}} \label{strichartz}
(Please see proposition 3.1 of \cite{strichartz}, here we use the Sobolev version in $\Rm^3$)
Let $2 \leq q_1,q_2 \leq \infty$, $2 \leq r_1, r_2 < \infty$ and $\rho_1, \rho_2, s \in \Rm$ with
\[
 1/{q_i} + 1/{r_i} \leq 1/2; \,\,\,\, i=1,2.
\]
\[
 1/{q_1} + 3/{r_1} = 3/2 - s + \rho_1.
\]
\[
 1/{q_2} + 3/{r_2} = 1/2 + s + \rho_2.
\]
In particular, if $(q_1, r_1, s, \rho_1) = (q, r, m, 0)$ satisfies the conditions above, we say
$(q,r)$ is an $m$-admissible pair.\\
Let $u$ be the solution of the following linear wave equation
\begin{equation}
\left\{\begin{array}{l} \partial_t^2 u - \Delta u = F (x,t), \,\,\,\,\, (x,t)\in \Rm^3 \times \Rm;\\
u |_{t=0} = u_0 \in \dot{H}^s (\Rm^3);\\
\partial_t u |_{t=0} = u_1 \in \dot{H}^{s-1}(\Rm^3).\end{array}\right.
\end{equation}
Then we have
\begin{eqnarray*}
 \lefteqn{\|(u(T), \partial_t u(T))\|_{\displaystyle \dot{H}^{s} \times \dot{H}^{s-1}} +
 \|D_x^{\rho_1} u\|_{\displaystyle L^{q_1} L^{r_1} ([0,T]\times \Rm^3)}}\\
 &\leq& C \left( \|(u_0,u_1)\|_{\displaystyle \dot{H}^{s} \times \dot{H}^{s-1}} +
 \|D_x^{-\rho_2} F(x,t)\|_{\displaystyle L^{\bar{q}_2} L^{\bar{r}_2}([0,T]\times \Rm^3)} \right).
\end{eqnarray*}
The constant $C$ does not depend on $T$.
\end{proposition}
\paragraph{Definition of Norms} Fix $3<p<5$. We define the following norms with $s_p \leq s < 1$
\begin{eqnarray*}
 \|v(x,t)\|_{S(I)} &=& \|v(x,t)\|_{L^{2(p-1)} L^{2(p-1)} (I \times \Rm^3)};\\
 \|v(x,t)\|_{W(I)} &=& \|v(x,t)\|_{{L^4} L^{4} (I \times \Rm^3)};\\
 \|v(x,t)\|_{Z_s (I)} &=& \|v(x,t)\|_{\displaystyle L^{\frac{2}{s+1}} L^{\frac{2}{2-s}}(I \times \Rm^3)};\\
 \|v(x,t)\|_{Y_s (I)} &=& \|v(x,t)\|_{\displaystyle L^{\frac{2p}{s + 1 -(2p-2) (s -s_p)}}L^{\frac{2p}{2-s}}(I \times \Rm^3)}.
\end{eqnarray*}
\paragraph{Remark} By the Strichartz estimates, we have if $u(x,t)$ is the solution of
\[
\left\{\begin{array}{l} \partial_t^2 u - \Delta u = F(x,t), \,\,\,\, (x,t)\in \Rm^3 \times \Rm;\\
u |_{t=0} = u_0 \in \dot{H}^{s} (\Rm^3);\\
\partial_t u |_{t=0} = u_1 \in \dot{H}^{s-1}(\Rm^3).\end{array}\right.
\]
then
\begin{eqnarray*}
 \lefteqn{\|(u(T), \partial_t u(T))\|_{\displaystyle \dot{H}^{s} \times \dot{H}^{s-1}} + \|u\|_{Y_s([0,T])}}\\
 &\leq&
 C \left( \|(u_0,u_1)\|_{\displaystyle \dot{H}^{s} \times \dot{H}^{s-1}} + \|F(x,t)\|_{Z_s ([0,T])} \right).
\end{eqnarray*}
\paragraph{Definition of Solutions} We say $u(t) (t\in I)$ is a solution of (\ref{eqn}),
if $(u,\partial_t u) \in C(I;{\dot{H}^{s_p}}\times{\dot{H}^{s_p-1}})$, with finite norms $\|u\|_{S(J)}$ and $\|D_x^{s_p - {1/2}} u\|_{W(J)}$
for any bounded closed interval $J \subseteq I$ so that the integral equation
\[
 u(t) = S(t)(u_0,u_1) + \int_0^t \frac{\sin ((t-\tau)\sqrt{-\Delta})}{\sqrt{-\Delta}} F(u(\tau)) d\tau
\]
holds for all time $t \in I$. Here $S(t)(u_0,u_1)$ is the solution of the linear wave equation with initial data $(u_0,u_1)$ and
\[
 F(u) = \pm |u|^{p-1}u.
\]
\paragraph{Remark} We can take another way to define the solutions by substituting $S(I)$ and $W(I)$ norms
by a single $Y_{s_p}(I)$ norm. Using the Strichartz estimates, these two definitions are equivalent to each other.
\paragraph{Local Theory} By the Strichartz estimate and a fixed-point argument, we have the following theorems.
(Please see \cite{ls} for more details)
%\paragraph{Strichartz estimate} Let $\rho_1, \rho_2, \mu \in \Rm$ and $2 \leq q_1$
\begin{theorem} \emph{\textbf{(Local solution)}}
For any initial data $(u_0,u_1) \in \dot{H}^{s_p} \times \dot{H}^{s_p-1}$, there is a maximal interval $(-T_{-}(u_0,u_1), T_{+}(u_0,u_1))$
in which the equation has a solution.
\end{theorem}
\begin{theorem}
\emph{\textbf{(Scattering with small data)}} There exists $\delta = \delta(p) > 0$ such that if the norm of the initial data
$\|(u_0,u_1)\|_{\dot{H}^{s_p} \times \dot{H}^{s_p-1}} < \delta$, then the Cauchy problem (\ref{eqn}) has a global-in-time solution $u$ with
$\|u\|_{S(-\infty,+\infty)} < \infty$.
\end{theorem}
\begin{lemma} \emph{\textbf{(Standard finite blow-up criterion)}}
If $T_{+} < \infty$, %under the uniform boundedness condition (\ref{add1})
then
\[
 \|u\|_{S([0,T_{+}))} = \infty.
\]
\end{lemma}
\begin{theorem} \emph{\textbf{(Long time perturbation theory)}}(See \cite{ctao, kenig, kenig1, km})
Let $M, A, A'$ be positive constants. There exists $\eps_0 = \eps_0(M,A,A')>0$ and $\beta>0$ such that if $\eps < \eps_0$,
then for any approximation solution $\tilde{u}$ defined on $\Rm^3 \times I$ ($0\in I$)
and any initial data $(u_0,u_1) \in \dot{H}^{s_p} \times \dot{H}^{s_p-1}$ satisfying
\[
 (\partial_t^2 - \Delta) (\tilde{u}) - F(\tilde{u}) = e(x,t), \,\,\,\,\, (x,t) \in \Rm^3 \times I;
\]
\begin{equation}
\left\{ \begin{array}{l}
\sup_{t \in I} \|(\tilde{u}(t), \partial_t \tilde{u}(t))\|_{\dot{H}^{s_p} \times \dot{H}^{s_p-1}} \leq A,\\
 \|\tilde{u}\|_{S(I)} \leq M,\\
\|D_x^{s_p -1/2}\tilde{u}\|_{W(J)} < \infty \hbox{ for each } J\subset\subset I; \end{array} \right.
\end{equation}
\[
 \|(u_0-\tilde{u}(0), u_1-\partial_t \tilde{u}(0))\|_{\dot{H}^{s_p} \times \dot{H}^{s_p-1}} \leq A';
\]
\[
 \|D_x^{s_p - \frac{1}{2}} e\|_{L_I^{4/3} L_x^{4/3}} + \|S(t)(u_0-\tilde{u}(0),u_1 - \partial_t \tilde{u}(0))\|_{S(I)} \leq \eps;
\]
there exists a solution of (\ref{eqn}) defined in the interval $I$ with the initial data $(u_0,u_1)$ and satisfying
\[
  \|u\|_{S(I)} \leq C(M,A,A');
\]
\[
 \sup_{t\in I} \|((u(t), \partial_t u(t)) - ((\tilde{u}(t), \partial_t \tilde{u}(t))\|_{\dot{H}^{s_p} \times \dot{H}^{s_p-1}}
 \leq C(M,A,A')(A' +\eps +\eps^\beta).
\]
\end{theorem}
%\paragraph{Remark 1} We could use the $Y_{s_p}(I)$ norm instead of the $S(I)$ norm in this theorem.
\begin{theorem} \emph{\textbf{(Perturbation theory with $Y_{s_p}$ norm)}} \label{perturbation theory in Ysp}
Let $M$ be a positive constant. There exists a constant $\eps_0 = \eps_0 (M)>0$, such that if $\eps < \eps_0$,
then for any approximation solution $\tilde{u}$ defined on $\Rm^3 \times I$ ($0\in I$)
and any initial data $(u_0,u_1) \in \dot{H}^{s_p} \times \dot{H}^{s_p-1}$ satisfying
\[
 (\partial_t^2 - \Delta) (\tilde{u}) - F(\tilde{u}) = e(x,t), \,\,\,\,\, (x,t) \in \Rm^3 \times I;
\]
\[
 \|\tilde{u}\|_{Y_{s_p}(I)} < M;\;\;\; \|(\tilde{u}(0),\partial_t \tilde{u}(0))\|_{\dot{H}^{s_p}\times \dot{H}^{s_p-1}}< \infty;
\]
\[
 \|e(x,t)\|_{Z_{s_p}(I)}+ \|S(t)(u_0-\tilde{u}(0),u_1 - \partial_t \tilde{u}(0))\|_{Y_{s_p}(I)} \leq \eps;
\]
there exists a solution $u(x,t)$ of (\ref{eqn}) defined in the interval $I$ with the initial data $(u_0,u_1)$ and satisfying
\[
 \|u(x,t) - \tilde{u}(x,t)\|_{Y_{s_p}(I)} < C(M) \eps.
\]
\[
 \sup_{t \in I} \left\|\left(\begin{array}{c} u(t)\\ \partial_t u(t)\end{array}\right)
 - \left(\begin{array}{c} \tilde{u}(t)\\ \partial_t \tilde{u}(t)\end{array}\right)
 - S(t)\left(\begin{array}{c} u_0 - \tilde{u}(0)\\ u_1 -\partial_t \tilde{u}(0)\end{array}\right)
 \right\|_{\dot{H}^{s_p} \times \dot{H}^{s_p-1}} < C(M)\eps.
\]
\end{theorem}
\paragraph{Remark} If $K$ is a compact subset of the space $\dot{H}^{s_p} \times \dot{H}^{s_p-1}$, then there exists $T = T(K) > 0$
such that for any $(u_0,u_1) \in K$, $T_{+}(u_0,u_1) > T(K)$. This is a direct result from the perturbation theory.
\subsection{Local Theory with more regular initial data}
Let $s \in (s_p,1]$. By a similar fixed argument we can obtain the following results.
\begin{theorem} \emph{\textbf{(Local solution with $\dot{H}^s \times \dot{H}^{s-1}$ initial data)}} \label{local existence in Hs}
If $(u_0,u_1) \in \dot{H}^s \times \dot{H}^{s-1}$, then there is a maximal interval $(-T_{-}(u_0,u_1), T_{+}(u_0,u_1))$
in which the equation has a solution $u(x,t)$. In addition, we have
\[
 T_-(u_0,u_1), T_+ (u_0,u_1) > T_1 \doteq C_{s,p} (\|(u_0,u_1)\|_{\dot{H}^s \times \dot{H}^{s-1}})^{-1 /(s-s_p)};
\]
\[
 \|u(x,t)\|_{Y_s ([-T_1,T_1])} \leq C_{s,p} \|(u_0,u_1)\|_{\dot{H}^s \times \dot{H}^{s-1}}.
\]
\end{theorem}
\begin{theorem} \emph{\textbf{(Weak long-time perturbation theory)}} \label{perturbation theory in Hs}
Let $\tilde{u}$ be a solution of the equation (\ref{eqn}) in the time interval
$[0,T]$ with initial data $(\tilde{u}_0,\tilde{u}_1)$, so that
\[
 \|(\tilde{u}_0,\tilde{u}_1)\|_{\dot{H}^s \times \dot{H}^{s-1}} < \infty;\; \|\tilde{u}\|_{Y_s ([0,T])} < M.
\]
There exist two constants $\eps_0 (T,M), C(T,M)>0$, such that if $(u_0,u_1)$ is another pair of initial data with
\[
 \|(u_0 - \tilde{u}_0,u_1- \tilde{u}_1)\|_{\dot{H}^s \times \dot{H}^{s-1}} < \eps_0 (T,M),
\]
then there exists a solution $u$ of the equation (\ref{eqn}) in the time interval $[0,T]$ with initial data $(u_0,u_1)$ so that
\[
 \|u - \tilde{u}\|_{Y_s ([0,T])} \leq C(T,M) \|(u_0 - \tilde{u}_0,u_1- \tilde{u}_1)\|_{\dot{H}^s \times \dot{H}^{s-1}};
\]
\[
 \sup_{t \in [0,T]} \|(u(t) - \tilde{u}(t), \partial_t u(t) - \partial_t \tilde{u}(t))\|_{\dot{H}^s \times \dot{H}^{s-1}} \leq
 C(T,M) \|(u_0 - \tilde{u}_0,u_1- \tilde{u}_1)\|_{\dot{H}^s \times \dot{H}^{s-1}}.
\]
\end{theorem}
\subsection{Notations and Technical Results}
\begin{figure}[h]
\centering
\includegraphics{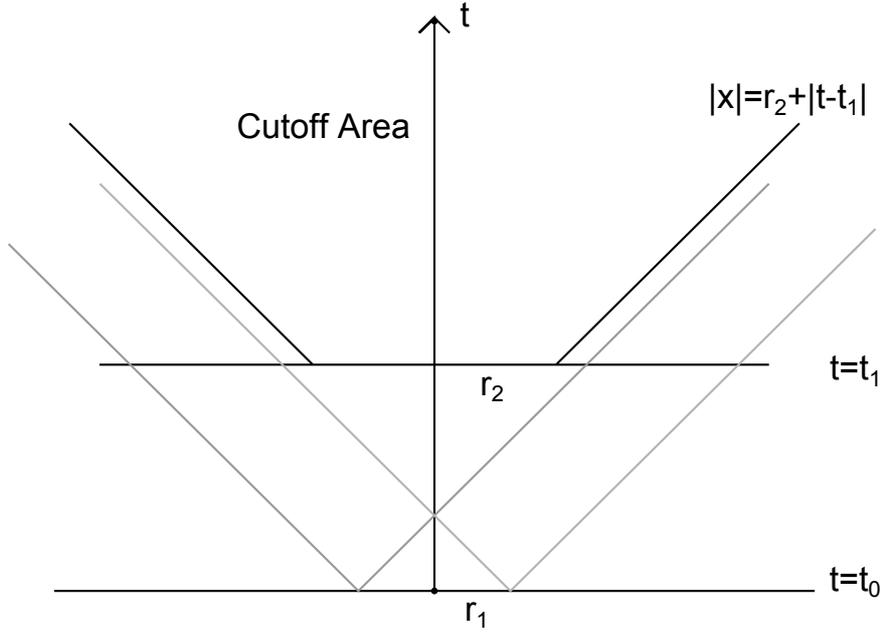}
\caption{Illustration of Proof}\label{drawing1}
\end{figure}

\paragraph{The Symbol $\lesssim$} Throughout this paper, the inequality $A \lesssim B$ means that there exists a constant $c$,
such that $A \leq c B$. In particular, a subscript of the symbol $\lesssim$ implies that the constant $c$ depends on the parameter(s) mentioned
in the subscript but nothing else.
\paragraph{The Smooth Frequency Cutoff} In this paper we use the notations $P_{<A}$ and $P_{>A}$ for the standard smooth frequency cutoff operators.
In particular, we use the following notation on $u$ for convenience.
\[
 u_{<A} \doteq P_{<A}u;\;\;\; u_{>A} \doteq P_{>A}u.
\]
\paragraph{Notation for Radial Functions} If $u(x,t)$ is radial in the space, then $u(r,t)$ represents the value $u(x,t)$ when $|x|=r$.
\paragraph{Linear Wave Evolution} Let $(u_0,u_1) \in \dot{H}^s \times \dot{H}^{s-1}(\Rm^3)$ be a pair of initial data. Suppose $u(x,t)$ is the solution
of the following linear wave equation
\[
\left\{\begin{array}{l} \partial_t^2 u - \Delta u = 0, \,\,\, (x,t)\in \Rm^3 \times \Rm;\\
u|_{t=0} = u_0;\\
\partial_t u |_{t=0} = u_1.\end{array}\right.
\]
We will use the following notations to represent this solution $u$.
\[
 S(t_0)(u_0,u_1) = u(t_0);
\]
\[
 S(t_0)\left(\begin{array}{l} u_0\\ u_1\end{array}\right) = \left(\begin{array}{c} u(t_0)\\ \partial_t u(t_0)\end{array}\right).
\]
\paragraph{Method of Center Cutoff} Let $(v_0, v_1) \in \dot{H}^1 \times L^2(\Rm^3 \setminus B(0,r))$ be a pair of radial functions. We define ($R>r$)
\[
 (\Psi_R v_0)(x) = \left\{\begin{array}{l} v_0 (x),\,\, \hbox{if}\,\, |x| > R;\\
  v_0 (R)\,\, \hbox{if}\,\, |x|\leq R. \end{array}\right.
\]
\[
 (\Psi_R v_1)(x) = \left\{\begin{array}{l} v_1 (x),\,\, \hbox{if}\,\, |x| > R;\\
  0,\,\, \hbox{if}\,\, |x|\leq R. \end{array} \right.
\]
\begin{lemma} \emph{\textbf{Glue of $\dot{H}^s$ Functions}} \label{glue}
Let $-1 \leq s \leq 1$. Suppose $f(x)$ is a tempered distribution defined on $\Rm^3$ such that ($R>0$)
\[
 f(x) = \left\{\begin{array}{l} f_1 (x),\; x \in B(0,2R)\\
 f_2 (x),\; x \in \Rm^3 \setminus B(0,R) \end{array}\right.
\]
with $f_1, f_2 \in \dot{H}^s (\Rm^3)$. Then $f$ is in the space $\dot{H}^s(\Rm^3)$ and
\[
 \|f\|_{\dot{H}^s (\Rm^3)} \leq C(s)\left( \|f_1\|_{\dot{H}^s (\Rm^3)} + \|f_2\|_{\dot{H}^s (\Rm^3)}\right).
\]
\end{lemma}
\paragraph{Proof} By a dilation we can always assume $R = 1$. Let $\phi (x)$ be a smooth, radial, nonnegative function such that
\[
 \phi(x) = \left\{\begin{array}{l} 1,\; x \in B(0,1)\\
 0,\; x \in \Rm^3 \setminus B(0,2) \end{array} \right.
\]
Let us define a linear operator: $P (f) = \phi(x) f$. We know this operator is bounded from $\dot{H}^1(\Rm^3)$ to $\dot{H}^1(\Rm^3)$, and
from $L^2 (\Rm^3)$ to $L^2 (\Rm^3)$. Thus by an interpolation, this is a bounded operator from $\dot{H}^s$ to itself if $0 < s < 1$.
By duality $P$ is also bounded from $\dot{H}^s$ to itself if $-1 \leq s \leq 0$.
In summary, $P$ is a bounded operator from $\dot{H}^s$ to itself for each $-1\leq s \leq 1$. Now we have
\[
 f = P f_1 + f_2 - P f_2
\]
as a tempered distribution. Thus
\[
 \|f\|_{\dot{H}^s} \leq \|P f_1\|_{\dot{H}^s} + \|f_2\|_{\dot{H}^s} + \|P f_2\|_{\dot{H}^s} \leq (\|P\|_s + 1) (\|f_1\|_{\dot{H}^s} + \|f_2\|_{\dot{H}^s}).
\]
\begin{lemma} Let $u(x,t)$ be a solution of the non-linear wave equation (\ref{eqn}) with the condition (\ref{add1}), then
\begin{equation}
\left\| \left(\begin{array}{l}  \displaystyle\int_{t_1}^{t_2} \frac{\sin((\tau -t)\sqrt{-\Delta})}{\sqrt{-\Delta}} F(u(\tau)) d\tau\\
 \displaystyle -\int_{t_1}^{t_2} \cos((\tau -t)\sqrt{-\Delta}) F(u(\tau)) d\tau\end{array} \right)\right\|_{\dot{H}^{s_p} \times \dot{H}^{s_p-1}} \lesssim 1.
 \label{uni hs es}
\end{equation}
\end{lemma}
\paragraph{Proof} Directly from the following identity.
\begin{equation}
 \left(\begin{array}{l} \displaystyle \int_{t_1}^{t_2} \frac{\sin((\tau -t)\sqrt{-\Delta})}{\sqrt{-\Delta}} F(u(\tau)) d\tau\\
 \displaystyle -\int_{t_1}^{t_2} \cos((\tau -t)\sqrt{-\Delta}) F(u(\tau)) d\tau\end{array} \right) =
 S(t-t_1)\left( \begin{array}{c} u(t_1)\\ \partial_t u (t_1) \end{array}\right) -
 S(t-t_2)\left( \begin{array}{c} u(t_2)\\ \partial_t u (t_2) \end{array}\right).
\end{equation}
\begin{lemma} (Please see lemma 3.2 of \cite{km}) Let $1/2 < s < 3/2$. If $u(y)$ is a radial $\dot{H}^s(\Rm^3)$ function, then
\begin{equation}
 |u(y)| \lesssim_s \frac{1}{|y|^{\frac{3}{2} -s}} \|u\|_{\dot{H}^{s}}. \label{estimate of radial hs function}
\end{equation}
\end{lemma}
\paragraph{Remark} This actually means that a radial $\dot{H}^s$ function is uniformly continuous in $\Rm^3 \setminus B(0,R)$ if $R > 0$.
\begin{lemma}\label{Strong Huygens in Duhamel}
Let $r_1, r_2 > 0$ and $t_0, t_1 \in \Rm$ so that $r_1 + r_2 \leq t_1 -t_0$.
Suppose $(u_0,u_1)$ is a weak limit in the space $\dot{H}^{s_p} \times \dot{H}^{s_p-1}$ as below
\begin{eqnarray*}
 u_0 &=& \lim_{T\rightarrow +\infty} \int_{t_1}^{T} \frac{\sin((t -t_0)\sqrt{-\Delta})}{\sqrt{-\Delta}} F(t) dt;\\
 u_1 &=& -\lim_{T\rightarrow +\infty} \int_{t_1}^{T} \cos((t -t_0)\sqrt{-\Delta}) F(t) dt.
\end{eqnarray*}
Here $F(x,t)$ is a function defined in $[t_1,\infty)\times \Rm^3$ with a finite $Z_{s_p}([t_1,T])$ norm for each $T>t_1$.
In addition, we have ($1/2 < s_1 \leq 1$, $\chi$ is a characteristic function of the region indicated)
\begin{equation}
 S= \|\chi_{|x|>r_2 + |t-t_1|}(x,t)F(x,t)\|_{L^1 L^\frac{6}{5-2 s_1}([t_1,\infty)\times \Rm^3)} < +\infty. \label{lin23}
\end{equation}
Then there exists a pair $(\tilde{u}_0,\tilde{u}_1)$ with
$\|(\tilde{u}_0, \tilde{u}_1)\|_{\dot{H}^{s_1} \times \dot{H}^{s_1 -1}(\Rm^3)} \leq C_{s_1} S$ and
\[
 (u_0, u_1) = (\tilde{u}_0, \tilde{u}_1)\; \hbox{in the ball}\; B(0,r_1).
\]
\end{lemma}
\paragraph{Proof} Let us define
\begin{eqnarray*}
 u_{0,T} &=& \int_{t_1}^{T} \frac{\sin((t -t_0)\sqrt{-\Delta})}{\sqrt{-\Delta}} F(t) dt;\\
 u_{1,T} &=& -\int_{t_1}^{T} \cos((t -t_0)\sqrt{-\Delta}) F(t) dt.\\
 \tilde{u}_{0,T} &=& \int_{t_1}^{T} \frac{\sin((t -t_0)\sqrt{-\Delta})}{\sqrt{-\Delta}} (\chi F(t)) dt;\\
 \tilde{u}_{1,T} &=& -\int_{t_1}^{T} \cos((t -t_0)\sqrt{-\Delta}) (\chi F(t)) dt.
\end{eqnarray*}
By the Strichartz estimates and the assumption (\ref{lin23}), we know the pair $(\tilde{u}_{0,T},\tilde{u}_{1,T})$ converges strongly in
$\dot{H}^{s_1}\times \dot{H}^{s_1-1}$ to a pair $(\tilde{u}_0, \tilde{u}_1)$ as $T \rightarrow +\infty$ so that
\[
 \|(\tilde{u}_0, \tilde{u}_1)\|_{\dot{H}^{s_1} \times \dot{H}^{s_1 -1}(\Rm^3)} \leq C_{s_1} S.
\]
In addition, we know the pair $(\tilde{u}_{0,T},\tilde{u}_{1,T})$ is the same as $(u_{0,T},u_{1,T})$ in the ball $B(0,r_1)$ by strong Huygens'
principal. The figure \ref{drawing1} shows the region where the value of $F(x,t)$ may affect the value of the integrals in the ball $B(0,r_1)$.
This region is disjoint with the cutoff area if $r_1 + r_2 \leq t_1 - t_0$.
As a result, the pair $(\tilde{u}_{0,T},\tilde{u}_{1,T})$ converges to $(u_0,u_1)$ weakly in the ball $B(0,r_1)$ as the pair $(u_{0,T},u_{1,T})$ does.
Considering both the strong and weak convergence, we conclude
\[
 (u_0, u_1) = (\tilde{u}_0, \tilde{u}_1)\; \hbox{in the ball}\; B(0,r_1).
\]
\section{Compactness Process}
As we stated in the first section, the standard technique here is to show if the main theorem failed, there would be a special minimal blow-up solution.
In addition, this solution is almost periodic modulo symmetries.
\paragraph{Definition} A solution $u(x,t)$ of (\ref{eqn})
is almost periodic modulo symmetries if there exists a positive function
$\lambda(t)$ defined on its maximal lifespan $I$ such that the set
\[
 \left\{\left(\frac{1}{\lambda(t)^{3/2-s_p}} u\left(\frac{x}{\lambda(t)}, t\right),
 \frac{1}{\lambda(t)^{5/2-s_p}} \partial_t u\left(\frac{x}{\lambda(t)}, t\right)\right): t \in I\right\}
\]
is precompact in the space $\dot{H}^{s_p} \times \dot{H}^{s_p-1}(\Rm^3)$. The function $\lambda(t)$ is called the frequency scale function,
because the solution $u(t)$ at time $t$ concentrates around the frequency $\lambda(t)$ by the compactness.\\
Please note that here we use the radial condition, thus the only available symmetries are scalings. If we did not assume the radial condition, similar
results would still hold but the symmetries would include translations besides scalings.
\subsection{Existence of Minimal Blow-up Solution}
\begin{theorem} \emph{\textbf{(Minimal blow-up solution)}}
Assume that the main theorem failed. Then there would exist
a solution $u(x,t)$ with a maximal lifespan $I$ such that
\[
\sup_{t\in I} \|(u(t), \partial_t u(t))\|_{\dot{H}^{s_p} \times \dot{H}^{s_p -1}} < \infty,
\]
$u$ blows up in the positive direction at time $T_+ \leq + \infty$ with
\[
 \|u\|_{S([0, T_+))} = \infty.
\]
In addition, $u$ is almost periodic modulo scalings with a frequency scale function $\lambda(t)$.
It is minimal in the following sense, if $v$ is another solution with a maximal lifespan $J$ and
\[
 \sup_{t \in J} \|(v(t), \partial_t v(t))\|_{\dot{H}^{s_p} \times \dot{H}^{s_p -1}} <
 \sup_{t\in I} \|(u(t), \partial_t u(t))\|_{\dot{H}^{s_p} \times \dot{H}^{s_p -1}},
\]
then $v$ is a global solution in time and scatters.
\end{theorem}
\noindent The main tool to obtain this result is the profile decomposition. One could follow the argument in \cite{kenig2} in order to find a proof.
In that paper C.E.Kenig and F.Merle deal with the cubic defocusing NLS under similar assumptions.
\subsection{Three enemies}
Since the frequency scale function $\lambda(t)$ plays an important role in the further discussion, it is
helpful if we could make additional assumptions on this function. It turns out that we could reduce
the whole problem into the following three special cases. This method of three enemies was
introduced in R.Killip, T.Tao and M.Visan's paper \cite{tao}.

\begin{theorem} \emph{\textbf{(Three enemies)}}
Suppose our main theorem failed, then there would exist a minimal blow-up solution $u$ satisfying all the conditions we mentioned
in the previous theorem, so that one of the following three assumptions on its lifespan $I$ and frequency scale function $\lambda(t)$ holds
\begin{itemize}
\item (I) \emph{\textbf{(Soliton-like case)}} $I = \Rm$ and $\lambda(t) = 1$.
\item (II) \emph{\textbf{(High-to-low frequency cascade)}} $I=\Rm$, $\lambda(t) \leq 1$ and
\[
 \liminf_{t \rightarrow \pm \infty} \lambda(t) = 0.
\]
\item(III) \emph{\textbf{(Self-similar case)}} $I = \Rm^{+}$ and $\lambda(t) = 1/t$.
\end{itemize}
\end{theorem}
\noindent Please note that the minimal blow-up solution $u$ here could be different from the one we found in the previous theorem.
But we can always manufacture a minimal blow-up solution in one of these three cases from the original one.
One can follow the method used in the paper \cite{tao} to verify this theorem.
\subsection{Further Compactness Results}
Fix a radial cutoff function $\varphi(x) \in C^{\infty}(\Rm^3)$ with the following properties.
\[
 \varphi(x) \left\{\begin{array}{l} =0, \,\,\,\, |x| \leq 1/2;\\
 \in [0,1],\,\,\,\, 1/2 \leq |x| \leq 1;\\
 = 1,\,\,\,\, |x| \geq 1.\end{array}\right.
\]
Given a minimal blow-up solution $u$ mentioned above and its frequency scale function $\lambda(t)$, we have the following
propositions by a compactness argument.
\begin{proposition} \label{local compactness 1}
Let $u$ be a minimal blow-up solution with a maximal lifespan $I$ as above.
There exist constants $d, C'>0$ and $C_1>1$ independent of $t$ such that \\
(i) The interval $[t-d \lambda^{-1}(t), t+ d \lambda^{-1}(t)] \subseteq I$ for all $t \in I$.
In addition, for each $t' \in [t-d \lambda^{-1}(t), t+ d \lambda^{-1}(t)]$, we have
\begin{equation}
 \frac{1}{C_1} \lambda(t) \leq \lambda(t') \leq C_1 \lambda(t).
\label{choiced}
\end{equation}
(ii) The following estimate holds for each $s_p$-admissible pair $(q,r)$ and each $t \in I$.
\[
 \|u\|_{\displaystyle L^q  L^{r} ([t-d \lambda^{-1}(t), t+ d \lambda^{-1}(t)] \times \Rm^3 )} \leq C'.
\]
\end{proposition}
\begin{proposition} \label{local compactness 2}
Given $\eps>0$, there exists $R_1 = R_1 (\eps) > 0$, such that the following inequality holds for each $t \in I$.
\[
 \left\|\left(\varphi\left(\frac{x}{R_1 \lambda^{-1}(t)}\right)u(t),\varphi\left(\frac{x}{R_1 \lambda^{-1}(t)}\right)\partial_t
 u(t) \right)\right\|_{\dot{H}^{s_p} \times \dot{H}^{s_p-1}(\Rm^3)} \leq \eps.
\]
\end{proposition}
\begin{proposition}\label{lower bound of int u p+1 over x}
There exists two constants $R_0, \eta_0>0$, such that the following inequality holds for each $t \in I$. (The constant $d$ is the same constant we used
in proposition \ref{local compactness 1})
\[
\int_{t}^{t+ d\lambda^{-1}(t)} \int_{|x|<R_0 \lambda^{-1}(t)} \frac{|u(x,\tau)|^{p+1}}{|x|} dx d\tau \geq \lambda(t)^{2 -2s_p} \eta_0.
\]
\end{proposition}
\paragraph{Proof} By a compactness argument we obtain that there exist $R_0, \eta_0>0$, so that for all $t \in I$,
\[
 \int_0^d \int_{|x| < R_0} \frac{(\frac{1}{\lambda(t)^{2/(p-1)}} |u(\lambda^{-1}(t)x, \lambda^{-1}(t)\tau + t)|)^{p+1}}{|x|}dxd\tau \geq \eta_0.
\]
This implies
\[
 \int_0^d \int_{|x| < R_0} \frac{ |u(\lambda^{-1}(t)x, \lambda^{-1}(t)\tau + t)|^{p+1}}{\lambda^{-1}(t)|x|}
 \frac{dxd\tau}{\lambda(t)^{\frac{2(p+1)}{p-1}+1}} \geq \eta_0.
\]
\[
 \frac{1}{\lambda(t)^{4/(p-1) -1}} \int_0^d \int_{|x| < R_0} \frac{ |u(\lambda^{-1}(t)x, \lambda^{-1}(t)\tau + t)|^{p+1}}{\lambda^{-1}(t)|x|}
 \frac{dxd\tau}{\lambda(t)^4} \geq \eta_0.
\]
\begin{equation}
\begin{array}{r}
  \displaystyle \int_{t}^{t + d \lambda^{-1}(t)} \int_{|x| < R_0 \lambda^{-1}(t)}
  \frac{|u(x,\tau)|^{p+1}}{|x|} dx d\tau \geq \lambda(t)^{4/(p-1) -1} \eta_0\\
   =  \lambda(t)^{2 -2s_p} \eta_0.
\end{array} \label{lowbound1}
\end{equation}
\subsection{The Duhamel Formula}
The following Duhamel formula will be frequently used in later sections.
\begin{proposition} \emph{\textbf{(The Duhamel formula)}} Let $u$ be a minimal blow-up solution described above with a maximal lifespan $I=(T_-,\infty)$.
Then we have
\begin{eqnarray*}
u(t) &=& \lim_{T\rightarrow +\infty} \int_{t}^{T} \frac{\sin((\tau -t)\sqrt{-\Delta})}{\sqrt{-\Delta}} F(u(\tau)) d\tau;\\% can generalize to F(\tau)
\partial_t u(t) &=& -\lim_{T\rightarrow +\infty} \int_{t}^{T} \cos((\tau -t)\sqrt{-\Delta}) F(u(\tau)) d\tau.
\end{eqnarray*}
\begin{eqnarray*}
u(t) &=& \lim_{T\rightarrow T_-} \int_{T}^{t} \frac{\sin((t-\tau)\sqrt{-\Delta})}{\sqrt{-\Delta}} F(u(\tau)) d\tau;\\
\partial_t u(t) &=& \lim_{T\rightarrow T_-} \int_{T}^{t} \cos((t-\tau)\sqrt{-\Delta}) F(u(\tau)) d\tau.
\end{eqnarray*}
Given a time $t \in I$, these limits are weak limits in the space $\dot{H}^{s_p} \times \dot{H}^{s_p -1}$. If $J$ is a closed interval compactly
supported in $I$, then one could also understand the formula for $u(t)$ as a strong limit in the space $L^q L^r (J\times \Rm^3)$,
as long as $(q,r)$ is an $s_p$-admissible pair with $q \neq \infty$.
\end{proposition}
\paragraph{Remark} Actually we have
\begin{equation}
 \left(\begin{array}{c} \displaystyle \int_{t}^{T} \frac{\sin((\tau -t)\sqrt{-\Delta})}{\sqrt{-\Delta}} F(u(\tau)) d\tau\\
 \displaystyle -\int_{t}^{T} \cos((\tau -t)\sqrt{-\Delta}) F(u(\tau)) d\tau\end{array} \right) =
 \left( \begin{array}{c} u(t)\\ \partial_t u (t) \end{array}\right) -
 S(t-T)\left( \begin{array}{c} u(T)\\ \partial_t u (T) \end{array}\right).
\end{equation}
Thus we only need to show the corresponding limit of the last term is zero in order to verify this formula. Please see
lemma \ref{proof of Duhamel} in the appendix for details.
\section{Energy Estimate near Infinity}
In this section, we will prove the following theorem for a minimal blow-up solution $u(x,t)$. The method was previously used in the
supercritical case of the equation. (please see \cite{km} for more details) In the supercritical case, by the Sobolev embedding,
the energy automatically exists at least locally in the space, for any given time $t \in I$.
In the subcritical case, however, we need to use the approximation techniques.
\begin{theorem} \emph{\textbf{(Energy estimate near infinity)}}\label{energy estimate near infinity}
Let $u(x,t)$ be a minimal blow-up solution as we found in the previous section. Then
$(u(x,t),\partial_t u(x,t)) \in \dot{H}^1 \times L^2 (\Rm \setminus B(0,r))$ for each $r>0$, $t \in I$. Actually we have
\begin{equation}
\int_{r < |x| < 4r} (|\nabla u(x,t)|^2+ |\partial_t u(x,t)|^2) dx \lesssim r^{-2(1-s_p)}.
\end{equation}
\end{theorem}
\subsection{Preliminary Results}
\paragraph{Introduction to $w(r,t)$}
Let $u(x,t)$ be a radial solution of the wave equation
\[
 \partial_t^2 u - \Delta u = F(x,t).
\]
If we define $w(r,t), h(r,t): \Rm^+ \times I \rightarrow \Rm$ so that
\[
 w(|x|,t) = |x| u(x,t);\;\;\; h(|x|,t) = |x| F(x,t);
\]
then we have $w(r,t)$ is the solution of the one-dimensional wave equation
\[
 \partial_t^2 w - \partial_r^2 w = h(r,t).
\]
\begin{lemma} \label{eqn between w and u}
Let $(u(x,t_0), \partial_t u(x,t_0))$ be radial and in the energy space $\dot{H}^1 \times L^2$ locally,
then for any $0 < a < b < \infty$, we have the identity
\[
 \frac{1}{4\pi}\int_{a < |x| < b} (|\nabla u|^2 + |\partial_t u|^2) dx = \left(\int_a^b [(\partial_r w)^2 + (\partial_t w)^2]dr\right) +
 \left(a u^2(a) - b u^2(b)\right)
\]
holds (if we take the value of the functions at time $t_0$).
\end{lemma}
\paragraph{Proof} By direct computation
\begin{eqnarray*}
 \int_a^b [(\partial_r w)^2 + (\partial_t w)^2]dr &=& \int_a^b [(r \partial_r u + u)^2 + (r \partial_t u)^2]dr\\
 &=& \int_a^b [r^2 (\partial_r u)^2 + u^2 + r^2 (\partial_t u)^2]dr + \int_a^b 2r u  \partial_r u dr\\
 &=& \int_a^b [r^2 (\partial_r u)^2 +  r^2 (\partial_t u)^2 + u^2]dr + \int_a^b r d(u^2)\\
 &=& \int_a^b r^2 [(\partial_r u)^2 + (\partial_t u)^2 ]dr + [r u^2]_a^b\\
 &=& \frac{1}{4\pi}\int_{a < |x| < b} (|\nabla u|^2 + |\partial_t u|^2) dx + b u^2(b) - a u^2(a).
\end{eqnarray*}
\begin{lemma} \label{derivative of int w}
Let $w(r,t)$ be a solution to the following equation for $(r,t) \in \Rm^+ \times I$ %w, \partial_t w is continuous in t; h is continuous
\[
  \partial_t^2 w - \partial_r^2 w = h (r,t),
\]
so that $(w, \partial_t w) \in C(I; \dot{H}^1 \times L^2 (R_1<r<R_2))$ for any $0 <R_1 < R_2< \infty$.
Let us define
\[
 z_1 (r,t) = \partial_t w(r, t) - \partial_r w(r, t);
\]
\[
 z_2 (r,t) = \partial_t w(r, t) + \partial_r w(r, t).
\]
Then we have ($M>0$)
\begin{eqnarray}
 \lefteqn{\left| \left(\int_{r_0}^{4r_0} |z_1 (r, t_0)|^2 dr\right)^{1/2} -
\left(\int_{r_0 + M}^{4r_0 + M} |z_1 (r, t_0 + M)|^2 dr\right)^{1/2} \right|}\nonumber\\
&\leq& \left(\int_{r_0}^{4 r_0} \left(\int_0^M h(r +t, t_0 +t) dt\right)^2 dr \right)^{1/2},
\end{eqnarray}
if $t_0, t_0 +M \in I$;
\begin{eqnarray}
  \lefteqn{\left| \left(\int_{r_0}^{4r_0} |z_2 (r, t_0)|^2 dr\right)^{1/2} -
\left(\int_{r_0 + M}^{4r_0 + M} |z_2 (r, t_0 - M)|^2 dr\right)^{1/2} \right|}\nonumber\\
&\leq& \left(\int_{r_0}^{4 r_0} \left(\int_0^M h(r +t, t_0 -t) dt\right)^2 dr \right)^{1/2},
\end{eqnarray}
if $t_0, t_0-M \in I$.
\end{lemma}
\paragraph{Proof} We will assume $w$ has sufficient regularity, otherwise we only need to use
the standard techniques of smooth approximation. Let us define
\[
 z(r,s) = (\partial_t - \partial_r) w(r + s, t_0 + s).
\]
We have
\[
 \partial_s z(r,s) = (\partial_t + \partial_r) (\partial_t - \partial_r) w(r + s, t_0 + s) = h(r+s, t_0 +s).
\]
Thus
\[
 z(r,M) = z(r,0) + \int_0^M h(r +t, t_0 +t) dt.
\]
Applying the triangle inequality, we obtain the first inequality. The second inequality could be proved in a similar way.
\subsection{Smooth approximation}
Let $u(x,t)$ be a minimal blow-up solution. Choose a smooth, nonnegative, radial function $\varphi(x,t)$
supported in the four-dimensional ball $B(0,1) \subset \Rm^4$ such that
\[
 \int_{\Rm^3 \times \Rm} \varphi(x,t) dx dt = 1.
\]
Let $d$ be the number given in proposition \ref{local compactness 1}. If $\eps < d$, we define
(both the functions $u$ and $F(u)$ are locally integrable)
\[
 \varphi_\eps (x,t) = \frac{1}{\eps^4} \varphi(x/\eps, t/\eps);
\]
\[
 u_\eps = u \ast \varphi_\eps;\;\;\; F_\eps = F(u) \ast \varphi_\eps.
\]
This makes $u_\eps (x,t)$ be a smooth solution of the linear wave equation
\[
 \partial_t^2 u_\eps (x,t) - \Delta u_\eps (x,t) = F_\eps (x,t).
\]
with the convergence (using the continuity of $(u(t), \partial_t u(t))$ in the space $\dot{H}^{s_p} \times \dot{H}^{s_p -1}$)
\[
 (u_\eps(t_0), \partial_t u_\eps(t_0)) \rightarrow (u(t_0),\partial_t u(t_0))\;
\hbox{in the space}\; \dot{H}^{s_p} \times \dot{H}^{s_p-1}\; \hbox{for each}\; t_0 \in I
\]
and
\[
 \|(u_\eps(t_0), \partial_t u_\eps(t_0))\|_{\dot{H}^{s_p} \times \dot{H}^{s_p-1}} \lesssim 1.
\]
In addition, if $a -\eps \in I$, we have
\[
 \|F_\eps (x,t)\|_{Z_{s_p}([a,b])} < \infty.
\]
\begin{lemma} \emph{\textbf{(Almost periodic property) %away from time $0$
}} The following set
\[
\left\{\left(\frac{1}{\lambda(t)^{3/2-s_p}} u_\eps \left(\frac{x}{\lambda(t)}, t\right), \frac{1}{\lambda(t)^{5/2-s_p}}
\partial_t u_\eps \left(\frac{x}{\lambda(t)}, t\right)\right): t \in [d+1,\infty)\right\}
\]
is precompact in the space $\dot{H}^{s_p} \times \dot{H}^{s_p -1}$ for each fixed $\eps <d$.
\end{lemma}
\paragraph{Stretch of the proof} Given a sequence $\{t_n\}$, WLOG, we could assume
\[
 \lambda(t_n) \rightarrow \lambda_0 \in [0,1];
\]
\[
 \left(\frac{1}{\lambda(t_n)^{3/2-s_p}} u \left(\frac{x}{\lambda(t_n)}, t_n\right), \frac{1}{\lambda(t_n)^{5/2-s_p}}
\partial_t u \left(\frac{x}{\lambda(t_n)}, t_n\right)\right) \rightarrow (u_0,u_1);
\]
by extracting a subsequence if necessary. Let $\tilde{u}(x,t)$ be the solution of the equation (\ref{eqn}) with initial data $(u_0,u_1)$.
By the long-time perturbation theory we know
\[
 \sup_{t\in [-d,d]} \left\|\left(\begin{array}{c} \frac{1}{\lambda(t_n)^{3/2-s_p}} u \left(\frac{x}{\lambda(t_n)}, t_n+ \frac{t}{\lambda(t_n)}\right)\\
 \frac{1}{\lambda(t_n)^{5/2-s_p}} \partial_t u \left(\frac{x}{\lambda(t_n)}, t_n + \frac{t}{\lambda(t_n)}\right) \end{array}\right)
 - \left(\begin{array}{c} \tilde{u}(t)\\ \partial_t \tilde{u}(t)\end{array} \right)\right\|_{\dot{H}^{s_p} \times \dot{H}^{s_p -1}}
\rightarrow 0.
\]
%Considering our assumption $\lambda(t_n)\rightarrow 0$, we obtain
%\[
% \sup_{t\in [-\eps \lambda(t_n),\eps \lambda(t_n)]}
% \left\|\left(\begin{array}{c} \frac{1}{\lambda(t_n)^{3/2-s_p}} u (\frac{x}{\lambda(t_n)}, t_n+ \frac{t}{\lambda(t_n)})\\
% \frac{1}{\lambda(t_n)^{5/2-s_p}} \partial_t u (\frac{x}{\lambda(t_n)}, t_n + \frac{t}{\lambda(t_n)}) \end{array}\right)
% - \left(\begin{array}{c} u_0\\ u_1\end{array} \right)\right\|_{\dot{H}^{s_p} \times \dot{H}^{s_p -1}}
%\rightarrow 0.
%\]
This implies
\begin{eqnarray*}
 \left(\begin{array}{c}\frac{1}{\lambda(t_n)^{3/2-s_p}} u_\eps (\frac{x}{\lambda(t_n)}, t_n)\\
 \frac{1}{\lambda(t_n)^{5/2-s_p}} \partial_t u_\eps (\frac{x}{\lambda(t_n)}, t_n)\end{array} \right) &=&
 \left[\varphi_{\eps \lambda(t_n)} \ast
 \left(\begin{array}{c} \frac{1}{\lambda(t_n)^{3/2-s_p}} u (\frac{\cdot}{\lambda(t_n)}, t_n+ \frac{\cdot}{\lambda(t_n)})\\
 \frac{1}{\lambda(t_n)^{5/2-s_p}} \partial_t u (\frac{\cdot}{\lambda(t_n)}, t_n + \frac{\cdot}{\lambda(t_n)}) \end{array}\right)\right](x,0)\\
 &=& \left[\varphi_{\eps \lambda(t_n)} \ast \left(\begin{array}{c} \tilde{u}\\ \partial_t \tilde{u}\end{array} \right)\right](x,0) + o(1)\\
 &=& \left\{\begin{array}{l} \left[\varphi_{\eps \lambda_0} \ast \left(\begin{array}{c} \tilde{u}\\ \partial_t \tilde{u}\end{array}\right)\right](x,0)
 + o(1)\;\hbox{if}\; \lambda_0 > 0;\\
 \left(\begin{array}{c} u_0\\ u_1\end{array}\right) + o(1)
 \;\hbox{if}\; \lambda_0 = 0;\end{array}\right.
\end{eqnarray*}
\paragraph{Remark} The error $o(1)$ tends to zero as $n \rightarrow \infty$ in the sense of the $\dot{H}^{s_p} \times \dot{H}^{s_p -1}$ norm.
%\paragraph{Case 2 $\lambda_0 >0$} Let
%\[
% v_0(x) = \lambda_0^{3/2 -s_p} u_0 (\lambda_0 x);\;\; v_1(x) = \lambda_0^{5/2 - s_p} (\lambda_0 x).
%\]
%and $v(x,t)$ be the solution of our equation (\ref{eqn}) with initial data $(v_0,v_1)$. By our assumptions we have
%\[
% (u(x,t_n), \partial_t u(x,t_n)) \rightarrow (v_0, v_1)
%\]
%in the space $\dot{H}^{s_p} \times \dot{H}^{s_p-1}$. Applying the long-time perturbation theory, we obtain
%\[
% \sup_{t\in [-\eps,\eps]}
% \left\|\left(\begin{array}{c} u (x, t_n+ t)\\ \partial_t u (x, t_n + t) \end{array}\right)
% - \left(\begin{array}{c} v(t)\\ \partial_t v(t)\end{array} \right)\right\|_{\dot{H}^{s_p} \times \dot{H}^{s_p -1}}
%\rightarrow 0.
%\]
%Thus
\paragraph{The Duhamel Formula} By the almost periodic property above we know the following Duhamel formula still holds for $u_\eps$ in the sense of
weak limit if $t_0 -\eps \in I$.
\begin{eqnarray*}
 u_\eps (t_0) &=& \int_{t_0}^{+\infty} \frac{\sin((\tau-t_0)\sqrt{-\Delta})}{\sqrt{-\Delta}} F_\eps (x,\tau) d\tau;\\
 \partial_t u_\eps (t_0) &=& - \int_{t_0}^{+\infty} \cos((\tau-t_0)\sqrt{-\Delta}) F_\eps (x,\tau) d\tau.
\end{eqnarray*}
In the soliton-like or high-to-low frequency cascade case, we can also verify the Duhamel formula in the negative
time direction. %The limit here is weak limit in $\dot{H}^{s_p} \times \dot{H}^{s_p -1}$.
%move convolution to its duality, then use dominated convergence theorem.
\paragraph{The idea to prove theorem \ref{energy estimate near infinity}} If we could obtain the following estimate
\begin{equation}
 \int_{r < |x| < 4r} (|\nabla u_\eps (x,t_0)|^2+ |\partial_t u_\eps (x,t_0)|^2) dx \leq C r^{-2(1-s_p)},
\end{equation}
so that the constant $C$ is independent of $r>0$, $t_0\in I$ and $\eps < \eps_0 (r, t_0)$,
then we would be able to prove theorem \ref{energy estimate near infinity}
by letting $\eps$ converge to zero. One could read lemma \ref{app a} if interested in the details of this argument.
\paragraph{Remark} We have to apply the smooth kernel on the whole non-linear term. Because if we just made
the initial data smooth, we would not resume the compactness conditions of the minimal blow-up solution.
\begin{lemma} \label{estimate on ueps and Feps}
If $|x| > 10 \eps$, we have
\[
 |u_\eps(x,t)| \leq \frac{C}{|x|^{2/(p-1)}};\;\;\;\; |F_\eps(x,t)| \leq \frac{C}{|x|^{2p/(p-1)}}.
\]
The constant $C$ depends only on the upper bound $\sup_{t \in I} \|(u,\partial_t u)\|_{\dot{H}^{s_p} \times \dot{H}^{s_p -1}}$.
\end{lemma}
\paragraph{Proof} This comes from the estimate (\ref{estimate of radial hs function}) and an easy computation.

\subsection{Uniform Estimate on $u_\eps$}
In this subsection, we will prove the following lemma. It implies theorem \ref{energy estimate near infinity}
immediately by our argument above. The functions $w_\eps (r,t)$ and $z_{i,\eps} (r,t)$ below are defined as described earlier using $u_\eps(x,t)$.
\begin{lemma} \label{uniform estimate on ueps}
Let $t_0 \in I$ and $r_0 >0$, then for sufficiently small $\eps$, we have
\begin{equation}
 \int_{r_0 < |x| < 4r_0} (|\nabla u_\eps (x,t_0)|^2+ |\partial_t u_\eps (x,t_0)|^2) dx \leq C r_0^{-2(1-s_p)}.
\end{equation}
The constant $C$ could be chosen independent of $t_0, r_0$ and $\eps$.
\end{lemma}
\paragraph{Step 1 Conversion to $w_\eps (r,t)$} First choose $\eps < \min\{r_0/{10},d\}$. If the minimal blow-up solution is a self-similar one,
we also require $\eps < t_0/2$.
Let us apply lemma \ref{eqn between w and u} and lemma \ref{estimate on ueps and Feps}. It is sufficient
to show
\[
 \int_{r_0}^{4 r_0} (|\partial_r w_\eps (r,t_0)|^2+ |\partial_t w_\eps (r,t_0)|^2) dr \leq C r_0^{-2(1-s_p)}.
\]
In other words,
\begin{equation}
 \int_{r_0}^{4 r_0} (|z_{1,\eps} (r,t_0)|^2+ |z_{2,\eps} (r,t_0)|^2) dr \leq C r_0^{-2(1-s_p)}.
\end{equation}
\paragraph{Step 2 Expansion of $z_{1,\eps}$} Let us break $(u_\eps(t), \partial_t u_\eps(t))$ into two pieces.
\begin{eqnarray*}
 u_\eps^{(1)}(t) &=& \int_{t}^{t_0 +100 r_0} \frac{\sin((\tau-t)\sqrt{-\Delta})}{\sqrt{-\Delta}} F_\eps (x,\tau) d\tau;\\
 \partial_t u_\eps^{(1)} (t) &=& - \int_{t}^{t_0 + 100 r_0} \cos((\tau-t)\sqrt{-\Delta}) F_\eps (x,\tau) d\tau.
\end{eqnarray*}
\begin{eqnarray*}
 u_\eps^{(2)}(t) &=& \int_{t_0 +100 r_0}^\infty \frac{\sin((\tau-t)\sqrt{-\Delta})}{\sqrt{-\Delta}} F_\eps (x,\tau) d\tau;\\
 \partial_t u_\eps^{(2)} (t) &=& - \int_{t_0 + 100 r_0}^\infty \cos((\tau-t)\sqrt{-\Delta}) F_\eps (x,\tau) d\tau.
\end{eqnarray*}
These are smooth functions and we have
\[
 (u_\eps (x,t_0), \partial_t u_\eps (x,t_0)) = (u_\eps^{(1)} (x,t_0), \partial_t u_\eps^{(1)} (x,t_0))
+ (u_\eps^{(2)} (x,t_0), \partial_t u_\eps^{(2)} (x,t_0)).
\]
Defining $w_\eps^{(j)}, z_{1,\eps}^{(j)}$ accordingly for $j=1,2$, we have
\[
 z_{1,\eps} (x,t_0) = z_{1,\eps}^{(1)} (x,t_0) + z_{1,\eps}^{(2)} (x,t_0).
\]
\paragraph{Step 3 Short-time Contribution} we have $u_{\eps}^{(1)}$ satisfies the wave equation
\[
 \left\{\begin{array}{l} \partial_t^2 u_{\eps}^{(1)} - \Delta u_{\eps}^{(1)} = F_\eps (x,t), \,\,\,\,
(x,t)\in \Rm^3 \times (t_0^-, +\infty);\\
u_{\eps}^{(1)} |_{t=t_0 + 100 r_0} = 0 \in \dot{H}^{s_p} (\Rm^3);\\
\partial_t u_{\eps}^{(1)} |_{t= t_0 + 100 r_0} = 0 \in \dot{H}^{s_p-1}(\Rm^3).\end{array}\right.
\]
Thus $w_\eps^{(1)}$ is a smooth solution of
\[
 \left\{\begin{array}{l} \partial_t^2 w_{\eps}^{(1)} - \partial_r^2 w_{\eps}^{(1)} = r F_\eps (r,t), \,\,\,\,
(r,t)\in \Rm^+ \times (t_0^-, +\infty);\\
w_{\eps}^{(1)} |_{t=t_0 + 100 r_0} = 0;\\
\partial_t w_{\eps}^{(1)} |_{t= t_0 + 100 r_0} = 0.\end{array}\right.
\]
Applying lemma \ref{derivative of int w} and lemma \ref{estimate on ueps and Feps}, we obtain
\begin{eqnarray*}
 \left(\int_{r_0}^{4r_0} |z_{1,\eps}^{(1)} (r,t_0)|^2 dr \right)^{1/2}
 &\leq& \left(\int_{r_0}^{4 r_0} \left(\int_{0}^{100r_0} (t + r) F_\eps (t + r, t+t_0) dt\right)^2 dr \right)^{1/2}\\
 &\lesssim& \left(\int_{r_0}^{4 r_0} \left(\int_{0}^{100r_0} (t + r)\frac{1}{(t+r)^\frac{2p}{p-1}} d t\right)^2 dr \right)^{1/2}\\
 &\lesssim& \left(\int_{r_0}^{4 r_0} \left(\int_{0}^{100r_0} \frac{1}{(t+ r)^{1 + \frac{2}{p-1}}} dt\right)^2 dr \right)^{1/2}\\
% &\lesssim& \left(\int_{r_0}^{4 r_0} (\frac{1}{r^{\frac{2}{p-1}}})^2 dr \right)^{1/2}\\
 &\lesssim& \left(\int_{r_0}^{4 r_0} \frac{1}{r^{4/(p-1)}} dr \right)^{1/2}\\
% &\lesssim& \left(\frac{1}{r_0^{\frac{4}{p-1}-1}} \right)^{1/2}\\
 &\lesssim& \frac{1}{r_0^{1-s_p}}.
\end{eqnarray*}
\paragraph{Step 4 Long-time Contribution} %The pair $(u_\eps^{(2)} (x,t_0), \partial_t u_\eps^{(2)} (x,t_0))$ is the weak
%$\dot{H}^{s_p} \times \dot{H}^{s_p -1}$ limit of
%\begin{eqnarray*}
%u_{\eps, T}^{(2)}(t_0) &=& \int_{t_0 +100 r_0}^T \frac{\sin((\tau-t_0)\sqrt{-\Delta})}{\sqrt{-\Delta}} F_\eps (x,\tau) d\tau;\\
% \partial_t u_{\eps, T}^{(2)} (t_0) &=& - \int_{t_0 + 100 r_0}^T \cos((\tau-t_0)\sqrt{-\Delta}) F_\eps (x,\tau) d\tau.
%\end{eqnarray*}
%We also have its center-cutoff version
%\begin{eqnarray*}
%\tilde{u}_{\eps, T}^{(2)}(t_0) &=& \int_{t_0 +100 r_0}^T \frac{\sin((\tau-t_0)\sqrt{-\Delta})}{\sqrt{-\Delta}} \chi(x,\tau) F_\eps (x,\tau) d\tau;\\
% \partial_t \tilde{u}_{\eps, T}^{(2)} (t_0) &=& - \int_{t_0 + 100 r_0}^T \cos((\tau-t_0)\sqrt{-\Delta}) \chi(x,\tau) F_\eps (x,\tau) d\tau.
%\end{eqnarray*}
Let us define a cutoff function $\chi(x,t)$ to be the characteristic function of the region $\{(x,t): |x|> t-t_0 -50r_0\}$.
By lemma \ref{estimate on ueps and Feps}, we know
%By strong Huygens's principal, we have the pair $(\tilde{u}_{\eps, T}^{(2)}(t_0), \partial_t \tilde{u}_{\eps, T}^{(2)} (t_0))$ is the same as
%the pair $(u_{\eps, T}^{(2)}(t_0), \partial_t u_{\eps, T}^{(2)} (t_0))$ in the ball $B(0,5r_0)$, thus it converges to
%$(u_{\eps}^{(2)}(t_0), \partial_t u_{\eps}^{(2)} (t_0))$ weakly in the ball $B(0,5r_0)$. On the other hand, it also converges strongly in the
%energy space. It is sufficient to consider
\begin{eqnarray*}
\|\chi F_\eps\|_{L^1 L^2 ([t_0 + 100 r_0, \infty)\times \Rm^3)} &=& \int_{t_0 +100 r_0}^\infty \left(
\int_{|x|>t-t_0 -50r_0}|F_\eps|^2 dx \right)^{1/2}dt\\
&\lesssim& \int_{t_0 +100 r_0}^\infty \left(\int_{|x|>t-t_0 -50r_0} \frac{1}{|x|^{4p/(p-1)}} dx \right)^{1/2}dt\\
&\lesssim& \int_{t_0 +100 r_0}^\infty \left(\frac{1}{|t-t_0 -50 r_0|^{1 + 4/(p-1)}}\right)^{1/2}dt\\
&\lesssim& \int_{t_0 +100 r_0}^\infty \frac{1}{|t-t_0 -50 r_0|^{\frac{1}{2} + \frac{2}{p-1}}}dt\\
&\lesssim& \frac{1}{r_0^{1-s_p}}
\end{eqnarray*}
%By Strichartz estimate, the pair $(\tilde{u}_{\eps, T}^{(2)}(t_0), \partial_t \tilde{u}_{\eps, T}^{(2)} (t_0))$ converges strongly to a pair
%$(\tilde{u}_{0,\eps}, \tilde{u}_{1,\eps})$ in the space $\dot{H}^1 \times L^2$ with
%\[
% \|\tilde{u}_{0,\eps}, \tilde{u}_{1,\eps})\|_{\dot{H}^{1} \times L^2 (\Rm^3)} \lesssim \frac{1}{r_0^{1-s_p}}.
%\]
%Considering both the strong and weak limit of $(\tilde{u}_{\eps, T}^{(2)}(t_0), \partial_t \tilde{u}_{\eps, T}^{(2)} (t_0))$, we have
%\[
% (u_{\eps}^{(2)}(t_0), \partial_t u_{\eps}^{(2)} (t_0)) = (\tilde{u}_{0,\eps}, \tilde{u}_{1,\eps})\; \hbox{in the ball}\; B(0,4r_0).
%\]
Applying lemma \ref{Strong Huygens in Duhamel}, we obtain
\[
 \int_{r_0 <|x|<4r_0} (|\nabla u_\eps^{(2)}(x,t_0)|^2 + |\partial_t u_\eps^{(2)}(x,t_0)|^2) dx \lesssim r_0^{2(s_p-1)}.
\]
Applying lemma \ref{eqn between w and u} and using the fact (plus (\ref{estimate of radial hs function}))
\begin{eqnarray*}
 &&\|(u_{\eps}^{(2)}(t_0), \partial_t u_{\eps}^{(2)} (t_0))\|_{\dot{H}^{s_p} \times \dot{H}^{s_p-1}}\\
 &=& \left\|S(-100r_0)\left(\begin{array}{c} u_{\eps} (t_0+100r_0)\\ \partial_t u_{\eps} (t_0+100 r_0)\end{array}
 \right)\right\|_{\dot{H}^{s_p} \times \dot{H}^{s_p-1}}\\
 &=& \|(u_{\eps} (t_0+100r_0), \partial_t u_{\eps} (t_0+100 r_0))\|_{\dot{H}^{s_p} \times \dot{H}^{s_p-1}}\\
 &\leq& \sup_I \|(u, \partial_t u)\|_{\dot{H}^{s_p} \times \dot{H}^{s_p-1}} \lesssim 1,
\end{eqnarray*}
we obtain
\[
 \int_{r_0}^{4r_0} (|\partial_r w_\eps^{(2)}(r,t_0)|^2 + |\partial_t w_\eps^{(2)}(r,t_0)|^2) dr \lesssim r_0^{2(s_p-1)}.
\]
\[
 \int_{r_0}^{4r_0} |z_{1,\eps}^{(2)} (r,t_0)|^2 dr \lesssim r_0^{2(s_p-1)}.
\]
Combining with the estimate for $z_{1,\eps}^{(1)}$, we have
\[
 \int_{r_0}^{4r_0} |z_{1,\eps} (r,t_0)|^2 dr \lesssim r_0^{2(s_p-1)}.
\]
\paragraph{Step 5 Estimate of $z_{2,\eps}$} We also need to consider $z_{2,\eps}$. In the soliton-like case or the high-to-low frequency cascade case,
this could be done in exactly the same way as $z_{1,\eps}$. Now let us consider the self-similar case.
\begin{lemma} Let $u$ be a self-similar minimal blow-up solution. If $t_0 \leq 0.3 r_0$,
Then $(u(t_0),\partial_t u(t_0))$ is in $\dot{H}^1 \times L^2 (|x|>0.9 r_0)$ with
\[
 \int_{|x|> 0.9 r_0} (|\nabla u(x,t_0)|^2 + |\partial_t u(x,t_0)|^2) dx \lesssim r_0^{2(s_p-1)}.
\]
\end{lemma}
\paragraph{Proof} We have (the Duhamel Formula)
\begin{eqnarray*}
 u(t_0) &=& \int_{0^+}^{t_0} \frac{\sin((t_0-t)\sqrt{-\Delta})}{\sqrt{-\Delta}} F (x,t) dt;\\
 \partial_t u(t_0) &=& \int_{0^+}^{t_0} \cos((t_0-t)\sqrt{-\Delta}) F (x,t) dt.
\end{eqnarray*}
\begin{eqnarray*}
 \tilde{u}_0 &=& \int_{0^+}^{t_0} \frac{\sin((t_0-t)\sqrt{-\Delta})}{\sqrt{-\Delta}} \chi(|x|>0.5 r_0) F (x,t) dt;\\
 \tilde{u}_1 &=& \int_{0^+}^{t_0} \cos((t_0-t)\sqrt{-\Delta}) \chi(|x|>0.5 r_0) F (x,t) dt.
\end{eqnarray*}
A straightforward computation shows $\|\chi F\|_{L^1 L^2 ((0^+, t_0)\times \Rm^3)} \lesssim r_0^{s_p-1}$. This means
$(\tilde{u}_0, \tilde{u}_1)$ is in the space $\dot{H}^1 \times L^2 (\Rm^3)$ with a norm $\lesssim r_0^{s_p-1}$. By strong Huygens's
principal we can repeat the argument we used in lemma \ref{Strong Huygens in Duhamel} and obtain
\[
 (u(t_0), \partial_t u(t_0)) = (\tilde{u}_0, \tilde{u}_1)\; \hbox{in the region}\; \Rm^3\setminus B(0, 0.9r_0).
\]
This completes the proof.
\begin{lemma} \label{close to 0}
Let $u$ be a self-similar solution. If $t_0 \leq 0.2 r_0$ and $\eps < t_0/2$, then we have
\[
  \int_{r_0 < |x| < 4 r_0} (|\nabla u_\eps (x,t_0)|^2 + |\partial_t u_\eps (x,t_0)|^2) dx \lesssim r_0^{2(s_p-1)}.
\]
\end{lemma}
\paragraph{Proof} We have $\nabla u_\eps = \varphi_\eps \ast \nabla u$, thus $|\nabla u_\eps| \leq \varphi_\eps \ast |\nabla u|$.
Thus (we have $\eps < 0.1 r_0$)
\[
 \int_{r_0 < |x| < 4 r_0} |\nabla u_\eps (x, t_0)|^2 dx \leq \sup_{t \in[t_0 -\eps, t_0 +\eps]} \int_{0.9r_0 < |x|< 4.1r_0} |\nabla u (x,t)|^2 dx
 \lesssim  r_0^{2(s_p-1)}
\]
by our previous lemma. The other term could be estimated using the same way.
\paragraph{Remark} By lemma \ref{eqn between w and u} and lemma \ref{estimate on ueps and Feps}, this lemma implies (if $t_0 \leq 0.2 r_0$)
\begin{equation}
 \int_{r_0}^{4r_0} (|\partial_r w_{\eps} (r,t_0)|^2 + |\partial_t w_\eps (r,t_0)|^2) dr \lesssim r_0^{2(s_p-1)}. \label{lin12}
\end{equation}
In the self-similar case, let us recall that we always choose $\eps < \min\{r_0/10, t_0/2, d\}$. By lemma \ref{close to 0} and its remark,
we only need to consider the case $t_0 > 0.2 r_0$ in order to estimate $z_{2,\eps}$. Applying lemma \ref{derivative of int w}, we have
\begin{eqnarray*}
 \left(\int_{r_0}^{4r_0} |z_{2,\eps} (r,t_0)|^2 dr\right)^{1/2} &\leq&
 \left(\int_{t_0 + 0.8 r_0}^{t_0 + 3.8 r_0} |z_{2,\eps} (r, 0.2 r_0)|^2 dr\right)^{1/2}\\
 && +  \left(\int_{r_0}^{4 r_0} \left(\int_{0}^{t_0 - 0.2 r_0} (t + r) F_\eps (t + r, t_0 - t) dt\right)^2 dr \right)^{1/2}
\end{eqnarray*}
The first term is dominated by $r_0^{s_p-1}$ because of (\ref{lin12}). We can gain the same upper bound for the second term by
a basic computation similar to the one we used for $z_{1,\eps}$.
\paragraph{Step 6 Conclusion} Now we combine the estimates for $z_{1,\eps}$ and $z_{2,\eps}$ thus conclude our lemma \ref{uniform estimate on ueps}
and theorem \ref{energy estimate near infinity}. Applying lemma \ref{eqn between w and u}, we obtain
\begin{proposition}\label{int w z estimate} Let $u(x,t)$ be a minimal blow-up solution as above, we have
\begin{eqnarray*}
 \int_{r_0}^{4r_0} (|\partial_r w(r,t_0)|^2 + |\partial_t w(r,t_0)|^2) dr &\lesssim& r_0^{2(s_p -1)}.\\
 \int_{r_0}^{4r_0} (|z_1(r,t_0)|^2 + |z_2(r,t_0)|^2) dr &\lesssim& r_0^{2(s_p -1)}.
\end{eqnarray*}
\end{proposition}
\section{Recurrence Process}
\paragraph{Starting Point} Let $u(x,t)$ be a minimal blow-up solution of (\ref{eqn}) as we obtained in the section of compactness
process with a frequency scale function $\lambda(t)$.
In addition, the following set is precompact in the space $\dot{H}^s \times \dot{H}^{s-1}(\Rm^3)$ for some $s \in [s_p,1)$.
\[
 \left\{\left(\frac{1}{\lambda(t)^{3/2-s_p}} u\left(\frac{x}{\lambda(t)}, t\right), \frac{1}{\lambda(t)^{5/2-s_p}}
 \partial_t u\left(\frac{x}{\lambda(t)}, t\right)\right): t \in I\right\}
\]
In this section we will try to gain higher regularity than $\dot{H}^s \times \dot{H}^{s-1}$ by assuming the conditions above.
\subsection{Setup and Technical Lemmas}
\paragraph{Definition} Let us define
\begin{eqnarray*}
 S(A) &=& \sup_{t \in I} (\lambda(t))^{s_p -s} \|u_{>\lambda(t) A}\|_{Y_s ([t, t+ d \lambda^{-1}(t)])};\\
 N(A) &=& \sup_{t \in I} (\lambda(t))^{s_p -s} \|P_{>\lambda(t) A}F(u)\|_{Z_s ([t, t+ d \lambda^{-1}(t)])}.
\end{eqnarray*}
By our assumptions on compactness and proposition \ref{local compactness 1}, we have
\[
 \left\{\left(\frac{1}{\lambda(t)^{3/2-s_p}} u\left(\frac{x}{\lambda(t)}, t + \frac{\tau}{\lambda(t)}\right), \frac{1}{\lambda(t)^{5/2-s_p}}
 \partial_t u\left(\frac{x}{\lambda(t)}, t+ \frac{\tau}{\lambda(t)}\right)\right): \tau \in [0,d], t \in I \right\}
\]
is uniformly bounded (and precompact) in the space $\dot{H}^s \times \dot{H}^{s-1}$. By the local theory of our equation with initial data in
$\dot{H}^s \times \dot{H}^{s-1}$(If $s=s_p$, please see proposition \ref{local compactness 2} and the long-time perturbation theory, otherwise see
theorem \ref{local existence in Hs} and \ref{perturbation theory in Hs}), we have
\[
 \left\{\frac{1}{\lambda(t)^{3/2-s_p}} u\left(\frac{x}{\lambda(t)}, t + \frac{\tau}{\lambda(t)}\right), \tau \in [0,d]: t \in I \right\}
\]
is precompact in the space $Y_s ([0,d])$. Thus we have
\[
 \left\|\frac{1}{\lambda(t)^{3/2-s_p}} u\left(\frac{x}{\lambda(t)}, t + \frac{\tau}{\lambda(t)}\right)\right\|_{Y_s ([0,d])} \lesssim 1,
\]
and
\[
 \left\|P_{> A} \frac{1}{\lambda(t)^{3/2-s_p}} u\left(\frac{x}{\lambda(t)}, t + \frac{\tau}{\lambda(t)}\right)\right\|_{Y_s ([0,d])}
 \rightrightarrows 0
\]
as $A \rightarrow \infty$. If we rescale the first inequality back, we obtain
\[
 (\lambda(t))^{s_p -s} \|u\|_{Y_s ([t, t+ d \lambda^{-1}(t)])} \lesssim 1 \Rightarrow
 (\lambda(t))^{s_p -s} \|F(u)\|_{Z_s ([t, t+ d \lambda^{-1}(t)])} \lesssim 1,
\]
which implies that $S(A)$ and $N(A)$ are uniformly bounded. In the similar way we can show $S(A)$ converges to zero as $A \rightarrow \infty$, using
the uniform convergence above.
\begin{lemma} \emph{\textbf{Bilinear Estimate}} \label{bilinear estimate}
Suppose $u_i$ satisfies the following linear wave equation on the time interval $I = [0,T]$, $i=1,2$,
\[
 \partial_t^2 u_i - \Delta u_i = F_i (x,t),
\]
with the initial data $(u_i|_{t=0}, \partial_t u_i|_{t=0}) = (u_{0,i}, u_{1,i})$. Then
\begin{eqnarray*}
 S &=& \|(P_{>R} u_1) (P_{<r}u_2)\|_{\displaystyle L^\frac{p}{s + 1 -(2p-2) (s -s_p)} L^{\frac{p}{2-s}} (I \times \Rm^3)} \\
 & \lesssim & (\frac{r}{R})^\sigma \left(\|(u_{0,1}, u_{1,1})\|_{\dot{H}^{s} \times \dot{H}^{s-1}} + \|F_1\|_{Z_s (I)}\right) \times \left(\|(u_{0,2}, u_{1,2})\|_{\dot{H}^{s} \times \dot{H}^{s-1}} + \|F_2\|_{Z_s (I)}\right).
\end{eqnarray*}
Here the number $\sigma$ is an arbitrary positive constant satisfying
\begin{equation}
 \sigma \leq 3 \left(\frac{1}{2} - \frac{s+ 1 -(2p-2)(s-s_p)}{2p} - \frac{2 -s}{2p}\right),\;\;
\sigma < 3\times \frac{2-s}{2p}.
\label{sigma}
\end{equation}
\end{lemma}
\paragraph{Remark} We can always choose
\[
 \sigma = \sigma(p) = \frac{3 \min\{p-3, 1\}}{2p} > 0.
\]
\paragraph{Proof} By the Strichartz estimate
\begin{eqnarray*}
 &&\|(P_{>R}) u_1\|_{\displaystyle L^\frac{2p}{s + 1 -(2p-2) (s -s_p)} L^{1/{(\frac{2-s}{2p} + \frac{\sigma}{3})}}}\\
 &\lesssim& \|(D_x^{-\sigma} P_{>R} u_{0,1}, D_x^{-\sigma} P_{>R} u_{1,1})\|_{\dot{H}^{s} \times \dot{H}^{s-1}} + \|D_x^{-\sigma} P_{>R} F_1\|_{Z_s(I)}.
\end{eqnarray*}
\begin{eqnarray*}
 &&\|(P_{<r}) u_2\|_{\displaystyle L^\frac{2p}{s + 1 -(2p-2) (s -s_p)} L^{1/{(\frac{2-s}{2p} - \frac{\sigma}{3})}}}\\
 &\lesssim& \|(D_x^{\sigma} P_{<r} u_{0,2}, D_x^{\sigma} P_{<r} u_{1,2})\|_{\dot{H}^{s} \times \dot{H}^{s-1}} + \|D_x^{\sigma} P_{<r} F_2\|_{Z_s(I)}.
\end{eqnarray*}
Our choice of $\sigma$ makes sure that the pairs above are admissible. Thus we have
\begin{eqnarray*}
 && \|(P_{>R} u_1) (P_{<r}u_2)\|_{\displaystyle L^\frac{p}{s + 1 -(2p-2) (s -s_p)} L^{\frac{p}{2-s}}} \\
&\lesssim& \|(P_{>R}) u_1\|_{\displaystyle L^\frac{2p}{s + 1 -(2p-2) (s -s_p)} L^{1/{(\frac{2-s}{2p} + \frac{\sigma}{3})}}}
\|(P_{<r}) u_2\|_{\displaystyle L^\frac{2p}{s + 1 -(2p-2) (s -s_p)} L^{1/{(\frac{2-s}{2p} - \frac{\sigma}{3})}}}\\
&\lesssim& \left(\|(D_x^{-\sigma} P_{>R} u_{0,1}, D_x^{-\sigma} P_{>R} u_{1,1})\|_{\dot{H}^{s} \times \dot{H}^{s-1}} + \|D_x^{-\sigma} P_{>R} F_1\|_{Z_s(I)}\right)\\
&& \times \left(\|(D_x^{\sigma} P_{<r} u_{0,2}, D_x^{\sigma} P_{<r} u_{1,2})\|_{\dot{H}^{s} \times \dot{H}^{s-1}} + \|D_x^{\sigma} P_{<r} F_2\|_{Z_s(I)}\right)\\
&\lesssim& (\frac{1}{R})^{\sigma} \left(\|( P_{>R} u_{0,1}, P_{>R} u_{1,1})\|_{\dot{H}^{s} \times \dot{H}^{s-1}} + \| P_{>R} F_1\|_{Z_s(I)}\right)\\
&& \times r^\sigma \left(\|( P_{<r} u_{0,2}, P_{<r} u_{1,2})\|_{\dot{H}^{s} \times \dot{H}^{s-1}} + \|P_{<r} F_2\|_{Z_s (I)}\right)\\
& \lesssim & \hbox{the right hand}.
\end{eqnarray*}
\begin{lemma} \label{lemma2}
Let $u(x,t)$ be a function defined on $I \times \Rm^3$, such that $\hat{u}$ is supported in the ball $B(0,r)$ for each $t \in I$, then
\[
 \|P_{>R} F(u(x,t))\|_{\displaystyle L^\frac{2}{s + 1 -(2p-2) (s -s_p)} L^{\frac{2}{2-s}} (I \times \Rm^3)}
\lesssim (\frac{r}{R})^2 \|u\|_{Y_s (I)}^p.
\]
\end{lemma}
\paragraph{Proof}
\begin{eqnarray*}
&& \|P_{>R} F(u(x,t))\|_{\displaystyle L^\frac{2}{s + 1 -(2p-2) (s -s_p)} L^{\frac{2}{2-s}}(I \times \Rm^3)}\\
&\lesssim& \frac{1}{R^2} \|P_{>R} \Delta_x F(u(x,t))\|_{\displaystyle L^\frac{2}{s + 1 -(2p-2) (s -s_p)} L^{\frac{2}{2-s}}(I \times \Rm^3)}\\
&\lesssim& \frac{1}{R^2} \|\Delta_x F(u(x,t))\|_{\displaystyle L^\frac{2}{s + 1 -(2p-2) (s -s_p)} L^{\frac{2}{2-s}}(I \times \Rm^3)}\\
&\lesssim& \frac{1}{R^2} \|p (\Delta_x u)|u|^{p-1} + p(p-1)|\nabla_x u|^2 |u|^{p-3}u \|_{\displaystyle L^\frac{2}{s + 1 -(2p-2) (s -s_p)} L^{\frac{2}{2-s}}}\\
&\lesssim& \frac{1}{R^2} \left(\|\Delta_x u\|_{Y_s (I)} \|u\|_{Y_s (I)}^{p-1} + \|\nabla_x u\|_{Y_s (I)}^2 \|u\|_{Y_s (I)}^{p-2} \right)\\
&\lesssim& \frac{r^2}{R^2} \|u\|_{Y_s (I)}^p.
\end{eqnarray*}
\begin{lemma} \label{L infinity estimate}
Let $v(t)$ be a long-time contribution in the Duhamel formula as below
\[
v(t_0) = \int_{T_1}^{T_2} \frac{\sin((t-t_0)\sqrt{-\Delta})}{\sqrt{-\Delta}} F(u(t))dt,
\]
then for any $t_0 < T_1 < T_2$, we have
\[
\|v(t_0)\|_{L^\infty (\Rm^3)} \lesssim (T_1 -t_0)^{-2/(p-1)}.
\]
\end{lemma}
\paragraph{Proof} Using the explicit expression of the wave kernel in dimension $3$, we obtain
\begin{eqnarray*}
 && \left|\left(\int_{T_1}^{T_2} \frac{\sin((t-t_0)\sqrt{-\Delta})}{\sqrt{-\Delta}} F(u(t))dt\right)(x)\right|\\
&=& \left|\int_{T_1}^{T_2} \int_{|y-x| = t-t_0} \frac{1}{4\pi (t-t_0)} F(u(y,t)) dS(y)dt\right|\\
&\lesssim& \int_{T_1}^{T_2} \int_{|y-x| = t-t_0} \frac{1}{4\pi (t-t_0)} |u(y,t)|^{p} dS(y)dt\\
&\lesssim& \int_{T_1}^{T_2} \int_{|y-x| = t-t_0} \frac{1}{(t-t_0)} \frac{1}{|y|^\frac{2p}{p-1}} dS(y)dt.
\end{eqnarray*}
In the last step, we use the estimate (\ref{estimate of radial hs function}) for radial $\dot{H}^{s_p}$ functions.
If $|x| \leq \frac{1}{2} (T_1 - t_0)$, then on the sphere for the integral we have
\[
 |y| \geq |t-t_0| - |x| \geq \frac{1}{2} (t-t_0).
\]
Thus for these small $x$,
\begin{eqnarray*}
 && \left|\left(\int_{T_1}^{T_2} \frac{\sin((t-t_0)\sqrt{-\Delta})}{\sqrt{-\Delta}} F(u(t))dt\right)(x)\right|\\
&\lesssim& \int_{T_1}^{T_2} \int_{|y-x| = t-t_0} \frac{1}{(t-t_0)} \frac{1}{(t-t_0)^{2p/(p-1)}} dS(y)dt\\
&\lesssim& \int_{T_1}^{T_2} \int_{|y-x| = t-t_0} \frac{1}{(t-t_0)^{3 + 2/(p-1)}} dS(y)dt\\
&\lesssim& \int_{T_1}^{T_2} \frac{(t-t_0)^2}{(t-t_0)^{3 + 2/(p-1)}} dt\\
&\lesssim& \int_{T_1}^{T_2} \frac{1}{(t-t_0)^{1 + 2/(p-1)}} dt\\
&\lesssim& (T_1 - t_0)^{-2/(p-1)}.
\end{eqnarray*}
On the other hand, if $x \geq \frac{1}{2} (T_1 -t_0)$, by (\ref{uni hs es}) we have
\begin{eqnarray*}
 && \left|\left(\int_{T_1}^{T_2} \frac{\sin((t-t_0)\sqrt{-\Delta})}{\sqrt{-\Delta}} F(u(t))dt\right)(x)\right|\\
 &\lesssim& \frac{1}{|x|^{2/(p-1)}} \left\|\int_{T_1}^{T_2} \frac{\sin((t-t_0)\sqrt{-\Delta})}{\sqrt{-\Delta}} F(u(t)) dt\right\|_{\dot{H}^{s_p}}\\
&\lesssim& \frac{1}{(T_1 - t_0)^{2/(p-1)}}.
\end{eqnarray*}
Combining these two cases, we finish our proof.
\begin{lemma} \label{gain decay from recurrence}
Suppose $S(A)$ is a nonnegative function defined in $\Rm^+$ satisfying $S(A) \rightarrow 0$ as $A \rightarrow \infty$.
In addition, there exist $\alpha,\beta \in (0,1)$ and $l, \omega > 0$ with
\[
 l \alpha + \beta > 1,
\]
such that
\begin{equation}
 S(A) \lesssim S(A^\beta)S^l (A^\alpha) + A^{-\omega} \label{recurrence a}
\end{equation}
is true for each sufficiently large $A$. Then
\[
 S(A) \lesssim A^{-\omega}
\]
for each sufficiently large $A$.
\end{lemma}
\paragraph{Proof} Let us first choose two constants $l^-$ and $\omega^-$, which are slightly smaller than $l$ and $\omega$ respectively, such that
the inequality $l^- \alpha + \beta >1$ still holds. By the conditions given, we can find a constant $A_0\gg 1$, such that
the following inequalities hold
\begin{equation}
 S(A) \leq \frac{1}{2} S(A^\beta) S^{l^-} (A^\alpha) + \frac{1}{2} A^{-\omega^-}\; \hbox{if}\; A \geq A_0. \label{recurrence a1}
\end{equation}
\[
 S(A) < 1/2\; \hbox{if}\; A\geq A_0^{1/\alpha}.
\]
Using the second inequality above, we know the following inequality holds for all $A \in [A_0^{\alpha}, A_0]$ if $\omega_1$ is sufficiently small
\begin{equation}
 S(A) \leq A^{-\omega_1}. \label{decay of s a1}
\end{equation}
Fix such a small constant $\omega_1 \leq \omega^-$.
We will show that the inequality (\ref{decay of s a1}) above holds for each $A \geq A_0^{\alpha}$ by an induction.
We already know this is true for $A \in [A_0^{\alpha}, A_0]$. If $A \in [A_0, A_0^{1/\beta}]$, the inequality (\ref{recurrence a1}) implies
\begin{eqnarray*}
 S(A) &\leq& \frac{1}{2} S(A^\beta) S^{l^-} (A^\alpha) + \frac{1}{2} A^{-\omega^-}\\
 &\leq& \frac{1}{2}(A^\beta)^{-\omega_1} \left((A^\alpha)^{-\omega_1}\right)^{l^-} + \frac{1}{2} A^{-\omega^-}\\
 &\leq& \frac{1}{2} \left(A^{-\omega_1}\right)^{\beta + l^- \alpha} + \frac{1}{2} A^{-\omega_1}\\
 &\leq& A^{-\omega_1}.
\end{eqnarray*}
Here we use the fact that $A^\alpha, A^\beta \in [A_0^\alpha, A_0]$ if $A$ satisfies our assumption. Conducting an induction,
we can show the inequality holds for each $A \in [A_0^{(1/\beta)^n}, A_0^{(1/\beta)^{n+1}}]$ if $n$ is a nonnegative integer. In summary,
the inequality (\ref{decay of s a1}) is true for each $A \geq A_0^\alpha$. Plugging this back in the original recurrence formula (\ref{recurrence a}),
we obtain for sufficiently large $A$,
\[
 S(A) \lesssim A^{-\omega_1 (\beta + l \alpha)} + A^{-\omega} \lesssim A^{-\min\{\omega_1(\beta + l\alpha), \omega\}},
\]
which indicates faster decay than $A^{-\omega_1}$. Iterating the argument if necessary, we gain the decay $S(A) \lesssim A^{-\omega}$ and finish the
proof.
\subsection{Recurrence Formula}
Under our setting in this section, given $0 < \alpha < \beta < 1$ and $\eps_1 > 0$, we have the following recurrence formula for sufficiently large $A$
\begin{eqnarray}
 N(A) &\lesssim& S(A^\beta) S^{p-1}(A^\alpha) + A^{-(\beta - \alpha)\sigma(p)} + A^{-2(1-\beta)};\label{recurrence1}\\
 S(A) &\lesssim& N(A^{1-\eps_1}) + A^{-\sigma_1(p)}.\label{recurrence2}
\end{eqnarray}
The constants $\sigma(p), \sigma_1 (p)$ depend on $p$ but nothing else.
\paragraph{Proof of the first inequality} In the following argument, all the space-time norms are taken in $[t, t+d \lambda^{-1}(t)] \times \Rm^3$.
\begin{eqnarray*}
 \lefteqn{\|P_{> \lambda(t)A}(F(u))\|_{Z_s}}\\
  &\lesssim& \lambda(t)^{-(p-1)(s-s_p)} \|P_{> \lambda(t)A}F(u)\|_{\displaystyle L^\frac{2}{s + 1 -(2p-2) (s -s_p)}
 L^{\frac{2}{2-s}}}\\
 &\leq& \lambda(t)^{-(p-1)(s-s_p)} \|P_{> \lambda(t)A}(F(u_{\leq A^\beta \lambda(t)}))\|_{\displaystyle L^\frac{2}{s + 1 -(2p-2) (s -s_p)}
 L^{\frac{2}{2-s}}} \\
 &+&  \lambda(t)^{-(p-1)(s-s_p)} \|P_{> \lambda(t)A}(F(u) - F(u_{\leq A^\beta \lambda(t)}))\|_{\displaystyle L^\frac{2}{s + 1 -(2p-2) (s -s_p)}
 L^{\frac{2}{2-s}}} \\
 &\leq& \lambda(t)^{-(p-1)(s-s_p)}(I_1 + I_2).
\end{eqnarray*}
By lemma \ref{lemma2}, we have
\[
 I_1 \lesssim (\frac{A^\beta}{A})^2 \|u\|_{Y_s}^p \lesssim (\lambda(t))^{p(s-s_p)} A^{-2(1-\beta)}.
\]
For $I_2$, we have (all norms unmarked are $L^\frac{2}{s + 1 -(2p-2) (s -s_p)} L^{\frac{2}{2-s}}([t, t+d \lambda^{-1}(t)] \times \Rm^3)$ norms)
\begin{eqnarray*}
 I_2 &\leq& \left\|P_{>\lambda(t)A}\left[u_{>A^\beta \lambda(t)} \int_0^1 F'(u_{\leq A^\beta \lambda(t)} + \tau u_{> A^\beta \lambda(t)})d\tau\right]
 \right\|\\
&\lesssim& \left\|u_{>A^\beta \lambda(t)} \int_0^1 F'(u_{\leq A^\beta \lambda(t)} + \tau u_{> A^\beta \lambda(t)})d\tau\right\|\\
&\lesssim& \left\| \begin{array}{l} u_{>A^\beta \lambda(t)} \displaystyle \int_0^1 F'(u_{\leq A^\beta \lambda(t)} + \tau u_{> A^\beta \lambda(t)})
d\tau\\
 -u_{>A^\beta \lambda(t)} \displaystyle \int_0^1 F'(u_{ A^\alpha \lambda(t)< \cdotp\leq A^\beta \lambda(t)} + \tau u_{> A^\beta \lambda(t)})
 d\tau \end{array} \right\|\\
&& + \left\|u_{>A^\beta \lambda(t)} \int_0^1 F'(u_{ A^\alpha \lambda(t)<\cdotp\leq A^\beta \lambda(t)} + \tau u_{> A^\beta \lambda(t)})
d\tau\right\|\\
&\lesssim& \left\| u_{>A^\beta \lambda(t)} u_{\leq A^\alpha \lambda(t)}
 \displaystyle \int_0^1 \int_0^1  F''( \tilde{\tau} u_{\leq A^\alpha \lambda(t)} + u_{ A^\alpha \lambda(t)< \cdotp\leq A^\beta \lambda(t)}
 + \tau u_{> A^\beta \lambda(t)} ) d\tau d\tilde{\tau} \right\|\\
&& + \left\|u_{>A^\beta \lambda(t)} \int_0^1 F'(u_{ A^\alpha \lambda(t)<\cdotp\leq A^\beta \lambda(t)} + \tau u_{> A^\beta \lambda(t)})
d\tau\right\|\\
&\lesssim& \left\|u_{>A^\beta \lambda(t)} u_{\leq A^\alpha \lambda(t)}\right\|_{\displaystyle L^\frac{p}{s+1-(2p-2)(s-s_p)} L^\frac{p}{2-s}}\\
&& \times \left\|\begin{array}{r} \displaystyle \int_0^1 \int_0^1 F''(\tilde{\tau} u_{\leq A^\alpha \lambda(t)} + u_{ A^\alpha \lambda(t)<
\cdotp\leq A^\beta \lambda(t)}\\
+ \tau u_{> A^\beta \lambda(t)}) d\tau d\tilde{\tau}\end{array} \right\|_{\displaystyle L^\frac{2p}{(p-2)(s+1-(2p-2)(s-s_p))}
L^\frac{2p}{(p-2)(2-s)}}\\
&+& \left\| u_{>A^\beta \lambda(t)} \right\|_{L^\frac{2p}{s+1-(2p-2)(s-s_p)} L^\frac{2p}{2-s}}\\
&& \times \left\|\int_0^1 F'(u_{ A^\alpha \lambda(t)<\cdotp\leq A^\beta \lambda(t)} + \tau u_{> A^\beta \lambda(t)})
d\tau\right \|_{L^\frac{2p}{(p-1)(s+1-(2p-2)(s-s_p))} L^\frac{2p}{(p-1)(2-s)}}\\
&\lesssim& (\lambda(t))^{p(s-s_p)} \left[\left(\frac{A^\alpha \lambda(t)}{A^\beta \lambda(t)}\right)^{\sigma(p)}
+ S(A^\beta) S^{p-1} (A^\alpha)\right]\\
&\lesssim& (\lambda(t))^{p(s-s_p)} \left(A^{-(\beta-\alpha)\sigma(p)} + S(A^\beta)S^{p-1}(A^\alpha)\right).
\end{eqnarray*}
The bilinear estimate is used here to estimate the term $u_{>A^\beta \lambda(t)} u_{\leq A^\alpha \lambda(t)}$. Collecting both terms and taking sup for
all $t \in I$, we obtain the first inequality.
\paragraph{Proof of the second inequality} To prove the inequality (\ref{recurrence2}) we first define $t_i$ for $i \geq 1$ given $t_0 \in I$.
%This definition will be used in the further argument.
\begin{equation}
 t_i = t_{i-1} + d \lambda^{-1} (t_{i-1}).
\label{boxdef}
\end{equation}
By the choice of $d$ (please see proposition \ref{local compactness 1}), all $t_i$'s are in the maximal lifespan $I$.
By the Strichartz estimate and the Duhamel formula, we have
\begin{eqnarray*}
 &&\|u_{> \lambda(t_0) A}\|_{Y_s ([t_0,t_1])}\\
 &=& \left\|\int_t^{\infty} \frac{\sin((\tau-t)\sqrt{-\Delta})}{\sqrt{-\Delta}} P_{> \lambda(t_0)A} F(u(\tau)) d\tau\right\|_{Y_s ([t_0,t_1])}\\
 &\leq& \left\|\int_t^{t_2} \frac{\sin((\tau-t)\sqrt{-\Delta})}{\sqrt{-\Delta}} P_{> \lambda(t_0)A} F(u(\tau)) d\tau\right\|_{Y_s ([t_0,t_1])}\\
&& + \liminf_{T \rightarrow \infty} \left\|\int_{t_2}^{T} \frac{\sin((\tau-t)\sqrt{-\Delta})}{\sqrt{-\Delta}} P_{> \lambda(t_0)A} F(u(\tau))
d\tau\right\|_{Y_s([t_0,t_1])}\\
 &\lesssim& \|P_{> \lambda(t_0) A} F(u)\|_{Z_s ([t_0,t_2] \times \Rm^3)}\\
&& + \liminf_{T \rightarrow \infty} \left\|\int_{t_2}^{T} \frac{\sin((\tau-t)\sqrt{-\Delta})}{\sqrt{-\Delta}}
P_{> \lambda(t_0)A} F(u(\tau))d\tau\right\|_{Y_s ([t_0,t_1])}\\
&=& I_1 + I_2.
\end{eqnarray*}
The first term can be dominated by
\begin{eqnarray*}
 I_1 &\lesssim& \|P_{> \lambda(t_0) A} F(u)\|_{Z_s ([t_0,t_1] \times \Rm^3)}
+ \|P_{> \lambda(t_0) A} F(u)\|_{Z_s ([t_1,t_2] \times \Rm^3)}\\
&\lesssim& (\lambda(t_0))^{s-s_p}N(A) + (\lambda(t_1))^{s-s_p} N\left(\frac{\lambda(t_0)}{\lambda(t_1)}A\right)\\
&\lesssim& (\lambda(t_0))^{s-s_p} N(A^{1-\eps_1}).
\end{eqnarray*}
for any small positive number $\eps_1$ and sufficiently large $A > A_0 (u,\eps_1)$, because $\lambda(t_0)$ and $\lambda(t_1)$ are comparable
to each other by the local compactness result (\ref{choiced}).\\
Now let us consider the term $I_2$. First of all, by (\ref{uni hs es}), we have
\begin{eqnarray*}
 && \left\|\int_{t_2}^{T} \frac{\sin((\tau-t)\sqrt{-\Delta})}{\sqrt{-\Delta}} P_{> \lambda(t_0)A} F(u(\tau))
 d\tau\right\|_{L^\infty L^2([t_0,t_1] \times \Rm^3)}\\
 &\lesssim& \frac{1}{(\lambda(t_0)A)^{s_p}} \left\|\int_{t_2}^{T} \frac{\sin((\tau-t)\sqrt{-\Delta})}{\sqrt{-\Delta}} F(u(\tau))
 d\tau\right\|_{L_{[t_0,t_1]}^\infty \dot{H}^{s_p}(\Rm^3)}\\
 &\lesssim& \frac{1}{(\lambda(t_0)A)^{s_p}}
\end{eqnarray*}
Using lemma \ref{L infinity estimate}, we also obtain
\begin{eqnarray*}
\lefteqn{\left\|\int_{t_2}^{T} \frac{\sin((\tau-t)\sqrt{-\Delta})}{\sqrt{-\Delta}}
P_{> \lambda(t_0)A} F(u(\tau)) d\tau\right\|_{L^\infty L^\infty([t_0,t_1]\times \Rm^3)}}\\
&\lesssim& \left\|\int_{t_2}^{T} \frac{\sin((\tau-t)\sqrt{-\Delta})}{\sqrt{-\Delta}} F(u(\tau))
d\tau\right\|_{L^\infty L^\infty ([t_0,t_1]\times \Rm^3)}\\
&\lesssim& (t_2 -t_1)^{-2/(p-1)}\\
&\lesssim& (\lambda(t_0))^{2/(p-1)}.
\end{eqnarray*}
Using the interpolation between $L^2$ and $L^\infty$, we have
\begin{eqnarray*}
 && \left\|P_{> \lambda(t_0)A} \int_{t_2}^T \frac{\sin((\tau-t)\sqrt{-\Delta})}{\sqrt{-\Delta}}
F(u(\tau))d\tau\right\|_{L^\infty L^{\frac{6}{3-2s}}([t_0,t_1]\times \Rm^3)}\\
&\leq& \|\cdotp\|_{L^\infty L^\infty ([t_0,t_1]\times \Rm^3)}^{2s/3} \|\cdotp\|_{L^\infty L^2([t_0,t_1]\times \Rm^3)}^{(3-2s)/{3}}\\
&\lesssim& [\lambda(t_0)^{2/(p-1)}]^{2s/{3}} [(\lambda(t_0)A)^{-s_p}]^{(3-2s)/{3}}\\
&=& (\lambda(t_0))^{s -s_p} A^{\frac{-s_p(3-2s)}{3}}.
\end{eqnarray*}
Next we will use the interpolation again to gain the estimate of $Y_s$ norm. There are two technical lemmas.
\begin{lemma} \label{interpolation} There exists a constant $\kappa = \kappa(p) \in (0,1)$ that depends only on $p$, so that for each $s\in [s_p, 1)$,
there exists an $s$-admissible pair $(q,r)$, with $q \neq \infty$ and
\[
 \frac{s + 1 - (2p-2)(s-s_p)}{2p} = \kappa \cdot 0 + (1-\kappa) \frac{1}{q};\;\;\; \frac{2-s}{2p} = \kappa \frac{3-2s}{6}
 + (1-\kappa) \frac{1}{r}.
\]
\end{lemma}
\paragraph{Proof} This is just a basic and boring computation. Please see the Appendix.
\begin{lemma} \label{LqLr}
Given any $s$-admissible pair $(q,r)$ with $q < \infty$, we have
\[
 \lim_{T \rightarrow \infty} \left\|\int_{t_2}^{T} \frac{\sin((\tau-t_0)\sqrt{-\Delta})}{\sqrt{-\Delta}}
 F(u(\tau))d\tau\right\|_{L^q L^r ([t_0,t_1]\times \Rm^3)}
 \leq C (\lambda(t_0))^{s-s_p}.
\]
The constant $C$ does not depend on $t_0$.
\end{lemma}
\paragraph{Proof} By lemma \ref{Duhamel in Hs}, we have
\[
 \lim_{T \rightarrow \infty} \int_{t_2}^{T} \frac{\sin((\tau-t)\sqrt{-\Delta})}{\sqrt{-\Delta}} F(u(\tau))d\tau = S(t-t_2) (u(t_2),\partial_t u(t_2))
\]
in the space $L^q L^r ([t_0,t_1] \times \Rm^3)$. Thus
\begin{eqnarray*}
&& \lim_{T \rightarrow \infty} \left\|\int_{t_2}^{T_i} \frac{\sin((\tau-t_0)\sqrt{-\Delta})}{\sqrt{-\Delta}}
(u(\tau))d\tau\right\|_{L^q L^r ([t_0,t_1])}\\
&=& \|S(t-t_2) (u(t_2),\partial_t u(t_2))\|_{L^q L^r ([t_0,t_1] \times \Rm^3)}\\
&\lesssim& \|(u(t_2),\partial_t u(t_2))\|_{\dot{H}^s \times \dot{H}^{s-1}}\\
&\lesssim& (\lambda(t_2))^{s-s_p}\\
&\lesssim& (\lambda(t_0))^{s-s_P}.
\end{eqnarray*}
Now let us apply the two lemmas
\begin{eqnarray*}
I_2 &=&\liminf_{T \rightarrow \infty} \left\|\int_{t_2}^{T} \frac{\sin((\tau-t)\sqrt{-\Delta})}{\sqrt{-\Delta}} P_{> \lambda(t_0)A} F(u(\tau))
d\tau\right\|_{Y_s([t_0,t_1])}\\
&\leq&\liminf_{T \rightarrow \infty} \left( \begin{array}{l} \left\|\int_{t_2}^{T} \frac{\sin((\tau-t)\sqrt{-\Delta})}{\sqrt{-\Delta}}
P_{> \lambda(t_0)A} F(u(\tau)) d\tau\right\|_{L^\infty L^{\frac{6}{3-2s}([t_0,t_1]\times \Rm^3)}}^{\kappa(p)}\\
 \times\left\|\int_{t_2}^{T} \frac{\sin((\tau-t)\sqrt{-\Delta})}{\sqrt{-\Delta}} P_{> \lambda(t_0)A} F(u(\tau))
d\tau\right\|_{L^q L^r([t_0,t_1]\times \Rm^3)}^{1-\kappa(p)}\end{array} \right)\\
&\lesssim& \left[(\lambda(t_0))^{s -s_p} A^{\frac{-s_p(3-2s)}{3}}\right]^{\kappa(p)} \times \lim_{T \rightarrow \infty}\left\|\int_{t_2}^{T}
\frac{\sin((\tau-t)\sqrt{-\Delta})}{\sqrt{-\Delta}} F(u(\tau)) d\tau\right\|_{L^q L^r}^{1-\kappa(p)}\\
&\lesssim& \left[(\lambda(t_0))^{s -s_p} A^{\frac{-s_p(3-2s)}{3}}\right]^{\kappa(p)} (\lambda(t_0))^{(s-s_p)(1-\kappa(p))}\\
&\lesssim& (\lambda(t_0))^{s -s_p} A^{\displaystyle \frac{-s_p \kappa(p) (3-2s)}{3}}\\
&\lesssim& (\lambda(t_0))^{s -s_p} A^{-\sigma_1(p)}
\end{eqnarray*}
Here $\sigma_1 (p) = \kappa(p)/6$. It depends only on $p$.\\
Combining $I_1$ and $I_2$ and then taking the sup for all $t\in I$ , we finish the proof of the second inequality.
\subsection{Decay of $S(A)$ and $N(A)$}
Plugging the first recurrence formula into the second one, we gain
\begin{eqnarray*}
 S(A) \lesssim S(A^{(1-\eps_1)\beta}) S^{p-1}(A^{(1-\eps_1)\alpha}) &+&  A^{-\sigma(p)(1-\eps_1)(\beta - \alpha)}\\
 &+& A^{-2(1-\eps_1)(1-\beta)} + A^{-\sigma_1(p)}.
\end{eqnarray*}
Choose $\alpha$, $\beta$ and $\eps_1$ so that
\begin{equation}
 (1-\eps_1)\beta = 2/3;\,  (1-\eps_1)\alpha = 1/3;\,  \eps_1=1/10000.
\end{equation}
Then we have
\[
 S(A) \lesssim S(A^{2/3})S^{p-1} (A^{1/3}) + A^{-\sigma_2(p)}.
\]
for sufficiently large $A$. Here the positive number $\sigma_2 (p)$ depends on $p$ only.
\[
 \sigma_2 = \min\{\sigma(p)/3, \sigma_1(p), 0.6\}.
\]
Applying lemma \ref{gain decay from recurrence}, we have $S(A) \lesssim A^{-\sigma_2 (p)}$ for sufficiently large $A$. Plugging this in the
first recurrence formula, we have $N(A) \lesssim A^{-\sigma_2 (p)}$ for large $A$. Observing that both $S(A)$ and $N(A)$ is uniformly bounded,
we know these two decay estimates are actually valid for each $A > 0$. Now let us choose
\[
 s_1 = \min\{1, s + \frac{99}{100}\sigma_2(p)\};
\]
and define (local contribution of the Duhamel Formula)
\begin{eqnarray*}
v_{t'}(t) &=& \int_{t'}^{t' +d \lambda(t')^{-1}} \frac{\sin((\tau-t)\sqrt{-\Delta})}{\sqrt{-\Delta}} F(u(\tau)) d\tau;\\
\partial_t v_{t'} (t) &=& -\int_{t'}^{t' +d \lambda(t')^{-1}} \cos((\tau-t)\sqrt{-\Delta}) F(u(\tau)) d\tau.
\end{eqnarray*}
We obtain for any $t \leq t'$ and integer $k \geq 0$
\begin{eqnarray*}
&& \left\|P_{ \lambda(t') 2^k < \cdot < \lambda(t')2^{k+1}}(v_{t'}(t), \partial_t v_{t'}(t)) \right\|_{\dot{H}^{s_1} \times \dot{H}^{s_1 -1}}\\
&\lesssim& (\lambda(t') 2^k)^{s_1 -s} \left\|P_{ \lambda(t') 2^k < \cdot < \lambda(t') 2^{k+1}}(v_{t'}(t), \partial_t v_{t'}(t))
\right\|_{\dot{H}^{s} \times \dot{H}^{s-1}}\\
&\lesssim& (\lambda(t') 2^k)^{s_1 -s} \left\|P_{> \lambda(t') 2^k}(v_{t'}(t), \partial_t v_{t'}(t)) \right\|_{\dot{H}^{s} \times \dot{H}^{s-1}}\\
&\lesssim& (\lambda(t') 2^k)^{s_1 -s} \left\|P_{> \lambda(t') 2^k} F(u)\right\|_{Z_s([t', t'+d \lambda(t')^{-1}])}\\
&\lesssim& (\lambda(t') 2^k)^{s_1 -s} (\lambda(t'))^{s-s_p} N(2^k)\\
&\lesssim& (\lambda(t'))^{s_1 -s_p} (2^k)^{s_1 -s -\sigma_2(p)}.
\end{eqnarray*}
Summing for all $k \geq 0$, we have
\[
\left\|P_{>\lambda(t')}(v_{t'}(t), \partial_t v_{t'}(t)) \right\|_{\dot{H}^{s_1} \times \dot{H}^{s_1 -1}} \lesssim (\lambda(t'))^{s_1 -s_p}.
\]
Combining this with
\begin{eqnarray*}
&&\left\|P_{\leq\lambda(t')}(v_{t'}(t), \partial_t v_{t'}(t)) \right\|_{\dot{H}^{s_1} \times \dot{H}^{s_1 -1}}\\
&\lesssim& (\lambda(t'))^{s_1 -s_p} \left\|P_{\leq\lambda(t')}(v_{t'}(t), \partial_t v_{t'}(t)) \right\|_{\dot{H}^{s_p} \times \dot{H}^{s_p -1}}\\
&\lesssim& (\lambda(t'))^{s_1 -s_p},
\end{eqnarray*}
we obtain
\begin{equation}
\left\|(v_{t'}(t), \partial_t v_{t'}(t))\right\|_{\dot{H}^{s_1} \times \dot{H}^{s_1 -1}}\lesssim (\lambda(t'))^{s_1 -s_p}.
 \label{local contribution in Duhamel}
\end{equation}
\subsection{Higher Regularity}
In this section we will show $(u(x,t),\partial_t u(x,t)) \in \dot{H}^{s_1} \times \dot{H}^{s_1 -1}(\Rm^3)$ for each $t \in I$.
\paragraph{Center estimate} Let us break the Duhamel formula into two pieces.
\begin{eqnarray*}
 u^{(1)} (t) &=& \int_{t}^{t_1} \frac{\sin((\tau-t)\sqrt{-\Delta})}{\sqrt{-\Delta}} F(u(\tau)) d\tau;\\
 u^{(2)} (t) &=& \int_{t_1}^{\infty} \frac{\sin((\tau-t)\sqrt{-\Delta})}{\sqrt{-\Delta}} F(u(\tau)) d\tau.
\end{eqnarray*}
%and consider the cutoff version of $u_3 (t)$
%\begin{eqnarray*}
% \tilde{u}_{3,T} (t) &=& \int_{t_1}^{T} \frac{\sin((\tau-t)\sqrt{-\Delta})}{\sqrt{-\Delta}} \chi(|x|>\frac{d\lambda(t_0)^{-1}}{2}+|\tau-t_1|)
% F(u(\tau)) d\tau;\\
% \tilde{u}_3 (t)     &=& \lim_{T \rightarrow \infty} \tilde{u}_{3,T} (t).
%\end{eqnarray*}
Let $\chi$ be the characteristic function of the region $\{(x,t): |x|>\frac{d\lambda^{-1}(t_0)}{2}+|t-t_1|\}$. We have
\begin{eqnarray*}
&& \left\|\chi F(u(t))\right\|_{L^1 L^\frac{6}{5-2s_1}([t_1,\infty)\times \Rm^3)}\\
&=& \int_{t_1}^\infty \left(\int_{|x|>\frac{d\lambda^{-1}(t_0)}{2}+|t-t_1|} (F(u))^\frac{6}{5-2s_1} dx \right)^{\frac{5-2s_1}{6}} dt\\
&\lesssim& \int_{t_1}^\infty \left(\int_{|x|>\frac{d\lambda^{-1}(t_0)}{2}+|t-t_1|} \left(\frac{1}{|x|^\frac{2p}{p-1}}\right)^{\frac{6}{5-2s_1}}
dx \right)^{\frac{5-2s_1}{6}} dt\\
&\lesssim& \int_{t_1}^\infty \left(\frac{1}{\left|\frac{d \lambda^{-1}(t_0)}{2} + t -t_1\right|^{\frac{2p}{p-1}\frac{6}{5-2s_1} -3}}
\right)^{\frac{5-2s_1}{6}} dt\\
&\lesssim& \int_{t_1}^\infty \left(\frac{d \lambda^{-1}(t_0)}{2} + t -t_1\right)^{s_p -s_1 -1} dt\\
&\lesssim& \lambda(t_0)^{s_1 -s_p}.
\end{eqnarray*}
%This implies $(\tilde{u}_{3,T}(t_0), \partial_t \tilde{u}_{3,T}(t_0))$ converges to $(\tilde{u}_3(t_0), \partial_t \tilde{u}_3(t_0))$ strongly in
%$\dot{H}^{s_1} \times \dot{H}^{s_1 -1}(\Rm^3)$ with
By lemma \ref{Strong Huygens in Duhamel}, there exists a pair $(\tilde{u}_0,\tilde{u}_1)$ so that
%\[
% \|(\tilde{u}_3(t_0), \partial_t \tilde{u}_3(t_0))\|_{\dot{H}^{s_1} \times \dot{H}^{s_1 -1}(\Rm^3)} \lesssim \lambda(t_0)^{s_1 -s_p}.
%\]
\[
 \|(\tilde{u}_0,\tilde{u}_1)\|_{\dot{H}^{s_1} \times \dot{H}^{s_1 -1}(\Rm^3)} \lesssim \lambda(t_0)^{s_1 -s_p},
\]
%In addition, by strong Huygens' principal, the pair $(\tilde{u}_{3,T}(t_0), \partial_t \tilde{u}_{3,T}(t_0))$
%is the same as $(u_{3,T}(t_0), \partial_t u_{3,T}(t_0))$ in the ball $B(0,d\lambda(t_0)^{-1}/2)$. Thus
%the pair $(\tilde{u}_{3,T}(t_0), \partial_t \tilde{u}_{3,T}(t_0))$ converges to $(u_3(t_0), \partial_t u_3(t_0))$ weakly in the ball
%$B(0,d\lambda(t_0)^{-1}/2)$ by the weak convergence of $(u_{3,T}(t_0), \partial_t u_{3,T}(t_0))$. Considering both the strong and weak convergence,
%we have
\[
 (u^{(2)}(t_0), \partial_t u^{(2)}(t_0)) = (\tilde{u}_0, \tilde{u}_1)\; \hbox{in}\; B\left(0, \frac{d \lambda^{-1}(t_0)}{2}\right).
\]
This implies
\begin{equation}
 (u(t_0), \partial_t u(t_0)) = (\tilde{u}_0 + u^{(1)}(t_0), \tilde{u}_1+ \partial_t u^{(1)}(t_0))\;
 \hbox{in}\; B\left(0, \frac{d \lambda^{-1}(t_0)}{2}\right).
 \label{center identity}
\end{equation}
By (\ref{local contribution in Duhamel}), we have
\[
 \|(u^{(1)}(t_0),\partial_t u^{(1)}(t_0))\|_{\dot{H}^{s_1} \times \dot{H}^{s_1 -1}} \lesssim \lambda(t_0)^{s_1 -s_p}.
\]
Combining this with the $\dot{H}^{s_1} \times \dot{H}^{s_1-1}$ bound of $(\tilde{u}_0, \tilde{u}_1)$, we have
\begin{equation}
\|(\tilde{u}_0 + u^{(1)}(t_0), \tilde{u}_1 + \partial_t u^{(1)}(t_0))\|_{\dot{H}^{s_1} \times \dot{H}^{s_1-1}}
\lesssim \lambda(t_0)^{s_1 -s_p}. \label{Hs1 estimate center}
\end{equation}
\paragraph{Tail Estimate}
Let $(u'_0, u'_1) = \Psi_{d\lambda^{-1}(t_0)/4} (u(t_0),\partial_t u(t_0))$, and
\[
 \frac{1}{q} = \frac{1}{2} + \frac{1-s_1}{3}.
\]
By theorem \ref{energy estimate near infinity}, if $r \geq d \lambda^{-1}(t_0)/4$, we have
\begin{eqnarray*}
 \left(\int_{r<|x|<4r} (|\nabla u'_0|^q + |u'_1|^q) dx\right)^{1/q} &\lesssim& \left(\int_{r<|x|<4r} (|\nabla u'_0|^2 + |u'_1|^2) dx\right)^{1/2}
 (r^3)^{\frac{1}{q} - \frac{1}{2}}\\
 &\lesssim& r^{-(1-s_p)} (r^3)^{(1-s_1)/3}\\
 &\lesssim& r^{-(s_1-s_p)}.
\end{eqnarray*}
Letting $r = 2^k d \lambda^{-1}(t_0)/4$ and summing for all $k \geq 0$,
we obtain that the pair $(u'_0,u'_1)$ is in the space $\dot{W}^{1,q} \times L^q(\Rm^3)$ with
\[
 \|(u'_0,u'_1)\|_{\dot{W}^{1,q} \times L^q(\Rm^3)} \lesssim (d \lambda(t_0)^{-1}/4)^{-(s_1-s_p)} \lesssim (\lambda(t_0))^{s_1 -s_p}.
\]
By the Sobolev embedding, we have
\begin{equation}
 \|(u'_0,u'_1)\|_{\dot{H}^{s_1} \times \dot{H}^{s_1-1}(\Rm^3)} \lesssim (\lambda(t_0))^{s_1 -s_p}. \label{lin11}
\end{equation}
Considering (\ref{center identity}), (\ref{Hs1 estimate center}), (\ref{lin11}) and using lemma \ref{glue}, we have
\[
 \|(u(t_0), \partial_t u(t_0))\|_{\dot{H}^{s_1}\times\dot{H}^{s_1-1}(\Rm^3)} \lesssim (\lambda(t_0))^{s_1 -s_p}.
\]
\subsection{Conclusion}
There are two cases
\paragraph{Case 1 $(s_1 =1)$} Now we have finished our argument and obtain the energy estimate.
\paragraph{Case 2 $(s_1 <1)$} This means $s_1 = s + 0.99 \sigma_2(p)$. Now let us consider the set
\[
 \left\{\left(\frac{1}{\lambda(t)^{3/2-s_p}} u\left(\frac{x}{\lambda(t)}, t\right),
 \frac{1}{\lambda(t)^{5/2-s_p}} \partial_t u\left(\frac{x}{\lambda(t)}, t\right)\right): t \in I\right\}
\]
This is precompact in the space $\dot{H}^{s_p} \times \dot{H}^{s_p -1}$, and bounded in the space
$\dot{H}^{s + 0.99\sigma_2 (p)} \times \dot{H}^{s -1 + 0.99 \sigma_2 (p)}$, thus it is also precompact
in the space $\dot{H}^{s + 0.98\sigma_2 (p)} \times \dot{H}^{s -1 + 0.98 \sigma_2 (p)}$.

\section{Global Energy Estimate and its Corollary} Repeat the recurrence process we described in the previous section starting
from the space $\dot{H}^{s_p} \times \dot{H}^{s_p -1}$. Each time we either gain the global energy estimate below or
gain additional regularity by $0.98\sigma_2 (p)$. This number depends on $p$ only. As a result, the
process has to stop at $\dot{H}^1 \times L^2$ after finite steps.
\begin{proposition}\emph{\textbf{Global Energy Estimate}}
Let $u(x,t)$ be a minimal blow-up solution. Then $(u(t_0), \partial_t u(t_0))$ is in the energy space for each $t_0 \in I$ with
\begin{equation}
 \|(u(t_0), \partial_t u(t_0))\|_{\dot{H}^1 \times L^2 (\Rm^3)} \lesssim \lambda(t_0)^{1-s_p}.
\end{equation}
By the local theory, we actually obtain
\[
  (u(t), \partial_t u(t)) \in C(I; \dot{H}^1 (\Rm^3) \times L^2 (\Rm^3)).
\]
\end{proposition}
\paragraph{Remark} By lemma \ref{eqn between w and u}, we have the following holds for any $0<a<b<\infty$
\[
 (\partial_r w(t), \partial_t w(t)) \in C(I; L^2 \times L^2 ([a,b]) )
\]
\subsection{Self-similar and High-to-low Frequency Cascade Cases}
In both two cases, we can choose $t_i \rightarrow \infty$ such that $\lambda(t_i) \rightarrow 0$. This implies
\[
 \int_{\Rm^3} (|\nabla u (t_i)|^2 + |\partial_t u(t_i)|^2) dx \rightarrow 0.
\]
By the Sobolev embedding, we have
\begin{equation}
\|u\|_{L^{p+1}(\Rm^3)}^{p+1} \leq \|u\|_{L^{\frac{3}{2}(p-1)}(\Rm^3)}^{p-1} \|u\|_{L^{6}(\Rm^3)}^{2}
\lesssim \|u\|_{\dot{H}^{s_p}(\Rm^3)}^{p-1} \|u\|_{\dot{H}^1 (\Rm^3)}^{2}. \label{finite energy}
\end{equation}
This implies $\|u(t_i)\|_{L^{p+1}(\Rm^3)}^{p+1} \rightarrow 0$. Using the definition of energy we have $E(t_i) \rightarrow 0$.
On the other hand, we know the energy is a constant. Therefore the energy must be zero.
\paragraph{Defocusing Case} It is nothing to say, because in this case an energy zero means that the solution is identically zero.
\paragraph{Focusing Case} We can still solve the problem using the following theorem. By the fact that the energy is zero,
 we know $u$ blows up in finite time in both time directions. But this is a contradiction with our assumption $T_+ = \infty$.
\begin{theorem}\label{Non-positive energy implies blow-up}
(Please see theorem 3.1 in \cite{kv1}, Non-positive energy implies blowup)\\
Let $(u_0,u_1) \in (\dot{H}^1 \times L^2) \cap (\dot{H}^{s_p} \times \dot{H}^{s_p -1})$ be initial data. Assume that
$(u_0,u_1)$ is not identically zero and satisfies $E(u_0, u_1) \leq 0$. Then the maximal life-span solution to the non-linear wave equation
blows up both forward and backward in finite time.
\end{theorem}
\subsection{Soliton-like Solutions in the Defocuing Case}
Now let us consider the soliton-like solutions in the defocusing case. First we have a useful global integral estimate in the defocusing case.
\begin{lemma} (Please see \cite{benoit}) Let $u$ be a solution of (\ref{eqn}) defined in a time interval $[0,T]$
with $(u,\partial_t u) \in \dot{H}^1 \times L^2$ and a finite energy
\[
 E = \int_{\Rm^3} \left(\frac{1}{2}|\nabla_x u|^2 + \frac{1}{2}|\partial_t u|^2 + \frac{1}{p+1}|u(x)|^{p+1}\right) dx.
\]
For any $R>0$, we have
\begin{eqnarray*}
 \frac{1}{2R} \int_0^T \int_{|x|< R} (|\nabla u |^2 + |\partial_t u|^2) dx dt + \frac{1}{2R^2} \int_0^T \int_{|x|=R} |u|^2 d\sigma_R dt && \\
+ \frac{1}{2R} \frac{2p - 4}{p+1} \int_0^T \int_{|x|<R} |u|^{p+1} dx dt
+ \frac{p-1}{p+1} \int_0^T \int_{|x|>R} \frac{|u|^{p+1}}{|x|} dx dt &&\\ + \frac{2}{R^2} \int_{|x|<R} |u(T)|^2 dx
 &\leq& 2E.
\end{eqnarray*}
\end{lemma}
\noindent Observing that each term on the left hand is nonnegative, we can obtain a uniform upper bound for the last term in the second line above
\[
  \int_0^T \int_{|x|>R} \frac{|u|^{p+1}}{|x|} dx dt \leq \frac{2(p+1)}{p-1} E.
\]
Letting $R$ approach zero and $T$ approach $T_{+}$, we have
\begin{equation}
\int_0^{T_{+}} \int_{\Rm^3} \frac{|u|^{p+1}}{|x|} dx dt \leq \frac{2(p+1)}{p-1} E.
\label{upperbound}
\end{equation}
The energy $E$ here is finite by our estimate (\ref{finite energy}).
On the other hand, recalling our local compactness result(See lemma \ref{lower bound of int u p+1 over x}),
we obtain ($T_+ = \infty$)
\[
 \int_0^\infty \int_{\Rm^3} \frac{|u|^{p+1}}{|x|} dx dt = \infty.
\]
This finishes our discussion in this case.
\section{Further Estimates in the Soliton-like Case}
Let $u$ be a soliton-like minimal blow-up solution. We will find additional decay of $u$ at infinity. The method used here is similar to
the one C.E.Kenig and F.Merle used in their paper \cite{km} for super-critical case.
Throughout this section $w(r,t)$, $h(r,t)$, $z_1 (r,t)$ and $z_2(r,t)$ are defined as usual using $u(x,t)$.
\paragraph{Remark} The argument in this section works in both the defocusing and focusing case. But we are particularly interested in the focusing case,
because the soliton-like solutions in the focusing case are the only solutions that still survive.

\subsection{Setup} Let $\varphi(x)$ be a smooth cutoff function in $\Rm^3$.
\[
 \varphi(x) \left\{\begin{array}{l} =0, \,\,\,\, |x| \leq 1/2;\\
 \in [0,1],\,\,\,\, 1/2 \leq |x| \leq 1;\\
 = 1,\,\,\,\, |x| \geq 1.\end{array}\right.
\]
Then by the compactness of $u$(Please see lemma \ref{local compactness 2}), $\|\varphi(x/R) u(x,t)\|_{\dot{H}^{s_p}}$ converges to zero uniformly in $t$ as
$R \rightarrow \infty$. Thus we have a positive function $g(r)$ so that $g(r)$ decreases to zero as $r$ increases to infinity with
\[
 \|\varphi(x/R) u(x,t)\|_{\dot{H}^{s_p}} \leq g(R).
\]
This means for each $|x|\geq R$
\[
 |u(x,t)| = |\varphi(x/R) u(x,t)| \leq C \frac{\|\varphi(\cdot/R) u(\cdot,t)\|_{\dot{H}^{s_p}}}{|x|^{2/(p-1)}} \leq \frac{C g(R)}{|x|^{2/(p-1)}}.
\]
Let us define
\[
 f_\beta (r) = \sup_{t\in \Rm, |x|\geq r} |x|^\beta |u(x,t)|
\]
for $\beta \in [2/(p-1), 1)$ and $r>0$. This is a nonincreasing function of $r$ defined from $\Rm^+$ to $[0,\infty)\cup \{\infty\}$.
Consider the set
\[
 U =\{\beta \in [2/(p-1),1): f_\beta(r) \rightarrow 0\; \hbox{as}\; r \rightarrow \infty\}.
\]
This is not empty, since $2/(p-1)$ is in $U$. Due to the estimate
\[
 |x|^\beta |u(x,t)| \leq C_p |x|^{\beta - \frac{2}{p-1}} \|u(\cdot,t)\|_{\dot{H}^{s_p}},
\]
we know if $\beta \in U$, then $f_\beta(r)$ is a bounded function. By definition of $f_\beta$, we have for any time $t$ and $|x|\geq r$
\begin{equation}
  |u(x,t)| \leq \frac{f_\beta(r)}{|x|^{\beta}}. \label{def of f beta}
\end{equation}
This is a meaningful inequality as long as $\beta \in U$.
%Define
%\[
% f(r) = \sup_{|x| \geq r} |x|^{\frac{3}{2} -s_p} |u(x,t)| \leq C g(r).
%\]
%Thus $f(r)$ decreases to zero as $r \rightarrow \infty$ with
%\begin{equation}
%  |u(x,t)| \leq \frac{f(r)}{|x|^{(3/2)-s_p}}. \label{def of f}
%\end{equation}
%for all $|x|\geq r$. It is easy to check $f(r)$ is bounded.
\paragraph{Local Energy of $w$} Let $\beta \in U$. Applying lemma \ref{derivative of int w} to $w$ we have
\begin{eqnarray*}
\left(\int_{r_0}^{4r_0} |z_1(r,t_0)|^2 dr\right)^{1/2}
&\leq& \left(\int_{r_0 + M}^{4r_0 + M} |z_1(r,t_0+M)|^2 dr\right)^{1/2}\\
&& + \left(\int_{r_0}^{4 r_0} \left(\int_0^M h(r +t, t_0 +t) dt\right)^2 dr \right)^{1/2}.
\end{eqnarray*}
Let $M \rightarrow \infty$. Using proposition \ref{int w z estimate} we have
\begin{eqnarray*}
\left(\int_{r_0}^{4r_0} |z_1(r,t_0)|^2 dr\right)^{1/2}
&\leq& \limsup_{M \rightarrow \infty} \left(\int_{r_0}^{4 r_0} \left(\int_0^M (r+t)F(u(r +t, t_0 +t)) dt\right)^2 dr \right)^{1/2}\\
&\leq& \limsup_{M \rightarrow \infty} \left(\int_{r_0}^{4 r_0}
\left(\int_0^M (r+t)\left(\frac{f_\beta (r_0)}{(r+t)^{\beta}}\right)^p dt\right)^2 dr \right)^{1/2}\\
&\lesssim_p& \limsup_{M \rightarrow \infty} \left(\int_{r_0}^{4 r_0} \left(\frac{f_\beta^p (r_0)}{r^{p\beta -2}}\right)^2 dr \right)^{1/2}\\
&\leq& f_\beta^p (r_0) \left(\int_{r_0}^{4 r_0} \frac{1}{r^{2p\beta -4}} dr \right)^{1/2}\\
&\lesssim_p& f_\beta^p (r_0) \left( \frac{1}{r_0^{2p\beta -5}} \right)^{1/2}\\
&\leq& f_\beta^p (r_0) \frac{1}{r_0^{p\beta -5/2}}.
\end{eqnarray*}
Similarly we have
\[
\left(\int_{r_0}^{4r_0} |z_2(r,t_0)|^2 dr\right)^{1/2} \lesssim  \frac{f_\beta^p (r_0)}{r_0^{p\beta -5/2}}.
\]
In summary we obtain
\begin{equation}
\left(\int_{r_0}^{4r_0} |\partial_t w(r,t_0)|^2 + |\partial_r w(r,t_0)|^2 dr\right)^{1/2} \leq  C_p \frac{f_\beta^p (r_0)}{r_0^{p\beta -5/2}}.
\label{decay of energy with f}
\end{equation}
The constant $C_p$ depends on $p$ only.
\paragraph{Remark} The estimate above holds as long as $\beta \geq 2/(p-1)$ and the inequality
\[
 |u(x,t)| \leq \frac{f(r)}{|x|^\beta}
\]
holds for all $|x| \geq r> 0$.
\subsection{Recurrence Formula}
We know $w = ru$ is a solution to the one-dimensional wave equation
\[
 \partial_t^2 w - \partial_r^2 w = r |u|^{p-1}u.
\]
Using the explicit formula to solve this equation, we obtain
\begin{eqnarray*}
 r_0 u(r_0, t_0) &=& \frac{1}{2} \left[\left(r_0 + \frac{r_0}{2}\right) u\left(r_0 + \frac{r_0}{2}, t_0 - \frac{r_0}{2}\right)
 + \left(r_0 - \frac{r_0}{2}\right) u\left(r_0 - \frac{r_0}{2}, t_0 - \frac{r_0}{2}\right)\right]\\
 && + \frac{1}{2} \int_{r_0 - \frac{r_0}{2}}^{r_0 + \frac{r_0}{2}} \partial_t w\left(r, t_0 - \frac{r_0}{2}\right) dr\\
 && + \frac{1}{2} \int_0^\frac{r_0}{2} \int_{\frac{r_0}{2} + t}^{\frac{3 r_0}{2} -t} r |u|^{p-1} u\left(r,t_0 - \frac{r_0}{2} + t\right) dr dt\\
 &=& I_1 + I_2 + I_3.
\end{eqnarray*}
By Cauchy-Schwartz and (\ref{decay of energy with f}), we have
\begin{eqnarray*}
|I_2| &\leq& \frac{1}{2} \left(\int_{\frac{r_0}{2}}^{\frac{3r_0}{2}} \left|\partial_t w\left(r, t_0 -
\frac{r_0}{2}\right)\right|^2 dr\right)^{1/2} \left(\int_{\frac{r_0}{2}}^{\frac{3r_0}{2}} 1 dr\right)^{1/2}\\
&\leq& C_p \frac{f_\beta^p (\frac{r_0}{2})}{r_0^{p\beta -5/2}} r_0^{1/2}\\
&=& C_p f_\beta^p (\frac{r_0}{2}) r_0^{3-p\beta}.
\end{eqnarray*}
Next we estimate $I_3$ using (\ref{def of f beta})
\begin{eqnarray*}
 |I_3| &\leq& \frac{1}{2} \int_0^{r_0/2} \int_{r_0/2 + t}^{3 r_0/2 -t} r \left(\frac{f_\beta(r_0/2)}{r^{\beta}}\right)^p dr dt\\
 &\leq& C_p \int_0^{r_0/2} r_0^2 \frac{f_\beta^p(r_0/2)}{r_0^{p\beta}} dt\\
 &\leq& C_p f_\beta^p(\frac{r_0}{2}) r_0^{3- p\beta}.
\end{eqnarray*}
While
\begin{eqnarray*}
|I_1| &\leq& \frac{1}{2}\left[\frac{3r_0}{2} \frac{f_\beta (3r_0/2)}{(3r_0/2)^{\beta}} + \frac{r_0}{2} \frac{f_\beta(r_0/2)}{(r_0/2)^{\beta}}\right]\\
&=& \frac{1}{2}\left[\left(\frac{3}{2}\right)^{1-\beta} f_\beta(\frac{3r_0}{2}) +
\left(\frac{1}{2}\right)^{1-\beta} f_\beta (\frac{r_0}{2})\right] r_0^{1 - \beta}.
\end{eqnarray*}
Combining these three terms and dividing both sides of the inequality by $r_0^{1 -\beta}$, we obtain (replace $r_0$ by $r$)
\[
r^{\beta} |u(r,t_0)| \leq \frac{1}{2}\left[\left(\frac{3}{2}\right)^{1-\beta} f_\beta (\frac{3r}{2}) +
\left(\frac{1}{2}\right)^{1-\beta} f_\beta (\frac{r}{2})\right] + C_p f_\beta^p (\frac{r}{2}) r^{2-(p-1)\beta}.
\]
Observing that the right hand is a nonincreasing function of $r$, we apply $\sup_{r\geq r_0}$ on both sides and obtain
\begin{equation}
 f_\beta (r_0) \leq \frac{1}{2}\left[\left(\frac{3}{2}\right)^{1-\beta} f_\beta (\frac{3r_0}{2}) + \left(\frac{1}{2}\right)^{1-\beta} f_\beta
 (\frac{r_0}{2})\right] + C_p f_\beta^p (\frac{r_0}{2}) r_0^{2-(p-1)\beta}. \label{recurrence}
\end{equation}
Thus
\begin{equation}
 f_\beta (r_0) \leq \frac{1}{2}\left[\left(\frac{3}{2}\right)^{1-\beta} + \left(\frac{1}{2}\right)^{1-\beta} \right] f_\beta (\frac{r_0}{2})
 + C_p f_\beta^p (\frac{r_0}{2}) r_0^{2-(p-1)\beta}. \label{f beta recurrence}
\end{equation}
\subsection{Decay of $u(x,t)$}
Let
\[
 g(\beta) = \frac{1}{2}\left[\left(\frac{3}{2}\right)^{1-\beta} + \left(\frac{1}{2}\right)^{1-\beta} \right] < 1.
\]
Because $f_\beta(r) \rightarrow 0$ and $2-(p-1)\beta \leq 0$, we know that there exists a large constant $R$,
such that if $r_0 > R$,
\[
 C_p f_\beta^p (\frac{r_0}{2}) r_0^{2-(p-1)\beta} \leq \frac{1 -g(\beta)}{2} f_\beta (\frac{r_0}{2}).
\]
Thus we have for $r_0 > R$,
\[
 f_\beta (r_0) \leq \frac{g(\beta)+1}{2} f_\beta (r_0 /2).
\]
This implies
\[
 f_\beta (r) \leq C r^{\displaystyle \log_2 (\frac{g(\beta)+1}{2})}
\]
for sufficiently large $r > R'$. As a result, for each $\beta_1 < \beta - \log_2 (\frac{g(\beta)+1}{2})$, we have (Note that the logarithm is negative)
\[
 |x|^{\beta_1} |u(x,t)| \leq f_\beta(|x|)|x|^{\beta_1 -\beta} \leq C |x|^{\beta_1 -\beta + \log_2 (\frac{g(\beta)+1}{2})} \rightarrow 0
\]
as $|x| \rightarrow \infty$. This implies
\[
 \left[\beta, \beta + \log_2 \frac{2}{1 + g(\beta)} \right) \subseteq U.
\]
\paragraph{The Upper Bound of $U$} Now we are ready to show $\sup U = 1$, if this was false, we could assume $\sup U = \beta_0 < 1$. Then
we have for each $\beta \in U$,
\[
 g(\beta) \leq G_0 \doteq \max \left\{g(\beta_0),g\left(\frac{2}{p-1}\right)\right\} < 1
\]
using the convexity of the function $g$. Thus
\[
\log_2 \frac{2}{1 + g(\beta)} \geq \log_2 \frac{2}{1 + G_0} > 0.
\]
This means
\[
 \left[\beta, \beta + \log_2 \frac{2}{1 + G_0} \right) \subseteq U.
\]
This gives us a contradiction as $\beta \rightarrow \sup U$.
\paragraph{Decay of $u$} Let $\beta$ be a number slightly smaller than 1. We know $\beta \in U$.
By (\ref{decay of energy with f}), we have
\begin{eqnarray*}
\int_{r_0}^{4 r_0} |\partial_r w(r,t_0)| dr &\leq& \left(\int_{r_0}^{4 r_0} |\partial_r w(r,t_0)|^2 dr\right)^{1/2}
\left(\int_{r_0}^{4 r_0} 1 dr\right)^{1/2}\\
&\leq& \frac{C_p f_\beta^p (r_0)}{r_0^{p\beta -5/2}} r_0^{1/2}\\
&\leq& \frac{C_{p,\beta}}{r_0^{p\beta -3}}
\end{eqnarray*}
We can choose $\beta \in U$ so that $p\beta - 3 > 0$ by the fact $p>3$. Thus we have
\begin{equation}
 \int_{1}^{\infty} |\partial_r w(r,t_0)| dr \leq C_{p,\beta}. \label{L1 for partial w}
\end{equation}
In addition for $r \leq 1$,
\[
 |w(r,t_0)| = r |u (r,t_0)| \leq C \|u(t_0)\|_{\dot{H}^{s_p}} r^{1 - \frac{2}{p-1}} \leq C \|u(t_0)\|_{\dot{H}^{s_p}}.
\]
Combining the two estimates above, we know that $|w(r,t)|$ is bounded by a universal constant $C_1$ for each pair $(r,t)$. Thus
\begin{equation}
 |u(x,t)| \leq \frac{C_1}{|x|}. \label{decay of u}
\end{equation}
Plugging this in the definition of $f_\beta (r)$, we have
\[
 f_\beta (r_0) = \sup_{|x| \geq r_0} |x|^\beta |u(x,t)| \leq \sup_{|x| \geq r_0} C_1 |x|^{\beta -1} = C_1 r_0^{\beta -1}.
\]
Plugging this in (\ref{decay of energy with f}), we obtain
\begin{equation}
\left(\int_{r_0}^{4r_0} |\partial_t w(r,t_0)|^2 + |\partial_r w(r,t_0)|^2 dr\right)^{1/2} \lesssim \frac{1}{r_0^{p -5/2}}.
\label{strong decay of w energy}
\end{equation}
By lemma \ref{eqn between w and u}, the estimates (\ref{decay of u}) and (\ref{strong decay of w energy}) imply
\begin{equation}
 \int_{r <|x|<4r} (|\nabla u(x,t)|^2 + |\partial_t u(x,t)|^2 ) dx \lesssim r^{-1}.
\end{equation}
\section{Death of Soliton-like Solution}
\subsection{Solitons in the Focusing Case}
In order to kill the soliton-like minimal blow-up solutions, we need to consider the solitons of the wave equation. It turns out that
there does not exist any soliton for our equation. The elliptic equation
\begin{equation}
 -\Delta W(x) = |W(x)|^{p-1} W(x) \label{elliptic equation}
\end{equation}
does admit a lot of radial solutions. However, none of these solutions is in the space $\dot{H}^{s_p}$. Among these solutions we are particularly interested
in the solutions satisfying the same kind of property at infinity as (\ref{decay of u}).
\begin{proposition}
The elliptic equation (\ref{elliptic equation}) has a solution $W_0(x)$ so that
\begin{itemize}
\item $W_0(x)$ is a radial and smooth solution in $\Rm^3 \setminus \{0\}$.
\item The point $0$ is a singularity of $W_0(x)$.
\item The solution $W_0(x)$ is NOT in the space $\dot{H}^{s_p}(\Rm^3)$.
\item Its behavior near infinity is given by ($|x| > R_0$)
\begin{equation}
 \left|W_0 (x) - \frac{1}{|x|}\right| \leq \frac{C}{|x|^{p-2}};\;\;\; |\nabla W_0 (x)| \leq \frac{C}{|x|^2}. \label{behaviour of w at infinity}
\end{equation}
\end{itemize}
\end{proposition}
\noindent Please see the last section for a complete discussion of this solution.

\paragraph{Idea to deal with the soliton-like solutions} We will show there does not exist a soliton-like minimal blow-up solution in the
focusing case. This conclusion is natural because there is actually no soliton. However to prove this result is not an easy task.
We will use a method developed by T. Duyckaerts, C. E. Kenig and F. Merle as I mentioned at the beginning of this paper.
In their paper \cite{secret} they use this method to prove the soliton resolution conjecture for
radial solutions of the focusing, energy-critical wave equation. The idea is to show that our soliton-like solution has to be so close to
the solitons $\pm W_0 (x)$ or their rescaled versions that they must be the same. But the soliton we mentioned above is not in the right space.
This is a contradiction. In order to achieve this goal, we have to be able to understand the behaviour of a minimal blow-up solution
if it is close to our soliton $W_0 (x)$.

\subsection{Preliminary Results}
First of all, we recall a lemma proved in \cite{tkm1}.
\begin{lemma} \emph{\textbf{(Energy channel)}} \label{energy channel}
Let $(v_0 ,v_1) \in \dot{H}^1 \times L^2$ be a pair of radial initial data. Suppose $v(x,t)$ is the solution of the linear
wave equation with the given initial data $(v_0,v_1)$. Let $w(r,t) = r v(r,t)$ as usual, then for any $R > 0$ either the inequality
\[
 \int_{|x|> R + t} (|\nabla v(x,t)|^2 + |\partial_t v(x,t)|^2) dx \geq 2\pi \int_R^\infty (|\partial_r w(r,0)|^2 + |\partial_t w(r,0)|^2) dr
\]
holds for all $t > 0$, or
\[
 \int_{|x|> R - t} (|\nabla v(x,t)|^2 + |\partial_t v(x,t)|^2) dx \geq 2\pi \int_R^\infty (|\partial_r w(r,0)|^2 + |\partial_t w(r,0)|^2) dr
\]
holds for all $t < 0$.
\end{lemma}
\paragraph{Definition of $V_R (x,t)$} Let us define ($R>0$)
\begin{equation}
V_R (x,t) = \left\{\begin{array}{l} W_0 (R+|t|), \,\, \hbox{if}\,\, |x|\leq R+|t|;\\
W_0 (|x|),\,\, \hbox{if}\,\, |x|> R+|t|. \end{array}\right.
\end{equation}
Now let us consider the norms of $V_R$. By (\ref{behaviour of w at infinity}), we have
\[
 W_0 (x) \leq \frac{C_R}{|x|}.
\]
for each $|x|\geq R$. Thus if $\frac{3}{r} + \frac{1}{q} < 1$,
\begin{eqnarray*}
 \|V_R\|_{L^q L^r (\Rm \times \Rm^3)} &=& \left(\int_{\Rm}\left(\int_{\Rm^3} |V_R (x,t)|^r dx\right)^{q/r} dt \right)^{1/q}\\
 &\lesssim& \left(\int_{\Rm}\left((R+|t|)^3 |W_0(R+|t|)|^r + \int_{|x|>R+|t|} |W_0(x)|^r dx \right)^{q/r} dt \right)^{1/q}\\
 &\lesssim& C_R \left(\int_{\Rm}\left((R+|t|)^{3-r} + \int_{|x|>R+|t|} |x|^{-r} dx \right)^{q/r} dt \right)^{1/q}\\
 &\lesssim_r& C_R \left(\int_{\Rm}\left((R+|t|)^{3-r}\right)^{q/r} dt \right)^{1/q}\\
% &\lesssim_r& C_R \left(\int_{\Rm}(R+|t|)^{(3-r)q/r} dt \right)^{1/q}\\
 &\lesssim_{r,q}& C_R \left(R^{(3-r)q/r+ 1} \right)^{1/q}\\
 &\lesssim_{r,q}& C_R R^{\frac{3}{r} + \frac{1}{q} - 1}.
\end{eqnarray*}
Thus the following norms are all finite for $R>0$.
\[
 \|V_R\|_{Y_{s_p}(\Rm)} < \infty;\,\,\, \|V_R\|_{L^{2p/(p-3)} L^{2p} (\Rm \times \Rm^3)} < \infty.
\]
Furthermore, if $R$ is sufficiently large $R > R'$, we could choose $C_R = 2$, thus
\begin{equation}
 \|V_R\|_{Y_{s_p}(\Rm)} \lesssim R^{\frac{1}{2} -s_p};\,\,\, \|V_R\|_{L^{2p/(p-3)} L^{2p} (\Rm \times \Rm^3)} \lesssim R^{-1/2}. \label{X bound of VR}
\end{equation}

\subsection{Approximation Theory}
\begin{theorem} \label{approximation theory}
Fix $3<p<5$. There exists a constant $\delta_0>0$, such that if $\delta < \delta_0$ and we have\\
(i)A function $V(x,t) \in L^{2p/(p-3)} L^{2p} (I \times \Rm^3)$  with $\|V(x,t)\|_{Y_{s_p}(I)} < \delta$.
Here $I$ is a time interval containing $0$;\\
(ii) A pair of initial data $(h_0, h_1)$ with
\[
\|(h_0,h_1)\|_{\dot{H}^1 \times L^2 (\Rm^3)} < \delta,\,\,\, \|(h_0,h_1)\|_{\dot{H}^{s_p} \times \dot{H}^{s_p -1} (\Rm^3)} < \delta.
\]
Then the equation
\[
\left\{\begin{array}{l} \partial_t^2 h - \Delta h = F(V+h) -F(V), \,\,\,\, (x,t)\in \Rm^3 \times I;\\
h|_{t=0} = h_0;\\
\partial_t h|_{t=0} = h_1\end{array}\right.
\]
has a unique solution $h(x,t)$ on $I \times \Rm^3$ so that
\[
 \|h\|_{Y_{s_p}(I)} \leq C_p \delta;
\]
\[
 \sup_{t \in I} \|(h, \partial_t h) - (h_L,\partial_t h_L)\|_{\dot{H}^1 \times L^2} \leq C_p \delta^{p-1} \|(h_0,h_1)\|_{\dot{H}^1 \times L^2}.
\]
Here $(h_L, \partial_t h_L)$ is the solution of the linear wave equation with initial data $(h_0,h_1)$.
\end{theorem}
\paragraph{Stretch of Proof}
In this proof $C_p$ represents a constant that depends on $p$ only. In different places $C_p$ may represent different constants.
We will also write $Y$ instead of $Y_{s_p}(I)$ for convenience. By the Strichartz estimates, we have
\[
 \|F(V+h) -F(V)\|_{Z_{s_p}} \leq C_p \|h\|_Y (\|h\|_Y^{p-1} + \|V\|_Y^{p-1});
\]
\[
 \|F(V+h^{(1)}) - F(V+h^{(2)})\|_{Z_{s_p}} \leq C_p \|h^{(1)} - h^{(2)}\|_Y (\|h^{(1)}\|_Y^{p-1} +
 \|h^{(2)}\|_Y^{p-1}+ \|V\|_Y^{p-1}).
\]
By a fixed point argument, if $\delta$ is sufficiently small, we have a unique solution $h(x,t)$ defined on $I \times \Rm^3$,
so that $\|h\|_Y \leq C_p \delta$. Now by Strichartz estimates, if $\delta$ is sufficiently small
\begin{eqnarray*}
 \|h\|_{L^\frac{4p}{9-p} L^\frac{4p}{p-3}} &\leq&  C_p (\|(h_0,h_1)\|_{\dot{H}^1\times L^2} + \|F(V+h)-F(V)\|_{L^1 L^2})\\
 &\leq& C_p \left(\|(h_0,h_1)\|_{\dot{H}^1\times L^2}+ \|h\|_{L^\frac{4p}{9-p} L^\frac{4p}{p-3}} (\|h\|_Y^{p-1} + \|V\|_Y^{p-1})\right)\\
 &\leq& C_p \|(h_0,h_1)\|_{\dot{H}^1\times L^2} + C_p \delta^{p-1}\|h\|_{L^\frac{4p}{9-p} L^\frac{4p}{p-3}}\\
 &\leq& C_p \|(h_0,h_1)\|_{\dot{H}^1\times L^2} + (1/2)\|h\|_{L^\frac{4p}{9-p} L^\frac{4p}{p-3}}
\end{eqnarray*}
Thus
\[
 \|h\|_{L^\frac{4p}{9-p} L^\frac{4p}{p-3}} \leq C_p \|(h_0,h_1)\|_{\dot{H}^1\times L^2}.
\]
This gives us
\begin{eqnarray*}
 \sup_{t \in I} \|(h, \partial_t h) - (h_L,\partial_t h_L)\|_{\dot{H}^1 \times L^2} &\leq& C_p \|F(V+h)-F(V)\|_{L^1 L^2}\\
 &\leq& C_p \|h\|_{L^\frac{4p}{9-p} L^\frac{4p}{p-3}} (\|h\|_Y^{p-1} + \|V\|_Y^{p-1})\\
 &\leq& C_p \delta^{p-1} \|(h_0,h_1)\|_{\dot{H}^1\times L^2}.
\end{eqnarray*}
%\paragraph{Remark} We can make the theorem work for the whole interval $1/2 < s_p < 1$, if we use the $S(I)$ and $W(I)$ norms of $V(x,t)$
%instead of $X(I)$ norm. The proof is similar.
\subsection{Match with $W_0 (x)$}
Using the estimate (\ref{strong decay of w energy}), we have
\[
 \int_{r_0}^{4 r_0} |\partial_r w(r,t)| dr \lesssim \left(\int_{r_0}^{4 r_0} |\partial_r w(r,t)|^2 dr\right)^{1/2} r_0^{1/2}
 \lesssim \frac{1}{r_0^{p-3}}.
\]
This means
\begin{equation}
 \int_{r_0}^\infty |\partial_r w(r,t)| dr \lesssim  \frac{1}{r_0^{p-3}}. \label{lin6}
\end{equation}
Thus we know the limit $\lim_{r \rightarrow \infty} w(r,t)$ exists for each $t$.

\paragraph{Case 1} If $\lim_{r \rightarrow \infty} w(r,0) =0$. Then in the rest of this section,
set $W(x) = 0$. By (\ref{lin6}) we have
\[
 \left| w(r,0)\right| \lesssim \frac{1}{r^{p-3}}.
\]
Thus
\[
 \left| u_0(x) - W(x)\right| = \frac{1}{|x|} \left| w(|x|,0)\right| \lesssim \frac{1}{|x|^{p-2}}.
\]
%\[
% f_1 (r_0) = \sup_{|x| \geq r_0} |x||u(x,t)| = \sup_{r \geq r_0} |w(r,t)| \leq C_1.
%\]
%Then follow the remark after (\ref{decay of energy with f}), we have
%\[
% \left(\int_{r_0}^{4r_0} (\partial_t w(r,t))^2 + (\partial_r w(r,t))^2 dr\right)^{1/2} \leq C_p f_1 (r_0)^p \frac{1}{r_0^{p -5/2}}.
%\]
%Now by $\lim_{r \rightarrow \infty} w(r,t) =0$, we have for each $r'\geq r_0$ and time $t$
%\begin{eqnarray*}
% |w(r',t)| &\leq& \int_{r_0}^\infty |\partial_r w(r,t)| dr\\
% &\leq& \left(\int_{r_0}^\infty (\partial_r w(r,t))^2 dr\right)^{1/2}
%\left(\int_{r_0}^\infty \frac{1}{r^2} dr\right)^{1/2}\\
% &\leq& C_p f_1^p (r_0) \frac{1}{r_0^{p -5/2}} r_0^{1/2}\\
% &\leq& \frac{C_p f_1^p (r_0)}{r_0^{p-3}}.
%\end{eqnarray*}
%Take the sup, we have
%\[
% f_1 (r_0) &\leq& \frac{C_p f_1^p (r_0)}{r_0^{p-3}}.
%\]
%for each $r_0>0$.

\paragraph{Case 2} If $\lim_{r \rightarrow \infty} w(r,0) \neq 0$. WLOG, let
us assume the limit is equal to $1$. Otherwise we only need to apply some space-time dilation and/or multiplication by $-1$ on $u$.\\
In the rest of this section, set $W(x)=W_0 (x)$.
Thus by (\ref{lin6}) we have
\[
 |w(r,0) -1| \leq \int_{r}^\infty |\partial_r w(r,t)| dr \lesssim  \frac{1}{r^{p-3}}.
\]
Dividing this inequality by $r$, we have
\[
 \left|u_0 (x) - \frac{1}{|x|}\right| \lesssim \frac{1}{|x|^{p-2}}.
\]
Combining this with our estimate for $W_0 (x)$, we have for large $x$
\[
 \left| u_0(x) - W(x)\right| \lesssim \frac{1}{|x|^{p-2}}.
\]
%\subsection{Preliminary Results }
%\paragraph{Method of Center Cutoff} Let $(v_0, v_1) \in \dot{H}^1 \times L^2$ be a pair radial functions. We define $R>0$
%\[
% (\Psi_R v_0)(x) = \left\{\begin{array}{l} v_0 (x),\,\, \hbox{if}\,\, |x| > R;\\
%  v_0 (R)\,\, \hbox{if}\,\, |x|\leq R. \end{array}\right.
%\]
%\[
% (\Psi_R v_1)(x) = \left\{\begin{array}{l} v_1 (x),\,\, \hbox{if}\,\, |x| > R;\\
%  0,\,\, \hbox{if}\,\, |x|\leq R. \end{array} \right.
%\]
\subsection{Identity near infinity}
\begin{theorem} \label{identity near infinity}
Let $W(x)= W_0 (x)$ or $W(x)=0$. Suppose $u(x,t)$ is a global radial solution of the equation (\ref{eqn}) with initial data $(u_0,u_1)$
satisfying the following conditions.\\
(I) $(u_0, u_1) \in \dot{H}^{s_p}\times \dot{H}^{s_p-1}$.\\
(II) %The solution is in the energy space at least locally.
The following inequality holds for each $t\in \Rm$ and $r>0$.
\begin{equation}
 \int_{r <|x|<4r} (|\nabla u(x,t)|^2 + |\partial_t u(x,t)|^2 ) dx \leq C_1 r^{-1}. \label{con1}
\end{equation}
(III) We have $u_0(x)$ and $W(x)$ are very close to each other as $|x|$ is large.
\begin{equation}
 |u_0 (x) - W (x)| \lesssim \frac{1}{|x|^{p-2}}. \label{con2}
\end{equation}
Then there exists $R_0 = R_0 (C_1, p)$ such that the pair $(u_0(x) - W(x), u_1(x))$ is essentially supported in the ball $B(0,R_0)$.
\end{theorem}
\paragraph{Remark} There are actually two separate theorems, both could be proved in the same way. If
$W(x) = W_0(x)$(the primary case), then define $V_{R_0}$ as usual in the proof below. Otherwise if $W(x)=0$, just make
$V_{R_0} =0$.
%\paragraph{Remark 2} The following proof uses the $X(I)$ norm of $V_{R_0}$. This won't work if $s_p$ is close to $1/2$, but we can still obtain the
%theorem following the same argument but using $S(I)$, $W(I)$ norms instead for $1/2 < s_p < 1$.
\paragraph{Proof} Let us define for $R \geq R_0$
\[
 g_0 = \Psi_{R} (u_0 - W);\,\,\, g_1 = \Psi_{R} u_1.
\]
\[
 G(r) = u_0(r) - W(r).
\]
Choose a small constant $\delta = \delta(p)$, so that it is smaller than the constant $\delta_0$ in theorem \ref{approximation theory}
and guarantees the number $C_p \delta^{p-1}$ in the conclusion of that theorem is smaller than $\eps(p)$,
which is a small number determined later in the argument below. By the condition (\ref{con1}) and the properties of $W(x)$, we know ($R>1$)
\begin{eqnarray*}
 \int_{\Rm^3} (|\nabla g_0|^2 + g_1^2) dx &\lesssim_{C_1,p}& R^{-1};\\
 \int_{\Rm^3} \left(|\nabla g_0|^{3(p-1)/(p+1)} + g_1^{3(p-1)/(p+1)}\right) dx &\lesssim_{C_1,p}& R^{-3(p-3)/(p+1)}.
\end{eqnarray*}
As a result, if $R_0 = R_0 (C_1, p)$ is sufficiently large,
we have the following inequalities hold as long as $R\geq R_0$. (We use the Sobolev embedding in order to obtain the second inequality)
\[
 \|(g_0 , g_1)\|_{\dot{H}^{1} \times L^2} \leq \delta;\;\;\;\;
 \|(g_0 , g_1)\|_{\dot{H}^{s_p} \times \dot{H}^{s_p -1}} \leq \delta;
\]
\[
 \|V_{R_0}\|_{Y_{s_p}(\Rm)} \leq \delta.
\]
Let $g$ be the solution of
\[
 \partial_t^2 g - \Delta g = F (V_{R_0} + g) - F(V_{R_0})
\]
with the initial data $(g_0,g_1)$ and $\tilde{g}$ be the solution of the linear wave equation with the same initial data.
On the other hand, we know $u (x,t)- W (x)$ is the solution of the equation
\begin{equation}
 \partial_t^2 \tilde{u} - \Delta \tilde{u} = F (W + \tilde{u}) - F(W) \label{s eqn}
\end{equation}
% u- V'_{R_0} is the solution of
% Box \tilde{u} = F (V'_{R_0} + \tilde{u}) - \Box V'_{R_0};
% Here V'_{R_0} is the smooth center cutoff of W.
in the domain $\Rm \times (\Rm^3 \setminus \{0\})$ with the initial data $(u_0 - W, u_1)$. Let $K$ be the domain
\[
 K = \{(x,t): |x| > |t| + R\}.
\]
Considering the fact $W(x) = V_{R_0}(x,t)$ in the region $K$ and the construction of $(g_0,g_1)$, we have
\[
 u(x,t) - W (x) = g (x,t);\,\,\, \partial_t u(x,t) = \partial_t g(x,t)
\]
in the domain $K$ by the finite speed of propagation. Using our assumption (\ref{con1}) and the decay of $W (x)$ at infinity, we have
\begin{equation}
 \lim_{t \rightarrow \pm \infty} \int_{|x| > |t|+R} (|\nabla g(x,t)|^2 + |\partial_t g(x,t)|^2) dx \rightarrow 0. \label{lin1}
\end{equation}
Using lemma \ref{energy channel}, WLOG, let us assume for all $t > 0$
\[
 \int_{|x|> R + t} (|\nabla \tilde{g}(x,t)|^2 + |\partial_t \tilde{g}(x,t)|^2) dx \geq 2\pi \int_R^\infty
 \left(\left|\partial_r (r g_0(r,0))\right|^2 + r^2 |g_1 (r,0)|^2\right) dr.
\]
That is
\[
 \int_{|x|> R + t} (|\nabla \tilde{g}(x,t)|^2 + |\partial_t \tilde{g}(x,t)|^2) dx \geq \frac{1}{2} \left(\int_{|x|> R}
 (|\nabla g_0|^2 + g_1^2) dx \right)- 2 \pi R g_0^2 (R).
\]
Combining this with (\ref{lin1}), we have
\begin{eqnarray*}
 && \liminf_{t \rightarrow \infty} \|(g(x,t), \partial_t g(x,t)) - (\tilde{g}, \partial_t \tilde{g})\|_{\dot{H}^1 \times L^2 (|x|> R+ t)}\\
 &\geq& \left(\frac{1}{2}\int_{|x|> R} (|\nabla g_0|^2 + g_1^2) dx - 2 \pi R g_0^2 (R) \right)^{1/2}.
\end{eqnarray*}
On the other hand, we have the following inequality by theorem \ref{approximation theory}
\[
 \|(g(x,t), \partial_t g(x,t)) - (\tilde{g}, \partial_t \tilde{g})\|_{\dot{H}^1 \times L^2} \leq C_p \delta^{p-1} \|(g_0, g_1)\|_{\dot{H}^1 \times L^2}
 \leq \eps(p) \|(g_0, g_1)\|_{\dot{H}^1 \times L^2}.
\]
Considering both inequalities above, we have
\[
 \frac{1}{2} \int_{|x|> R} (|\nabla g_0|^2 + g_1^2) dx - 2 \pi R g_0^2 (R)  \leq \eps^2(p) \int_{|x|> R} (|\nabla g_0|^2 + g_1^2) dx.
\]
Thus
\begin{equation}
 \int_{|x|> R} (|\nabla g_0|^2 + g_1^2) dx \leq \frac{4\pi}{1 - 2\eps^2(p)} R g_0^2 (R). \label{lin2}
\end{equation}
We have
\begin{eqnarray*}
 |g_0 (m R) - g_0 (R)| &\leq& \int_R^{m R} |\partial_r g_0| dr\\
 &\leq& \left(\int_R^{m R} |r \partial_r g_0|^2 dr\right)^{1/2} \left(\int_R^{m R} \frac{1}{r^2} dr\right)^{1/2}\\
 &\leq& \left(\frac{1}{4\pi} \int_{|x|> R} (|\nabla g_0|^2 + g_1^2) dx\right)^{1/2} \left(\frac{1}{R} - \frac{1}{mR}\right)^{1/2}\\
 &\leq& \left(\frac{R g_0^2 (R)}{1 -2\eps^2 (p)}\right)^{1/2} \left(1 - \frac{1}{m}\right)^{1/2} R^{-1/2} \\
 &\leq& \left(\frac{1-1/m}{1 - 2\eps^2(p)}\right)^{1/2} |g_0 (R)|.
\end{eqnarray*}
By the fact $p-2 > 1$, we can choose $k = k(p) \in {\mathbb Z}^+$ such that $(k+1)/k < p-2$. Let $m = 2^k$. Since
\[
 (1 - 1/m)^{1/2} < 1 - \frac{1}{2m},
\]
we can choose $\eps(p)$ so small that
\[
 \left(\frac{1-1/m}{1 - 2\eps^2(p)}\right)^{1/2} \leq 1 - \frac{1}{2m} = 1 - \frac{1}{2^{k+1}}.
\]
Plugging this in our estimate above, we obtain
\[
 |g_0 (2^k R) - g_0 (R)| \leq (1 - \frac{1}{2^{k+1}})|g_0 (R)|.
\]
Thus
\[
 |g_0 (2^k R)| \geq \frac{1}{2^{k+1}}|g_0 (R)|.
\]
By the definition of $g_0$, this is the same as
\[
 |G(2^k R)| \geq \frac{1}{2^{k+1}} |G(R)|.
\]
This inequality holds for all $R \geq R_0$. Now let us consider the value of $G(R_0)$. If $G(R_0)=0$, let us choose $R = R_0$. Thus we have
$g_0(R)= 0$. Plugging this back in (\ref{lin2}), we have $(g_0,g_1) = (0,0)$. This means that $(u_0 -W, u_1)$ is supported in $B(0,R_0)$ and
finishes the proof. If $|G(R_0)| > 0$, then we have
\[
 |G(2^{kn} R_0)| \geq \frac{1}{(2^{kn})^{(k+1)/k}} |G(R_0)|>0
\]
for each positive integer $n$. This contradicts with the condition (\ref{con2}) because $(k+1)/k < p-2$ by our choice of $k$.
\paragraph{Remark} If one feels uncomfortable about the singularity at zero in the equation (\ref{s eqn}), we could use the following center-cutoff
version instead. Let $\varphi$ be a smooth, radial, nonnegative function satisfying
\[
 \varphi(x) = \left\{\begin{array}{l} 1,\;\hbox{if}\; |x| \geq 1;\\ \in [0,1],\; \hbox{if}\; |x| \in (1/2,1);\\
 0,\; \hbox{if}\;|x| \leq 1/2. \end{array}\right.
\]
Then $u(x,t)- \varphi(|x|/R_0) W_0(x)$ is a solution to the equation
\[
 \left\{\begin{array}{l} \partial_t^2 \tilde{u} - \Delta \tilde{u} = F(\varphi(|x|/R_0) W_0 + \tilde{u}) + \Delta (\varphi(|x|/R_0) W_0(x)),
 (x,t) \in \Rm^3 \times \Rm;\\
 \tilde{u}|_{t=0} = u_0 - \varphi(|x|/R_0) W_0 \in \dot{H}^{s_p}(\Rm^3);\\
 \partial_t \tilde{u}|_{t=0} = u_1 \in \dot{H}^{s_p-1}(\Rm^3).\end{array}\right.
\]
For any $T>0$, we know
\[
 \|\varphi(|x|/R_0) W_0(x)\|_{Y_{s_p}([-T,T])} < \infty;\;\;\; \|\Delta (\varphi(|x|/R_0) W_0(x))\|_{Z_{s_p}([-T,T])} < \infty.
\]
In addition, the function $\Delta (\varphi(|x|/R_0) W_0(x)) = - F(W_0(x))$ in the region $K$. We can do the argument as usual in the proof above but
avoid the singularity at zero with this new cutoff version of the equation (\ref{s eqn}).
This method also works in the proof of theorem \ref{behavior of compact solution}, which will be introduced in the next subsection.
\paragraph{Application of the theorem} Now apply theorem \ref{identity near infinity} to our soliton-like minimal blow-up solution.
All the conditions are satisfied by our earlier argument. Thus $(u_0(x)-W(x), u_1(x))$ is supported in the ball of radius $R_0$ centered at the origin.
In particular, because $R_0$ depends only on the constant $C_1$ and $p$, the same $R_0$ also works for other time $t$ as long as the condition
(\ref{con2}) is true at that time. But by the finite speed of propagation, we know $(u(x,t)-W(x), \partial_t u(x,t))$ is compactly supported in
$B(0, R_0 +|t|)$ at each time $t$. This means the condition (\ref{con2}) is always true at any time.
Thus the pair $(u(x,t)-W(x), \partial_t u)$ is supported in the cylinder $B(0,R_0) \times \Rm$.

\subsection{Local Radius Analysis}
Let us define the essential radius of the support of $(u(x,t) - W (x), \partial_t u(x,t))$ at time $t$.
\[
 R(t) = \min \{R \geq 0: (u(x,t) - W (x), \partial_t u(x,t))=(0,0)\, \hbox {holds for}\, |x|> R\}.
\]
This is well-defined for our minimal blow-up solution. Actually $R(t) \leq R_0$ holds for all $t$.

\begin{theorem} \emph{\textbf{(Behavior of compactly supported solutions)}} \label{behavior of compact solution}
Let $W(x) = W_0(x)$ or $W(x)=0$. Let $u(x,t)$ be a radial solution of the equation (\ref{eqn}) in a time interval $I$, so that\\
(I) The pair $(u(x,t), \partial_t u(x,t)) \in \dot{H}^1 (\Rm^3) \times L^2 (\Rm^3)$ for each $t \in I$.\\
(II) The pair $(u(x,0) - W (x), \partial_t u(x,0))$ is compactly supported with an essential radius
of support $R(0) > R_1 >0$.
Then there exists a constant $\tau=\tau (R_1,p)$, such that
\[
 R(t) = R(0) + |t|
\]
holds either for each $t \in [0,\tau]\cap I$ or for each $t \in [-\tau,0]\cap I$.
\end{theorem}
\paragraph{Remark}  If $W(x) = W_0(x)$(the primary case), then define $V_{R_1}$ as usual in the proof.
Otherwise if $W(x)=0$, just make $V_{R_1} =0$. In this case we can choose $\tau = \infty$.
\paragraph{Proof} By our previous argument, we have $\|V_{R_1}\|_{Y_{s_p}(\Rm)} < \infty$. Thus we can choose
$\tau = \tau(R_1,p) > 0$ such that $\|V_{R_1}\|_{Y_{s_p}([-\tau,\tau])} < \delta$. Here $\delta$ is a small constant so that
we can apply theorem \ref{approximation theory} and make the number $C_p \delta^{p-1} < 1/{100}$ in that theorem. If $\eps < R(0)-R_1$,
let us define a pair of initial data $(g_0,g_1)$ for each $R \in (R(0)-\epsilon, R(0))$
\[
 g_0 = \Psi_{R} (u_0 - W);\,\,\, g_1 = \Psi_{R} u_1.
\]
This pair $(g_0(x) ,g_1(x))$ is nonzero by the definition of $R(0)$.\\
By our assumptions on $(u_0,u_1)$, we know the following inequalities hold for each $R \in (R(0)-\epsilon, R(0))$
as long as $\eps$ is sufficiently small.
(In order to obtain the second inequality we use the Sobolev embedding)
\[
 \|(g_0 ,g_1)\|_{\dot{H}^1 \times L^2} < \delta;
\]
\[
 \|(g_0, g_1)\|_{\dot{H}^{s_p} \times \dot{H}^{s_p -1}} < \delta.
\]
Furthermore, we have
\begin{eqnarray*}
 |g_0 (R)| &=& \left|g_0 (R(0)) - \int_R^{R(0)} \partial_r g_0 (r) dr\right|\\
  &\leq& \int_R^{R(0)} |\partial_r g_0 (r)| dr\\
  &\leq& \left(\int_R^{R(0)} r^2 |\partial_r g_0 (r)|^2 dr \right)^{1/2}\left(\int_R^{R(0)} \frac{1}{r^2} dr\right)^{1/2}\\
  &\leq& \left(\int_R^{R(0)} r^2 |\partial_r g_0 (r)|^2 dr \right)^{1/2}\left(\frac{R(0)- R}{R(0) R}\right)^{1/2}\\
  &\leq& \left(\frac{\eps}{R(0) R}\int_R^{R(0)} r^2 |\partial_r g_0 (r)|^2 dr \right)^{1/2}.
\end{eqnarray*}
Thus
\[
 R g_0^2 (R) \leq \frac{\eps}{R(0)} \int_R^{R(0)} r^2 |\partial_r g_0 (r)|^2 dr \leq
 \frac{\eps}{4\pi R(0)} \int_{R<|x|<{R(0)}} (|\nabla g_0(x)|^2 +|g_1 (x)|^2) dx.
\]
If $\eps$ is sufficiently small, we can apply lemma \ref{eqn between w and u} to obtain
\[
 \int_R^{R(0)} \left[\left|\partial_r (r g_0(r))\right|^2 + r^2 g_1 (r)^2\right] dr \geq \frac{0.99}{4\pi}
\int_{R<|x|<R(0)} (|\nabla g_0(x)|^2 +|g_1 (x)|^2) dx.
\]
Let $\tilde{g} (x,t)$ be the solution to the linear wave equation with the initial data $(g_0,g_1)$.
By lemma \ref{energy channel},
\begin{eqnarray*}
 \int_{|x|> R +|t|} \left(|\nabla \tilde{g}(x,t)|^2 + |\partial_t \tilde{g}(x,t)|^2 \right) dx &\geq& 2\pi
 \int_R^{\infty} \left[\left|\partial_r (r g_0(r))\right|^2 + r^2 |g_1 (r)|^2\right] dr\\
 &=& 2\pi \int_R^{R(0)} \left[\left|\partial_r (r g_0(r))\right|^2 + r^2 |g_1 (r)|^2\right] dr\\
 &\geq& 0.49 \int_{R<|x|<R(0)} (|\nabla g_0(x)|^2 +|g_1 (x)|^2) dx
\end{eqnarray*}
holds either for each $t\geq 0$ or for each $t\leq 0$. WLOG, let us choose $t \geq 0$, then we have
\begin{equation}
 \|(\tilde{g}(x,t),\partial_t \tilde{g}(x,t))\|_{\dot{H}^1 \times L^2 (|x|> R+t)} \geq 0.7
\|(g_0 ,g_1)\|_{\dot{H}^1 \times L^2 (\Rm^3)}. \label{lin4}
\end{equation}
Let $g$ be the solution of the following equation
\[
\left\{\begin{array}{l} \partial_t^2 g - \Delta g = F(V_{R_1} +g) -F(V_{R_1}), \,\,\,\, (x,t)\in \Rm^3 \times I;\\
g|_{t=0} = g_0;\\
\partial_t g|_{t=0} = g_1.\end{array}\right.
\]
By lemma \ref{approximation theory}, we have
\[
 \|(g(x,t), \partial_t g(x,t)) - (\tilde{g}(x,t), \partial_t \tilde{g}(x,t))\|_{\dot{H}^1 \times L^2} \leq 0.01
 \|(g_0 ,g_1)\|_{\dot{H}^1 \times L^2 (\Rm^3)}
\]
for each $t \in [-\tau, \tau]$. Combining this with (\ref{lin4}), for $t\in [0,\tau]$ we obtain
\begin{equation}
 \|(g(x,t), \partial_t g(x,t))\|_{\dot{H}^1 \times L^2(|x|>R + t)} \geq 0.69  \|(g_0 ,g_1)\|_{\dot{H}^1 \times L^2 (\Rm^3)}.
\label{lin5}
\end{equation}
In addition, we know $u(x,t)- W(x)$ is the solution of equation
\[
\left\{\begin{array}{l} \partial_t^2 \tilde{u} - \Delta \tilde{u} = F(W(x) +\tilde{u}) -F(W(x)), \,\,\,\, (x,t)\in \Rm^3 \times I;\\
\tilde{u}|_{t=0} = u_0 - W;\\
\partial_t \tilde{u}|_{t=0} = u_1\end{array}\right.
\]
in $(\Rm^3 \setminus \{0\}) \times I$. The initial data of these two equations is the same in the region
$\{x: |x| \geq R\}$ and the nonlinear part is the same function in the region
\[
 K = \{(x,t): |x|>R+t,\, t\in [0,\tau]\cap I\}
\]
Thus by the finite speed of propagation, we have $g(x,t) = u(x,t) - W (x)$ and
$\partial_t g(x,t) = \partial_t u(x,t)$ in $K$. Plugging this in (\ref{lin5}), we obtain
\[
 \|(u(x,t) - W (x), \partial_t u(x,t))\|_{\dot{H}^1 \times L^2(|x|>R + t)} \geq 0.69  \|(g_0 ,g_1)\|_{\dot{H}^1 \times L^2 (\Rm^3)}
\]
for each $t \in I \cap [0,\tau]$. Since $R < R(0)$, we know the right hand of the inequality above is positive by the definition
of essential radius of support. Thus we have
\begin{equation}
 R(t) \geq R + |t| \label{inequality1}
\end{equation}
for all $t \in [0,\tau] \cap I$. Letting $R \rightarrow R(0)^-$, we obtain $R(t) \geq R(0) + |t|$. By the finite speed
of propagation, we have $R(t)= R(0) + |t|$.
\paragraph{Remark} For each $R \in (R(0)-\eps, R(0))$, we know that the inequality (\ref{inequality1}) above holds either in the positive
or negative time direction. It may work in different directions as we choose different $R$'s. However, we
can always choose a sequence $R_i \rightarrow R(0)^-$ such that the inequality works in the same time direction for
all $R_i$'s. This is sufficient for us to conclude the theorem.

\subsection{End of Soliton-like Solution}
Now let us show $R(0)= 0$. If it was not zero, let $R_1 = R(0)/2$, and then apply theorem \ref{behavior of compact solution}.
We have (WLOG) $R(t)= R(0) + t$ for each $t \in [0,\tau]$.
Applying theorem \ref{behavior of compact solution} again at $t = \tau$, we obtain
\[
 R(t) = R(0) + \tau + (t-\tau) = R(0) + t
\]
for $t\in [\tau,2\tau]$, because (i) The same constant $\tau$ works by $R(\tau) > R(0)> R_1$; (ii) The theorem may only work
in the positive time direction, since we know the radius of support $R(t)$ decreases in the other direction.
Repeating this process, we have for each $t>0$,
\[
 R(t) = R(0) + t.
\]
But it is impossible since $R(t)$ is uniformly bounded by $R_0$.
Therefore we must have $R(0)=0$. But this means either $u_0 = W_0 (x) \notin \dot{H}^{s_p}(\Rm^3)$ or
$(u_0,u_1)=(0,0)$. This is a contradiction.

\section{The Solution of the Elliptic Equation}
In this section we will consider the elliptic equation
\begin{equation}
 -\Delta W(x) = |W(x)|^{p-1} W(x).
\end{equation}
It has infinitely many solutions. For example,
\[
 W_1 (x) = C |x|^{-2/(p-1)}
\]
is a solution if we choose an appropriate constant $C$. We are interested in radial solutions of this elliptic equation. Let us assume
$W(x) = y (|x|)$. Here $y(r)$ is a function defined in $(0,\infty)$. The function $y(r)$ satisfies the equation
\begin{equation}
 y'' (r) + \frac{2}{r} y'(r) + |y|^{p-1} y(r) = 0.
\end{equation}
Let us show that the solution $W_0(x)$ we mentioned earlier in this paper exists.
\subsection{Existence of $W_0 (x)$}
We are seeking a solution with the property $W_0(x) \simeq 1/|x|$ as $x$ is large. That is equivalent to $y(r) \simeq 1/r$. Let us
define $\rho(r) = r y(r)$, then $\rho(r)$ satisfies
\[
 \rho''(r) = - \frac{F(\rho)}{r^{p-1}};\;\;\;\; F(\rho) = |\rho|^{p-1}\rho.
\]
We expect $\rho(r) \simeq 1$ for large $r$'s, thus let us assume $\rho(r) = \phi(r) + 1$. The corresponding equation for $\rho(r)$ is given as below
\[
 \phi''(r) = - \frac{F(\phi +1)}{r^{p-1}}.
\]
The idea is to show
\begin{itemize}
\item (I) This equation has a solution in the interval $[R, \infty)$ with boundary conditions at infinity
$\phi(+\infty) = \phi'(+\infty)= 0$, by a fixed point argument.
\item (II) We can expand the domain of this solution to $\Rm^+$.
\end{itemize}
\paragraph{The Fixed Point Argument} Let us consider the metric space
\[
 K = \{\phi: \phi \in C([R, \infty); [-1,1]),\, \lim_{r \rightarrow +\infty}\phi(r) = 0 \}.
\]
with the distance $d(\phi_1, \phi_2) = \sup_r |\phi_1(r)-\phi_2 (r)|$. One can check $K$ is complete.
Let us define a map $L: K \rightarrow K$ by
\[
 L(\phi)(r) = \int_r^\infty \left(\int_s^\infty \left(- \frac{F(\phi(t) +1)}{t^{p-1}}\right) dt\right) ds.
\]
We have
\[
 |L(\phi)(r)| \leq \int_r^\infty \left(\int_s^\infty \left(\frac{2^p}{t^{p-1}}\right) dt\right) ds \leq \frac{C_p}{r^{p-3}};
\]
\[
 |L(\phi_1)(r) - L(\phi_2)(r)|\leq C_p \int_r^\infty \left(\int_s^\infty \left(\frac{d(\phi_1,\phi_2)}{t^{p-1}}\right) dt\right)ds
 \leq C_p \frac{d(\phi_1,\phi_2)}{r^{p-3}}.
\]
Thus if $R > R(p)$ is a sufficiently large number, then $L$ is a contraction map from $K$ to itself. As a result, 
there exists a unique fixed point $\phi_0 (r)$.
This gives us a classic smooth solution of the ODE in $[R,\infty)$. We have $\phi_0 (r) \lesssim r^{3-p}$ and its derivative $\phi'_0 (r)$ satisfies
\[
|\phi'_0 (r)| = \left|\int_r^\infty \left(\frac{F(\phi_0(t) +1)}{t^{p-1}}\right) dt\right| \leq \frac{C_p}{r^{p-2}}.
\]
%This is a perfect solution at $r \geq R$.
\paragraph{Expansion of the Solution} Now let us solve the ODE backward from $r =R$.
We need to show it will never break down before we approach zero. Actually we have
\[
 \frac{d}{dr}\left(\frac{|\phi_0 + 1|^{p+1}}{p+1} + \frac{r^{p-1}|\phi'_0|^2}{2}\right) = \frac{p-1}{2} r^{p-2}|\phi'_0|^2 \geq 0.
\]
Thus we have the following inequality holds for all $0 < r \leq R$ as long as the solution still exists at $r$
\[
 \frac{|\phi_0(r) + 1|^{p+1}}{p+1} + \frac{r^{p-1}|\phi'_0 (r)|^2}{2} \leq C.
\]
But this implies the solution will never break down at a positive $r$. Let us define
\[
 W_0 (x) = \frac{\phi_0(|x|)+ 1}{|x|}.
\]
This is a $C^2$, radial solution of our elliptic equation for $|x|>0$. Furthermore, we have for large $x$
\[
 \left|W_0 (x) - \frac{1}{|x|}\right| = \frac{|\phi_0 (|x|)}{|x|} \leq \frac{C_p}{|x|^{p-2}};
\]
\[
 |\nabla W_0(x)| = \left|\frac{r \phi'_0 (r) - \phi_0 (r) - 1}{r^2}\right|_{r =|x|} \leq \frac{C_p}{|x|^2}.
\]
Now the remaining task is to show $W_0 (x)$ is not in the space $\dot{H}^{s_p}$.
This implies $W_0 (x)$ must have a singularity at $0$. It turns out that it is not trivial. For instance,
if we repeat the argument as above in the case $p =5$, then the solution we obtain will be a smooth function in the whole space, as below
\[
 W(x) = \frac{\sqrt{3}}{(1 + 3|x|^2)^{1/2}}.
\]
\subsection{Radial $\dot{H}^{s_p}$ Solution Does Not Exist}
\paragraph{Regularity} Let us first show any radial $\dot{H}^{s_p}$ solution of the elliptic equation must be in the space
$C^2(\Rm^3 \setminus \{0\})$.
We know a radial $\dot{H}^{s_p}$ function must be continuous except for $x =0$. Using this in the elliptic equation,
we have the solution is $C^2$ except for $x=0$.
\paragraph{Introduction to $r^\theta y(r)$} We assume $W(x)= y(|x|)$. The function $y(r)$ defined in $\Rm^+$ is a $C^2$ solution of
\[
 y'' (r) + \frac{2}{r} y'(r) + |y|^{p-1} y(r) = 0.
\]
Let us define another $C^2(\Rm^+)$ function
\[
 v(r) = r^\theta y(r),\,\,\,\, \theta = \frac{2}{p-1}.
\]
If $W(x)= y(|x|)$ is in the space $\dot{H}^{s_p}$, then we have
\[
\lim_{r \rightarrow 0^+} v(r) = \lim_{r \rightarrow +\infty} v(r) = 0.
\]
Plugging $y(r) = r^{-\theta} v(r)$ in the equation for $y(r)$, we obtain an equation for $v(r)$,
\[
 r^2 v''(r) + \frac{2(p-3)}{p-1} r v'(r) - \frac{2(p-3)}{(p-1)^2} v(r) + |v|^{p-1} v(r) = 0.
\]
Multiplying both sides by $v'(r)$, we obtain
\begin{equation}
 \frac{d}{dr} \left(r^2 \frac{|v'(r)|^2}{2} - \frac{p-3}{(p-1)^2} v^2(r) + \frac{|v(r)|^{p+1}}{p+1}\right) = \frac{5-p}{p-1} r |v'(r)|^2 \geq 0.
 \label{var of v}
\end{equation}
These identities imply
\paragraph{Claim 1} There exist $r_1, R_1 > 0$, such that the function $v(r)$ does not admit any positive local maximum
or negative local minimum in the set $K = (0,r_1) \cup (R_1, \infty)$.
Actually by the limits of the function $v(r)$ at $0^+$ and $+\infty$, we know there exist $r_1, R_1 > 0$, such that $v(r) \leq \eps$ for $r \in K$.
If $\eps$ is sufficiently small, then the sign of the sum (in the equation)
\[
 - \frac{2(p-3)}{(p-1)^2} v(r) + |v|^{p-1} v(r)
\]
is the same as $- v(r)$. If there was a positive local maximum or a negative local minimum in $K$,
we would find a contradiction by considering the sign of $v''(r)$, $v'(r)$ and $v(r)$.
\paragraph{Claim 2} Let $r_1, R_1$ be the constants in Claim 1. Then in each of the intervals $(0,r_1)$ and $(R_1, +\infty)$,
the function $v(r)$ is monotone. (Namely it is nondecreasing or nonincreasing) Suppose this was not true. WLOG, let us assume $s_1 < s_2 < s_3$
and $v(s_2) < v(s_3), v(s_1)$. If we like, we can choose $s_2$ to be a local minimum (a minimum in the interval $[s_1,s_3]$).
Then by Claim 1, $v(s_2)$ must be nonnegative. Now we obtain\\
(i) If $s_i < r_1$, then there must be a positive local maximum in $(0,s_2)$;\\
(ii) If $s_i > R_1$, then this yields a positive local maximum in $(s_2, \infty)$.\\
In both cases we have a contradiction.
\paragraph{Claim 3} If $v(r)$ is not the zero function, then at least one of the following inequality holds
\[
 \liminf_{r \rightarrow 0^+} r^2 |v'(r)|^2 > 0. \,\,\,\,\,\,\, \liminf_{r \rightarrow +\infty} r^2 |v'(r)|^2 > 0.
\]
If both of these failed, by considering the integral of (\ref{var of v}) in the interval $r = (\eps, M)$ and letting $\eps \rightarrow 0^+$ and
$M \rightarrow +\infty$, we would have
\[
 \frac{5-p}{p-1} \int_0^\infty r |v'(r)|^2 dr = 0
\]
This means $v'(r) = 0$ everywhere, so $v(r) = 0$. But we assume it is not the zero function.
\paragraph{Contradiction} If $v(r)$ is not identically zero, WLOG, let us assume
\[
 \liminf_{r \rightarrow 0^+} r^2 |v'(r)|^2 > 0.
\]
This means there exist $C>0$ and $r_2>0$, such that if $r \in (0, r_2)$, we have $r^2 |v'(r)|^2 > C$. This means $|v'(r)| > \sqrt{C} r^{-1}$.
By Claim 2 $v'(r)$ does not change its sign in the interval $(0,r_1)$. Combining this fact with the lower bound of $|v'(r)|$, we know the limit of
$v(r)$ does not exist at $0^+$. This gives us a contradiction. Therefore we have
\begin{theorem}
If $3<p<5$, then a radial $\dot{H}^{s_p}(\Rm^3)$ solution to the elliptic equation
\[
 -\Delta W(x) = |W(x)|^{p-1} W(x)
\]
must be the zero solution.
\end{theorem}
\paragraph{Conclusion} In summary, any radial nontrivial solution of our elliptic equation is not in the space $\dot{H}^{s_p}(\Rm^3)$. In particular,
$W_0(x)$ is not in the space $\dot{H}^{s_p}(\Rm^3)$. Actually we have $\limsup_{x \rightarrow 0^+} |x|^\theta |W_0 (x)|> 0$ by the argument above.
This gives us a singularity at zero.
\subsection{$W_0(x)$ is smooth in $\Rm^3\setminus \{0\}$}
In this subsection, we will discover some additional properties of the soliton $W_0(x)$. Assume that $y(r)$ and $v(r)$ are defined in the same manner
as the previous subsection.
\paragraph{$W_0(x)$ is a positive solution} If this was not true, we could assume that $v(r_0) = 0$ for some $r_0 > 0$. Then by (\ref{var of v}),
we obtain
\begin{equation}
 r^2 \frac{|v'(r)|^2}{2} - \frac{p-3}{(p-1)^2} v^2(r) + \frac{|v(r)|^{p+1}}{p+1} \geq r_0^2 \frac{|v'(r_0)|^2}{2} > 0. \label{inequality2}
\end{equation}
for each $r>r_0$. However the decay of $W_0(x)$ implies (if $r$ is large)
\[
 v(r) \lesssim r^{\theta-1};
\]
\[
 v'(r) = \theta r^{\theta-1} y(r) + r^\theta y'(r) \lesssim r^{\theta -2}.
\]
This gives us a contradiction if we consider the limit of the left hand in the inequality (\ref{inequality2}) using these estimates.
\paragraph{Remark} Due to the fact that the function $F$ is smooth in $\Rm^+$,
A direct corollary follows that the function $W_0(x)$ is smooth everywhere except for $x=0$.
\section{Appendix}
\subsection{The Duhamel Formula}
\begin{lemma} \label{compact LqLr}
Let $1/2 < s \leq 1$. If $K$ is a compact subset of $\dot{H}^s \times \dot{H}^{s-1}$ with an $s$-admissible pair $(q,r)$ so that
$q \neq \infty$, then for each $\eps > 0$, there exist two constants $M, \delta > 0$ such that
\begin{eqnarray*}
 \|S(t)(u_0,u_1)\|_{L^q L^r (J \times \Rm^3)} &+& \|S(t)(u_0,u_1)\|_{L^q L^r([M,\infty)\times \Rm^3)}\\
&+& \|S(t)(u_0,u_1)\|_{L^q L^r((-\infty, M]\times \Rm^3)} < \eps
\end{eqnarray*}
holds for any $(u_0,u_1) \in K$ and any time interval $J$ with a length $|J|\leq \delta$.
\end{lemma}
\paragraph{Proof} Given $(u_0,u_1) \in \dot{H}^s \times \dot{H}^{s-1}$, it is clear that we are able to find $M, \delta > 0$ so that the inequality
holds for this particular pair of initial data and any interval $J$ with a length $|J|\leq \delta$ by the fact $q < \infty$ and the Strichartz estimate
\[
 \|S(t)(u_0,u_1)\|_{L^q L^r (\Rm \times \Rm^3)} < \infty.
\]
If $K$ is a finite set, then we can find $M$ and $\delta$ so that they work for each pair in $K$ by taking a maximum among $M$'s and a minimum among
$\delta$'s. In the general case, we can just choose a finite subset $\{(u_{0,i},u_{1,i})\}_{i=1,2,\cdots,n}$ of $K$ such that for each
$(u_0,u_1) \in K$, there exists $1 \leq i \leq n$ with
\[
 \|S(t)(u_0-u_{0,i}, u_1 -u_{1,i})\|_{L^q L^r (\Rm\times \Rm^3)} \leq C \|(u_0-u_{0,i}, u_1 -u_{1,i})\|_{\dot{H}^s \times \dot{H}^{s-1}}
 < 0.01 \eps,
\]
and then use our result for a finite subset.
\begin{lemma} \emph{\textbf{(the Duhamel formula)}} \label{proof of Duhamel}
Let $u(x,t)$ be almost periodic modulo scaling in the interval $I = (T_-,\infty)$, namely the set
\[
 K= \left\{\left(\frac{1}{\lambda(t)^{3/2-s_p}} u\left(\frac{x}{\lambda(t)}, t\right), \frac{1}{\lambda(t)^{5/2-s_p}}
 \partial_t u\left(\frac{x}{\lambda(t)}, t\right)\right): t \in I\right\}
\]
is precompact in the space $\dot{H}^{s_p} \times \dot{H}^{s_p-1}(\Rm^3)$. Then for any time $t_0 \in \Rm$, any bounded closed interval $[a,b]$
and an $s_p$-admissible pair $(q,r)$ with $q <\infty$, we have
\[
 \lim_{T\rightarrow +\infty} \|S(t-T)(u(T), \partial_t u(T))\|_{L^q L^r ([a,b]\times \Rm^3)} =0.
\]
\[
 \hbox{weak}\,\lim_{T\rightarrow +\infty} S(t_0-T)\left(\begin{array}{c}u(T)\\ \partial_t u(T)\end{array}\right) = 0.
\]
\end{lemma}
\paragraph{Proof} We have
\begin{eqnarray*}
&&\|S(t-T)(u(T), \partial_t u(T))\|_{L^q L^r ([a,b]\times \Rm^3)}\\ &=& \|S(t)(u(T), \partial_t u(T))\|_{L^q L^r ([a-T,b-T]\times \Rm^3)}\\
&=& \left\|S(t)(u_0^{(T)}, u_1^{(T)})\right\|_{L^q L^r ([\lambda(T)(a-T),\lambda(T)(b-T)]\times \Rm^3)},
\end{eqnarray*}
here
\[
 (u_0^{(T)}, u_1^{(T)}) = \left(\frac{1}{\lambda(T)^{3/2-s_p}} u\left(\frac{\cdot}{\lambda(T)}, T\right), \frac{1}{\lambda(T)^{5/2-s_p}}
\partial_t u\left(\frac{\cdot}{\lambda(T)}, T\right)\right).
\]
Given $\eps>0$, let $M, \delta$ be the constants as in lemma \ref{compact LqLr}.
It is clear that if $T$ is sufficiently large, we have either ($\lambda(T)$ is small)
\[
 \lambda(T)(b-T) - \lambda(T)(a-T) = (b-a)\lambda(T) < \delta;
\]
or ($\lambda(T)$ is large)
\[
 \lambda(T)(b-T) < -M.
\]
In either case, by lemma \ref{compact LqLr} we have $\|S(t-T)(u(T), \partial_t u(T))\|_{L^q L^r ([a,b]\times \Rm^3)} < \eps$. This completes the
proof of the first limit.\\
In order to obtain the second limit, we only need to choose $t_1 \in (t_0,+\infty)$, set $[a,b]=[t_0,t_1]$ and apply lemma \ref{weak limit}
below using the first limit and the following identity.
\[
 S(t-t_0) \left[S(t_0-T)\left(\begin{array}{c}u(T)\\ \partial_t u(T)\end{array}\right)\right] =
 S(t-T) \left(\begin{array}{c}u(T)\\ \partial_t u(T)\end{array}\right).
\]
\paragraph{Remark} We can obtain the similar result in the negative time direction using exactly the same argument. This implies the corresponding
Duhamel formula in the negative time direction.
\begin{itemize}
 \item Soliton-like Case or High-to-low Frequency Cascade Case
\[
 \lim_{T\rightarrow -\infty} \|S(t-T)(u(T), \partial_t u(T))\|_{L^q L^r ([a,b]\times \Rm^3)} =0.
\]
\[
 \hbox{weak}\,\lim_{T\rightarrow -\infty} S(t_0-T)\left(\begin{array}{c}u(T)\\ \partial_t u(T)\end{array}\right) = 0.
\]
 \item Self-similar Case (let $a,t_0 > 0$)
\[
 \lim_{T\rightarrow 0^+} \|S(t-T)(u(T), \partial_t u(T))\|_{L^q L^r ([a,b]\times \Rm^3)} =0.
\]
\[
 \hbox{weak}\,\lim_{T\rightarrow 0^+} S(t_0-T)\left(\begin{array}{c}u(T)\\ \partial_t u(T)\end{array}\right) = 0.
\]
\end{itemize}
\begin{lemma} \label{weak limit}
Suppose that $\{(u_{0,n}, u_{1,n})\}_{n\in {\mathbb Z}}$ is a bounded subset of $\dot{H}^s \times \dot{H}^{s-1}$ so that
\[
 \lim_{n \rightarrow \infty} \|S(t)(u_{0,n},u_{1,n})\|_{L^q L^r ([0,\mu]\times \Rm^3)} = 0.
\]
Here $(q,r)$ is an $s$-admissible pair and $\mu$ is a positive constant. Then we have the following weak limit in
$\dot{H}^s \times \dot{H}^{s-1}(\Rm^3)$,
\[
 (u_{0,n}, u_{1,n}) \rightharpoonup 0.
\]
\end{lemma}
\paragraph{Proof} Let us suppose the conclusion was false. This means that there exists a subsequence (WLOG, let us use the same notation as the original
sequence) so that it converges weakly to a nonzero limit $(\tilde{u}_0,\tilde{u}_1)$. We know the operator
$P: \dot{H}^s \times \dot{H}^{s-1} \rightarrow L^q L^r([0,\mu]\times \Rm^3)$ defined by
\[
 P(u_0,u_1) = S(t)(u_0,u_1)
\]
is bounded by the Strichartz estimate. This implies that we have the weak limit below in $L^q L^r([0,\mu]\times \Rm^3)$
\[
 P(u_{0,n},u_{1,n}) \rightharpoonup P(\tilde{u}_0,\tilde{u}_1).
\]
On the other hand, we know $P(u_{0,n},u_{1,n})$ converges to zero strongly. Thus $P(\tilde{u}_0,\tilde{u}_1) = 0$. This means
$(\tilde{u}_0,\tilde{u}_1) = 0$, which is a contradiction.
\begin{lemma} \label{Duhamel in Hs}
Assume $s \in [s_p,1]$. Let $u(x,t)$ be defined on $I=(T_-, \infty)$ and
almost periodic modulo scalings in $\dot{H}^s \times \dot{H}^{s-1}(\Rm^3)$, namely the set
\[
 K= \left\{\left(\frac{1}{\lambda(t)^{3/2-s_p}} u\left(\frac{x}{\lambda(t)},t\right),
 \frac{1}{\lambda(t)^{5/2-s_p}} \partial_t u\left(\frac{x}{\lambda(t)}, t\right)\right): t \in I\right\}
\]
is precompact in the space $\dot{H}^{s} \times \dot{H}^{s-1}(\Rm^3)$. In addition $\lambda(t) \leq 1$ for each $t \geq 1$. Then for
%any time $t_0 \in I$,
any closed interval $[a,b] \subset I$
and any $s$-admissible pair $(q,r)$ with $q <\infty$, we have
\[
 \lim_{T\rightarrow +\infty} \|S(t-T)(u(T), \partial_t u(T))\|_{L^q L^r ([a,b]\times \Rm^3)} =0.
\]
%\[
% \hbox{weak}\,\lim_{T\rightarrow +\infty} S(t_0-T)\left(\begin{array}{l}u(T),\\ \partial_t u(T)\end{array}\right) = 0.
%\]
\end{lemma}
\paragraph{Proof} One could use the similar method as used in lemma \ref{proof of Duhamel} by observing
\begin{eqnarray*}
&&\|S(t-T)(u(T), \partial_t u(T))\|_{L^q L^r ([a,b]\times \Rm^3)}\\ &=& \|S(t)(u(T), \partial_t u(T))\|_{L^q L^r ([a-T,b-T]\times \Rm^3)}\\
&=& (\lambda(T))^{s-s_p}\left\|S(t)(u_0^{(T)},u_1^{(T)})\right\|_{L^q L^r ([\lambda(T)(a-T),\lambda(T)(b-T)]\times \Rm^3)}.
\end{eqnarray*}
Here
\[
 (u_0^{(T)},u_1^{(T)})  = \left(\frac{1}{\lambda(T)^{3/2-s_p}} u\left(\frac{\cdot}{\lambda(T)}, T\right), \frac{1}{\lambda(T)^{5/2-s_p}}
\partial_t u\left(\frac{\cdot}{\lambda(T)}, T\right)\right).
\]
\subsection{Perturbation Theory}
\paragraph{Proof of theorem \ref{perturbation theory in Ysp}}
Let us first prove the perturbation theory when $M$ is sufficiently small. Let $I_1$ be the maximal lifespan of the solution $u(x,t)$ to the
equation (\ref{eqn}) with the given initial data $(u_0,u_1)$ and assume $[0,T] \subseteq I\cap I_1$. By the Strichartz estimate, we have
\begin{eqnarray*}
\|\tilde{u}-u\|_{Y_{s_p}([0,T])} &\leq& \|S(t)(u_0-\tilde{u}(0), u_1-\tilde{u}(0))\|_{Y_{s_p}([0,T])} + C_p
\|e + F (\tilde{u})- F(u)\|_{Z_{s_p}([0,T])}\\
&\leq& \eps + C_p \|e\|_{Z_{s_p}([0,T])} + C_p \|F(\tilde{u}) - F(u)\|_{Z_{s_p}([0,T])}\\
&\leq& \eps + C_p \eps + C_p \|\tilde{u}-u\|_{Y_{s_p}([0,T])} (\|\tilde{u}\|_{Y_{s_p}([0,T])}^{p-1} + \|\tilde{u}-u\|_{Y_{s_p}([0,T])}^{p-1})\\
&\leq& C_p \eps + C_p \|\tilde{u}-u\|_{Y_{s_p}([0,T])} (M^{p-1} + \|\tilde{u}-u\|_{Y_{s_p}([0,T])}^{p-1}).
\end{eqnarray*}
By a continuity argument in $T$, there exist $M_0 = M_0 (p), \eps_0 = \eps_0 (p) >0$, such that if $M \leq M_0$ and $\eps < \eps_0$, we have
\[
 \|\tilde{u}-u\|_{Y_{s_p}([0,T])} \leq C_p \eps.
\]
Observing that this estimate does not depend on the time $T$, we are actually able to conclude $I \subseteq I_1$
by the standard blow-up criterion and obtain
\[
 \|\tilde{u}-u\|_{Y_{s_p}(I)} \leq C_p \eps.
\]
In addition, by the Strichartz estimate
\begin{eqnarray*}
 \lefteqn{\sup_{t \in I} \left\|\left(\begin{array}{c} u(t)\\ \partial_t u(t)\end{array}\right)
 - \left(\begin{array}{c} \tilde{u}(t)\\ \partial_t \tilde{u}(t)\end{array}\right)
 - S(t)\left(\begin{array}{c} u_0 - \tilde{u}(0)\\ u_1 -\partial_t \tilde{u}(0)\end{array}\right)
 \right\|_{\dot{H}^{s_p} \times \dot{H}^{s_p-1}}}\\
 &\leq& C_p\| F(u) - F(\tilde{u}) - e\|_{Z_{s_p}(I)}\\
 &\leq& C_p\left(\|e\|_{Z_{s_p}(I)} + \|F(u) - F(\tilde{u})\|_{Z_{s_p}(I)}\right)\\
 &\leq& C_p\left[\eps + \|u -\tilde{u}\|_{Y_{s_p}(I)}\left(\|\tilde{u}\|_{Y_{s_p}(I)}^{p-1} + \|u-\tilde{u}\|_{Y_{s_p}}^{p-1}\right)\right]\\
 &\leq& C_p \eps.
\end{eqnarray*}
This finishes the proof as $M$ is sufficiently small. To deal with the general case,
we can separate the time interval $I$ into finite number of subintervals $\{I_j\}$,
so that $\|\tilde{u}\|_{Y_{s_p}(I_j)} < M_0$, and then iterate our argument above.
\paragraph{Proof of theorem \ref{perturbation theory in Hs}}
Let us first prove the perturbation theory when $M$ and $T$ are sufficiently small. Let $I_1$ be the maximal lifespan of the solution $u(x,t)$ to the
equation (\ref{eqn}) with the given initial data $(u_0,u_1)$ and assume $[0,T_1] \subseteq [0,T]\cap I_1$. By the Strichartz estimate, we have
\begin{eqnarray*}
\|\tilde{u}-u\|_{Y_{s}([0,T_1])} &\leq& \|S(t)(u_0-\tilde{u}_0, u_1-\tilde{u}_1)\|_{Y_{s}([0,T_1])} + C_{s,p}
\|F (\tilde{u})- F(u)\|_{Z_{s}([0,T_1])}\\
&\leq& C_{s,p} \|(u_0 -\tilde{u}_0, u_1 -\tilde{u}_1)\|_{\dot{H}^s \times \dot{H}^{s-1}} + C_{s,p} \|F(\tilde{u}) - F(u)\|_{Z_{s}([0,T_1])}\\
&\leq& C_{s,p} \|(u_0 -\tilde{u}_0, u_1 -\tilde{u}_1)\|_{\dot{H}^s \times \dot{H}^{s-1}}\\
&& + C_{s,p} T_1^{(p-1)(s-s_p)} \|F(\tilde{u})-F(u)\|_{L^{\frac{2}{s+1 -(2p-2)(s-s_p)}}L^{\frac{2}{2-s}}([0,T_1]\times \Rm^3)}\\
&\leq& C_{s,p} \|(u_0 -\tilde{u}_0, u_1 -\tilde{u}_1)\|_{\dot{H}^s \times \dot{H}^{s-1}}\\
&& + C_{s,p} T_1^{(p-1)(s-s_p)} \|\tilde{u}-u\|_{Y_s([0,T_1])} \left(\|\tilde{u}-u\|_{Y_s([0,T_1])}^{p-1} + \|\tilde{u}\|_{Y_s([0,T_1])}^{p-1}\right)\\
&\leq& C_{s,p} \|(u_0 -\tilde{u}_0, u_1 -\tilde{u}_1)\|_{\dot{H}^s \times \dot{H}^{s-1}}\\
&& + C_{s,p} T_1^{(p-1)(s-s_p)} \|\tilde{u}-u\|_{Y_s([0,T_1])} \left(\|\tilde{u}-u\|_{Y_s([0,T_1])}^{p-1} + M^{p-1}\right).
\end{eqnarray*}
By a continuity argument in $T_1$, there exist $M_0 = M_0 (s,p), \eps_0 = \eps_0 (s,p) >0$, such that if $M \leq M_0$, $T \leq 1$ and
$\|(u_0 -\tilde{u}_0, u_1 -\tilde{u}_1)\|_{\dot{H}^s \times \dot{H}^{s-1}} \leq \eps_0$, we have
\[
 \|\tilde{u}-u\|_{Y_{s}([0,T_1])} \leq C_{s,p} \|(u_0 -\tilde{u}_0, u_1 -\tilde{u}_1)\|_{\dot{H}^s \times \dot{H}^{s-1}}.
\]
Observing that this estimate does not depend on the time $T_1$ as long as $T_1 \leq T \leq 1$,
we are actually able to conclude $[0,T] \subseteq I_1$ by theorem \ref{local existence in Hs} and obtain
\[
 \|\tilde{u}-u\|_{Y_{s}([0,T])} \leq C_{s,p} \|(u_0 -\tilde{u}_0, u_1 -\tilde{u}_1)\|_{\dot{H}^s \times \dot{H}^{s-1}}.
\]
In addition, by the Strichartz estimate
\begin{eqnarray*}
 \lefteqn{\sup_{t \in I} \left\|\left(\begin{array}{c} u(t)\\ \partial_t u(t)\end{array}\right)
 - \left(\begin{array}{c} \tilde{u}(t)\\ \partial_t \tilde{u}(t)\end{array}\right)
 \right\|_{\dot{H}^{s} \times \dot{H}^{s-1}}}\\
 &\leq& \left\|S(t)\left(\begin{array}{c} u_0 - \tilde{u}_0\\ u_1 -\tilde{u}_1\end{array}\right)\right\|_{\dot{H}^s \times \dot{H}^{s-1}} +
 C_{s,p}\| F(u) - F(\tilde{u})\|_{Z_{s}([0,T])}\\
 &\leq& C_{s,p} \|(u_0 -\tilde{u}_0, u_1 -\tilde{u}_1)\|_{\dot{H}^s \times \dot{H}^{s-1}}\\
 &&+ C_{s,p} T^{(p-1)(s-s_p)} \|\tilde{u}-u\|_{Y_s([0,T])} (\|\tilde{u}-u\|_{Y_s([0,T])}^{p-1} + \|\tilde{u}\|_{Y_s([0,T])}^{p-1})\\
 &\leq& C_{s,p} \|(u_0 -\tilde{u}_0, u_1 -\tilde{u}_1)\|_{\dot{H}^s \times \dot{H}^{s-1}}.
\end{eqnarray*}
This finishes the proof as $M$ and $T$ are sufficiently small. To deal with the general case,
we can separate the time interval $[0,T]$ into finite number of subintervals $\{I_j\}$,
so that $\|\tilde{u}\|_{Y_{s}(I_j)} \leq M_0$ and $|I_j| \leq 1$, then iterate our argument above.
\subsection{Technical Lemmas}
\begin{lemma} \label{app a}
 Suppose that $(u_{0,\eps}(x), u_{1, \eps}(x))$ are radial, smooth pairs defined in $\Rm^3$ and converge to
$(u_0(x), u_1(x))$ strongly in $\dot{H}^{s_p} \times \dot{H}^{s_p -1}(\Rm^3)$. In addition, we have
\[
 \int_{r_0< |x|<4r_0} (|\nabla u_{0,\eps} (x,t_0)|^2+ | u_{1,\eps} (x,t_0)|^2) dx \leq C
\]
for each $\eps < \eps_0$. Then $(u_0(x), u_1(x))$ is in the space $\dot{H}^1 \times L^2 (r<|x|<4r)$ and satisfies
\[
 \int_{r_0< |x|<4r_0} (|\nabla u_{0} (x)|^2+ |u_{1} (x)|^2) dx \leq C.
\]
\end{lemma}
\paragraph{Proof} By the uniform bound of the integral, we can extract a sequence $\eps_i \rightarrow 0$
so that $\partial_r u_{0,\eps_i}(r)$ converges to $\tilde{u}'_0(r)$ weakly in $L^2 (r_0,4r_0)$, and
$u_{1,\eps_i}$ converges to $\tilde{u}_1$ weakly in $L^2 (r_0 < |x| <4r_0)$. Define
\[
 \tilde{u}_0 (r) = u_0(r_0) + \int_{r_0}^{r} \tilde{u}'_0(\tau) d \tau.
\]
We have
\[
 \int_{r_0< |x|<4r_0} (|\nabla \tilde{u}_{0} (x)|^2+ |\tilde{u}_{1} (x)|^2) dx \leq C.
\]
By the strong and weak convergence, we have immediately $u_1 = \tilde{u}_1$ in the region $r_0 < |x| < 4 r_0$.
In order to conclude, we only need to show $u_0(r) = \tilde{u}_0 (r)$. Observing $\int_{r_0}^{r_1} f(\tau) d\tau$ is a bounded
linear functional in $L^2 (r_0,4r_0)$ for each $r_1 \in (r_0,4r_0)$, we have
\begin{eqnarray*}
 \tilde{u}_0 (r_1) &=& u_0(r_0) + \int_{r_0}^{r_1} \tilde{u}'_0(\tau) d \tau\\
                 &=& \lim_{i \rightarrow \infty} u_{0,\eps_i}(r_0) + \lim_{i \rightarrow \infty}
  \int_{r_0}^{r_1} \partial_r u_{0,\eps_i}(\tau) d \tau\\
&=& \lim_{i \rightarrow \infty} \left(u_{0,\eps_i}(r_0) + \int_{r_0}^{r_1} \partial_r u_{0,\eps_i}(\tau) d \tau\right)\\
&=& \lim_{i \rightarrow \infty} u_{0,\eps_i} (r_1)\\
&=& u_0 (r_1).
\end{eqnarray*}
This completes the proof.
\begin{lemma}
There exists a constant $\kappa = \kappa(p) \in (0,1)$ that depends only on $p$, so that for each $s\in [s_p, 1)$,
there exists an $s$-admissible pair $(q,r)$, with $q \neq \infty$ and
\[
 \frac{s + 1 - (2p-2)(s-s_p)}{2p} = \kappa \cdot 0 + (1-\kappa) \frac{1}{q};\;\;\; \frac{2-s}{2p} = \kappa \frac{3-2s}{6}
 + (1-\kappa) \frac{1}{r}.
\]
\end{lemma}
\paragraph{Proof} We will choose $\kappa = 1 - \frac{3}{p} \in (0,0.4)$. Basic Computation shows
\[
 \frac{1}{q} = \frac{s + 1 - (2p-2)(s-s_p)}{2p (1-\kappa)} = \frac{s + 1 - (2p-2)(s-s_p)}{6} \in (0,1/3);
\]
\begin{eqnarray*}
 \frac{1}{r} &=& \frac{2-s}{2p(1-\kappa)} - \frac{\kappa}{1-\kappa} \times \frac{3-2s}{6}\\
 &=& \frac{2-s}{6} -\frac{\kappa}{1-\kappa} \times  \frac{3-2s}{6}\\
 &\in& \left(\frac{2-s}{6} -\frac{2}{3} \times \frac{3-2s}{6},\;\frac{2-s}{6}\right)\\
 &\subseteq& \left(\frac{s}{18}, \frac{2-s}{6}\right)\\
 &\subseteq& \left(1/36, 1/4\right)
\end{eqnarray*}
Thus we can solve two positive real number $q,r$ so that the two identities hold. In addition, we have
$q \in (3,\infty)$ and $r \in (4,36)$. Furthermore, by adding the identities together, we obtain
\[
 \frac{3 - (2p-2)(s -s_p)}{2p} = \kappa \frac{3-2s}{6} + (1 -\kappa) (\frac{1}{q} + \frac{1}{r})
\]
This implies
\[
 \frac{1}{q} + \frac{1}{r} < \frac{3 - (2p-2)(s -s_p)}{2p(1-\kappa)} = \frac{3 - (2p-2)(s -s_p)}{6} \leq 1/2.
\]
Using the same method, one can show $1/q + 3/r = 3/2 - s$. In summary, $(q,r)$ is an $s$-admissible pair.

\end{document}